\documentclass[conference]{IEEEtran}

\usepackage{amsmath,amssymb,epsfig,verbatim}
\usepackage[T1]{fontenc}
\usepackage{ctable}
\usepackage{amsthm}
\usepackage{textcomp}
\usepackage{caption}
\usepackage{subcaption}

\def\reals{\mathbb{R}}

\def\comp{\raise 1pt \hbox{$\scriptstyle\circ$}}

\def\argmin{\mathop{\rm argmin}\limits}
\def\argmax{\mathop{\rm argmax}\limits}
\def\minimizeQuotes{\mathop{\rm \mbox{``}min\mbox{''}}\limits}
\def\minimize{\mathop{\rm min}\limits}

\def\st{\mathop{\rm subject\ to}}

\def\gph{\mathop{\rm gph}}

\def\upto{{\raise 1pt \hbox{$\scriptstyle \,\nearrow\,$}}}
\def\downto{{\raise 1pt \hbox{$\scriptstyle \,\searrow\,$}}}

\def\tos{\rightrightarrows}

\newtheorem{theorem}{Theorem}

\newtheorem{definition}[theorem]{Definition}

\begin{document}
\newcommand{\bolda}{\mbox{$\mathbf{a}$}}
\newcommand{\boldb}{\mbox{$\mathbf{b}$}}
\newcommand{\boldx}{\mbox{$\mathbf{x}$}}
\newcommand{\boldy}{\mbox{$\mathbf{y}$}}
\newcommand{\boldf}{\mbox{$\mathbf{f}$}}
\newcommand{\boldz}{\mbox{$\mathbf{z}$}}
\newcommand{\boldF}{\mbox{$\mathbf{F}$}}
\newcommand{\boldG}{\mbox{$\mathbf{G}$}}
\newcommand{\boldg}{\mbox{$\mathbf{g}$}}
\newcommand{\boldh}{\mbox{$\mathbf{h}$}}
\newcommand{\boldH}{\mbox{$\mathbf{H}$}}
\newcommand{\boldzero}{\mbox{$\mathbf{0}$}}
\newcommand{\Rbb}{\mbox{$\mathbb R$}}

\title{A Review on Bilevel Optimization: From Classical to Evolutionary Approaches and Applications}

\author{Ankur Sinha, Pekka Malo, Kalyanmoy Deb
\thanks{Ankur Sinha is an Assistant Professor at the Department of Production and Quantitative Methods, Indian Institute of Management, Ahmedabad 380015, India (asinha@iimahd.ernet.in).}
\thanks{Pekka Malo is an Assistant Professor at the Department of Information and Service Economy, Aalto University School of Business, PO Box 21220, FI-00076 AALTO, Finland (pekka.malo@aalto.fi).}
\thanks{Kalyanmoy Deb is a Professor and Koenig Endowed Chair at the  Department of Electrical and Computer Engineering, Michigan State University, East Lansing, MI 48824, USA (kdeb@egr.msu.edu).}
}

\maketitle

\begin{abstract}
Bilevel optimization is defined as a mathematical program, where an optimization problem contains another optimization problem as a constraint. These problems have received significant attention from the mathematical programming community. Only limited work exists on bilevel problems using evolutionary computation techniques; however, recently there has been an increasing interest due to the proliferation of practical applications and the potential of evolutionary algorithms in tackling these problems. This paper provides a comprehensive review on bilevel optimization from the basic principles to solution strategies; both classical and evolutionary. A number of potential application problems are also discussed. To offer the readers insights on the prominent developments in the field of bilevel optimization, we have performed an automated text-analysis of an extended list of papers published on bilevel optimization to date. This paper should motivate evolutionary computation researchers to pay more attention to this practical yet challenging area.
\end{abstract}


\IEEEdisplaynontitleabstractindextext

\IEEEpeerreviewmaketitle

\section{Introduction}

Many large-scale optimization and decision-making processes faced by public and private organizations are {\em hierarchical\/} in the sense that the realized outcome of any solution or decision taken by the upper level authority (leader) to optimize their goals is affected by the response of lower level entities (follower), who will seek to optimize their own outcomes. {\color{black} Figure~\ref{fig:sketch} illustrates a general bilevel problem solving structure involving interlinked optimization and decision-making tasks at both levels. 
\begin{figure}[hbt]
\epsfig{file=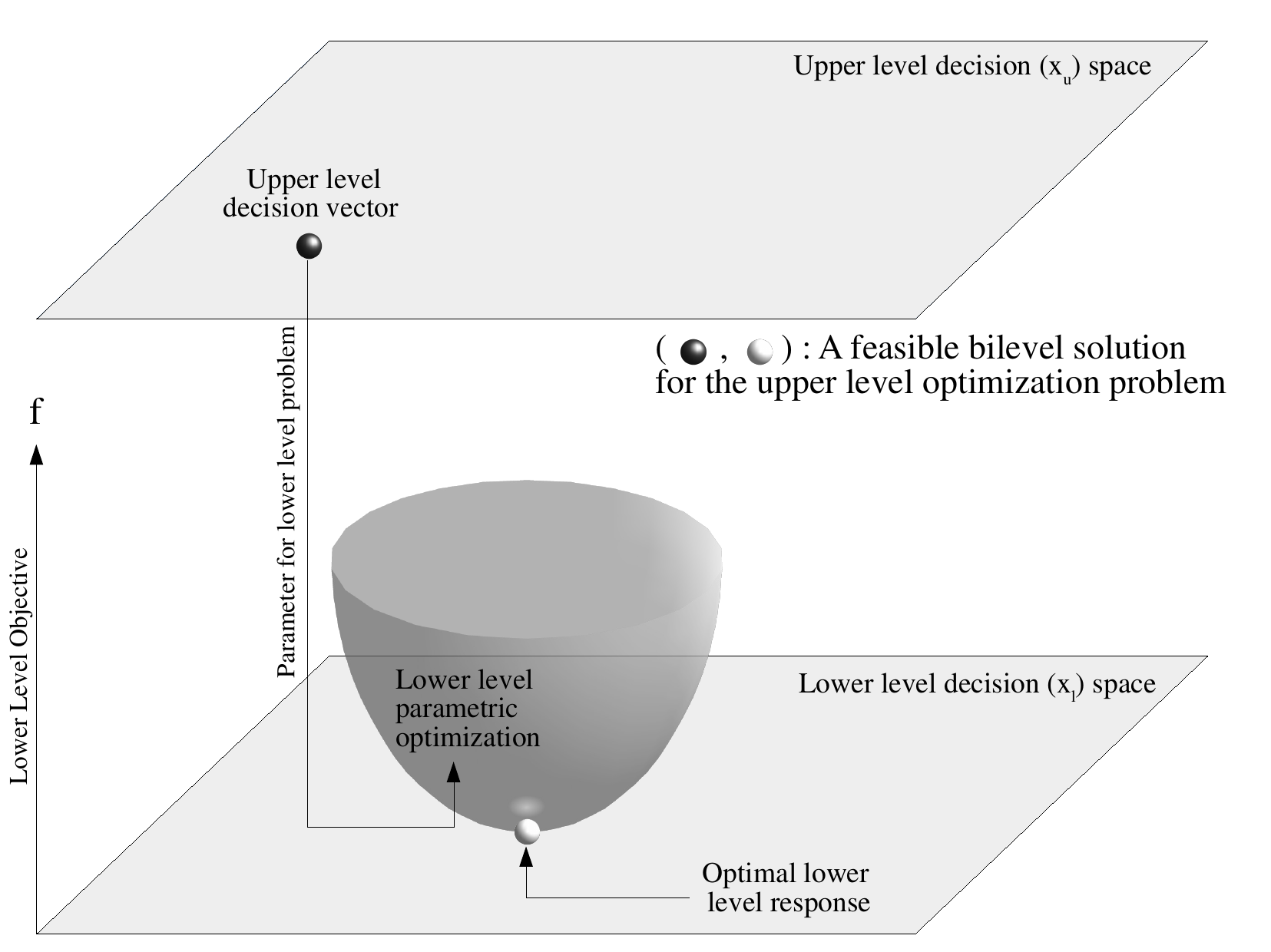,width=1.0\linewidth}
\caption{A general sketch of a bilevel problem.}
\vspace{-3mm}
\label{fig:sketch}
\end{figure}
The figure shows that for any given upper level decision vector, there is a corresponding (parametric) lower level optimization problem to be solved that provides the rational (optimal) response of the follower for the leader's decision. The leader's decision vector is represented by $x_u$ and the follower's decision vector is represented by $x_l$. An ($x_u$, $x_{l}^{\ast}$) pair where $x_{l}^{\ast}$ is an optimal response to $x_u$ represents a feasible solution to the upper level optimization problem provided that it also satisfies the constraints in the problem.} Each level has its own objectives and constraints. One aspect of bilevel problems is that it is not symmetric in terms of two levels. The upper level decision maker usually has complete knowledge of the lower level problem, while the lower level decision maker only observes the decisions of the leader and then optimizes its own strategies. Interestingly, an incomplete knowledge about the follower's optimization problem to the leader may lead to bilevel optimization problems involving uncertainties.

It is not uncommon that the objectives of generally profit-seeking private agents can well be in conflict with those of the controlling authority. {\color{black} What makes such bilevel problem-solving tasks highly relevant is that they are typically characterized by very large spillover effects to the economy as well as the surrounding environment. Given the far-reaching future impacts of the decisions, it is not surprising that the interest towards bilevel programming has grown strong especially among researchers and practitioners dealing with large-scale public sector decision-making problems. For instance, farmers often tend to overuse fertilizers to increase the productivity, which leads to negative externalities such as pollution. In \cite{whittaker2016spatial} authors use a bilevel model to design policy measures to control the overuse of fertilizers and its negative impact on the environment.}
Apart from problems concerned with environmentally sensitive decisions (such as allocation of mining permits or controlling the use of fertilizers), there has been widespread interest across a number of fields in operations research. A good example is homeland security, where bilevel as well as trilevel optimization frameworks have been utilized in problems ranging from interdiction of nuclear-weapons projects to defending critical infrastructure and solving border security problems~\cite{an13,wein09,brown05,brown09}. In addition to the public sector challenges, there is abundant research on bilevel decision-making problems in economics, logistics, as well as diverse areas of computer science.

The research on decision-making problems with hierarchical leader-follower structures (bilevel optimization) can be traced to two roots. The first root is in the domain of game theory, where Stackelberg~\cite{St52} used bilevel programming to build descriptive models of decision behavior and establish game-theoretic equilibria. The second root is in the domain of mathematical programming, where the problems appeared as bilevel optimization problems containing a nested inner optimization problem as a constraint of an outer optimization problem~\cite{BrMc73}. Since then a substantial body of mathematical literature on bilevel optimization has emerged. Given that the hierarchical optimization structure may introduce difficulties such as non-convexity and disconnectedness even for simpler instances of bilevel optimization, the problems have turned out to be surprisingly difficult to handle mathematically. {\color{black}Bilevel programming is known to be strongly NP-hard~\cite{HaJaSa92}, and it has been proven that merely evaluating a solution for optimality is also a NP-hard task~\cite{ViSaJu94}. Even in the simplest case of linear bilevel programs, where the lower level problem has a unique optimal solution for all the parameters, it is not likely to find a polynomial algorithm that is capable of solving the linear bilevel program to global optimality. The proof for the non-existence of a polynomial time algorithm for linear bilevel problems can be found in \cite{deng1998complexity}.}

Due to lack of well-established solution procedures, a complex practical problem is usually modified into simpler single level optimization task, which is solved to arrive at a satisficing instead of an optimal solution. For the complex bilevel problems, classical methods often fail due to real world difficulties such as non-linearity and discreteness. Under such circumstances, evolutionary methods can be useful tools to offset some of these difficulties.
Recent initiatives on bilevel optimization using evolutionary algorithms suggest that a coordinated effort on bilevel optimization by the evolutionary community could help make significant progress on this challenging class of optimization problems (e.g., \cite{my-ieeetec16,my-bleaq-arxiv13,angelo13,my-ecj10}). 

\begin{figure*}[hbt]
\begin{center}
\epsfig{file=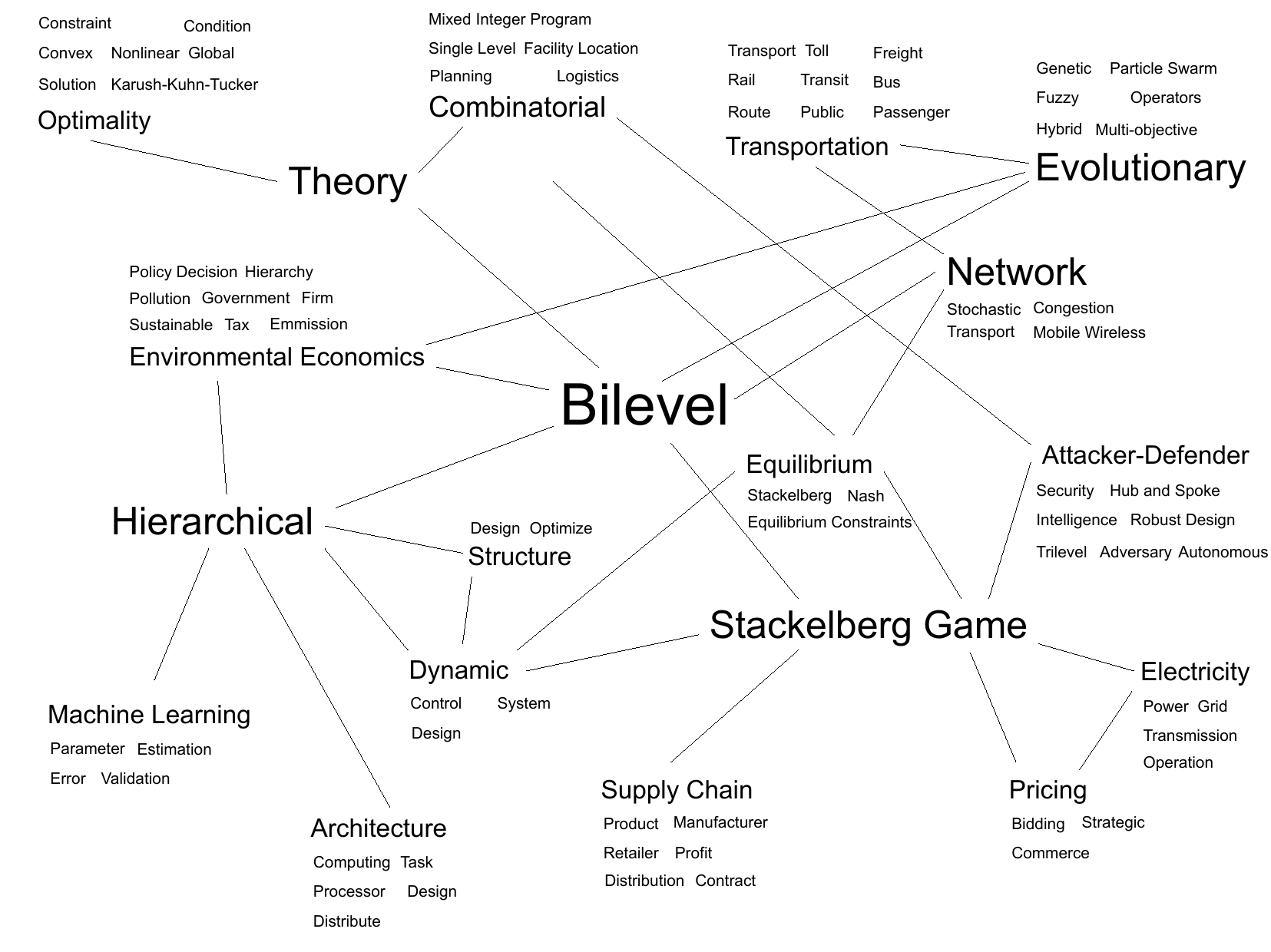,width=0.95\linewidth}
\end{center}
\caption{Bilevel network-map showing connections between various applications and theory since 1950s. Each connecting link represents either a topic connected with a subtopic, or an overlap between two subtopics.}
\label{fig:network}
\end{figure*}

{\color{black}Figure~\ref{fig:network} provides a network map of different themes on bilevel optimization that have been studied since 1950s. The network map shows different theoretical and application topics that have evolved under bilevel optimization. Each link in the map connects a subtopic with a higher level topic that are differentiated by font sizes. Subtopics connected with a link denote an overlap.} To provide a more comprehensive overview on the past as well as recent developments in the field of bilevel optimization, we have organized this review paper along the three lines; theory, applications and text-analysis of the entire bilevel literature body. First, to formalize the notion of bilevel programming, we begin by introducing a few central definitions and discuss the differences between optimistic and pessimistic formulations of bilevel problems. Once the common terminology has been established, we offer an overview on the algorithms that have been proposed for bilevel optimization. After a brief coverage of the commonly used classical approaches (e.g., descent methods, penalty function methods, and trust region methods), we move on to discuss the developments in the field of evolutionary computation, discrete bilevel optimization and multiobjective bilevel optimization. The method sections are followed by a review on the central application areas. Finally, we study the research topics, and the evolution of interest over time. The entire bilevel literature is divided into topics and a time series analysis across each topic is performed. The text-analysis performed in this paper is based on a recently developed non-parametric topic model~\cite{blei2012probabilistic,teh2006hierarchical} for analyzing unstructured information. The technical details on the automated text-analysis approach are provided in an appendix. The paper concludes with a brief discussion on the directions for future research.

\section{General Formulation and Definitions}
In this section, we provide a general formulation for bilevel optimization problem. These problems contain two levels of optimization tasks where one optimization task is nested within the other. The outer optimization problem is commonly referred as the leader's (upper level) optimization problem and the inner optimization problem is known as the follower's (or lower level) optimization problem. The two levels have their own objectives and constraints. Correspondingly, there are also two classes of decision vectors, namely, leader's (upper level) decision vectors and follower's (lower level) decision vectors. The lower level optimization is a parametric optimization problem that is solved with respect to the lower level decision vectors while the upper level decision vectors act as parameters. The lower level optimization problem is a constraint to the upper level optimization problem, such that, only those members are considered feasible that are lower level optimal and also satisfy the upper level constraints. A summary of the terminologies and notations used in the context of bilevel optimization is given in Table~\ref{tab:notations}.


\begin{table*}[!ht]
{\footnotesize
\caption{Summary of central notations}\label{tab:notations}
\begin{center}
\begin{tabular}{p{.30\textwidth} >{$\displaystyle}p{.15\textwidth}<{$}  p{.45\textwidth}}
\hline
Category & \text{Notation(s)} & Description\\
\toprule
Decision vectors & x_u\in X_U & Leader's (upper level) decision variable and decision space.\\
	 & x_l \in X_L & Follower's (lower level) decision variable and decision space.\\

\midrule
Objectives & F & Leader's (upper level) objective functions.\\
& f & Follower's (lower level) objective functions.\\
\midrule
Constraints & G_k, \ k=1,\dots,K  & Leader's (upper level) constraint functions.\\
 & g_j, \ j=1,\dots,J & Follower's (lower level) constraint functions.\\

\midrule
Lower level feasible region & \Omega: X_U \tos X_L & $\Omega(x_u)=\{x_l: g_j(x_u,x_l) \leq 0 \; \forall \; j\}$, represents the lower level feasible region for any given upper level decision vector\\

\midrule
Constraint region (Relaxed feasible set) & \Phi = \gph \Omega & $\Phi=\{(x_u,x_l) : G_k(x_u,x_l) \leq 0  \; \forall \; k, \; g_j(x_u,x_l) \leq 0  \; \forall \; j \}$, represents the region satisfying both upper and lower level constraints\\

\midrule
Lower level reaction set & \Psi: X_U \tos X_L & $\Psi(x_u)=\{x_l: x_l \in \argmin_{x_l \in X_L} \{f(x_u,x_l): x_l \in \Omega(x_u) \} \}$, represents the lower level optimal solution(s) for an upper level decision vector\\

\midrule
{\color{black} Inducible region (Feasible set)} & I=\gph \Psi & $I=\{(x_u,x_l): (x_u,x_l) \in \Phi, x_l \in \Psi(x_u) \}$, {\color{black} represents the set of upper level decision vectors and corresponding lower level optimal solution(s) belonging to feasible constraint region}\\

\midrule
Choice function & \psi: X_U \to X_L & $\psi(x_u)$ represents the solution chosen by the follower for any upper level decision vector. It becomes important in case of multiple lower level optimal solutions.\\

\midrule
Optimal value function & \varphi: X_U \to R & $\varphi(x_u)=\minimize_{x_l \in X_L} \{f(x_u,x_l): x_l \in \Omega(x_u) \} $ represents the minimum lower level function value corresponding to a given upper level decision vector.\\

\bottomrule		
\end{tabular}
\end{center}
\label{default}
}
\end{table*}%

\begin{definition}\label{def:bilevel1}
For the upper-level objective function $F:\reals^n\times\reals^m \to\reals$ and lower-level objective function $f:\reals^n\times\reals^m \to\reals$, the bilevel problem is given by 
\begin{align*}
\minimizeQuotes_{x_u\in X_U,x_l\in X_L} & F(x_u,x_l) \\
\st &\\  & \hspace{-12mm}x_l\in \argmin_{x_l \in X_L} 
	\lbrace
		f(x_u,x_l) : g_j(x_u,x_l)\leq 0, j=1,\dots,J
	\rbrace\\
 & \hspace{-12mm}G_k(x_u,x_l)\leq 0, k=1,\dots,K
\end{align*}
where $G_k:\reals^n\times\reals^m \to \reals$, $k=1,\dots,K$ denote the upper level constraints, and $g_j:\reals^n\times\reals^m \to \reals$ represent the lower level constraints, respectively. Equality constraints may also exist that have been avoided for brevity. The sets $X_U \subset \reals^n$ and $X_L \subset \reals^m$ in the definition may denote additional restrictions like integrality. It is common to assume these to be sets of reals, unless mentioned otherwise.
\end{definition}

An equivalent formulation of the above problem can be stated in terms of {\color{black}set-valued  mapping (multi-valued function)} as follows:

\begin{definition}\label{def:bilevel2}
Let $\Psi:\reals^n\tos\reals^m$ be a set-valued mapping, 
\begin{align*}
\Psi(x_u)=&\argmin_{x_l \in X_L}\{f(x_u,x_l) : g_j(x_u,x_l)\leq 0, j=1,\dots,J\},
\end{align*}
which represents the constraint defined by the lower-level optimization problem, i.e. $\Psi(x_u)\subset X_L$ for every $x_u\in X_U$. Then the bilevel optimization problem can be expressed as a constrained optimization problem as follows: 
\begin{align*}
\minimizeQuotes_{x_u\in X_U, x_l\in X_L}\quad & F(x_u,x_l) \\
\st\quad  & \\
 & \hspace{-12mm} x_l \in \Psi(x_u) \\
 & \hspace{-12mm} G_k(x_u,x_l)\leq 0, k=1,\dots,K
\end{align*}
where $\Psi$ can be interpreted as a parameterized range-constraint for the lower-level decision vector $x_l$. 
\end{definition}


In the above two definitions, quotes have been used while specifying the upper level minimization problem because of an ambiguity that arises in case of multiple lower level optimal solutions for any given upper level decision vector. In the presence of multiple lower level optimal solutions there is lack of clarity at the upper level as to which optimal solution from the lower level should be utilized. This ambiguity can be sorted by defining different positions that may be assumed by the leader. The two common positions that have been widely studied are optimistic (weak) position and pessimistic (strong) position, which we discuss next.

\subsection{Optimistic Position}
In an optimistic position, in the presence of multiple lower level optimal solutions, the leader expects the follower to choose that solution from the optimal set $\Psi^{\rm o}(x_u)$, which leads to the best objective function value at the upper level. The choice function of the follower in this case may be defined as follows:
$$
\Psi^{\rm o}(x_u) = \argmin_{x_l \in X_L}\{F(x_u,x_l): x_l \in \Psi(x_u)\}
$$
This formulation assumes some extent of cooperation between the two players. The bilevel optimization problem under an optimistic position  has been defined below:
\begin{align*}
\minimize_{x_u\in X_U, x_l\in X_L}\quad & F(x_u,x_l) \\
\st\quad  & \\
 & \hspace{-12mm} x_l = \Psi^{\rm o}(x_u) \\
 & \hspace{-12mm} G_k(x_u,x_l)\leq 0, k=1,\dots,K
\end{align*}

Optimistic position is more tractable as compared to the pessimistic position; therefore, most of the studies handle optimistic version of the bilevel optimization problem. The optimistic formulation is guaranteed to have an optimal solutions under reasonable assumptions of regularity and compactness that are stated in the theorem below:

\setcounter{theorem}{0}
\begin{theorem}
If the functions $F, f, G_k$ and $g_i$ are sufficiently smooth, the constraint region $\Phi$ of the bilevel optimization problem is non-empty and compact, and the Mangasarian-Fromowitz constraint qualification holds at all points, then the problem is guaranteed to have an optimistic bilevel optimum provided there exists a feasible solution.
\end{theorem}
See \cite{HaPa88,Ou93,LiMo95,LiMo02,dempe02,dempe2007new} for further discussion on existence of optimistic bilevel optimum and additional results on optimality conditions.

\subsection{Pessimistic Position}
In a pessimistic position, in the presence of multiple lower level optimal solutions, the leader optimizes for the worst case, i.e., she assumes that the follower may choose that solution from the optimal set which leads to the worst objective function value at the upper level. Such a worst case choice function of the follower may be defined as:
$$
\Psi^{\rm p}(x_u) = \argmax_{x_l \in X_L}\{F(x_u,x_l): x_l \in \Psi(x_u)\}
$$
This formulation does not assume any form of cooperation. The bilevel optimization problem under a pessimistic position  has been defined below:
\begin{align*}
\minimize_{x_u\in X_U, x_l\in X_L}\quad & F(x_u,x_l) \\
\st\quad  & \\
 & \hspace{-12mm} x_l = \Psi^{\rm p}(x_u) \\
 & \hspace{-12mm} G_k(x_u,x_l)\leq 0, k=1,\dots,K
\end{align*}
{\color{black}Pessimistic position is relatively less tractable when compared to optimistic position. In case of an optimisitc formulation with a convex lower level problem, it is possible to reduce the bilevel problem to single level using the variational inequality corresponding to the lower level problem. However, such a straightforward single level reduction is not possible in case of a pessimistic bilevel program. This poses significant challenges in designing methodologies that can handle pessimistic bilevel problems. For every lower level optimization problem solved one has to keep track of that lower level optimal solution that is worst for the upper level. This essentially makes a pessimistic bilevel optimization a three level task.} The pessimistic formulation is guaranteed to have an optimal solutions under stronger assumptions, as compared to the optimistic formulation, that are given below: 
\begin{theorem}
If the functions $F, f, G_k$ and $g_i$ are sufficiently smooth, the constraint region $\Phi$ of the bilevel optimization problem is non-empty and compact, and the set-value mapping, $\Psi^{\rm p}$, is lower semi-continuous for all upper level decision vectors, then the problem is guaranteed to have a pessimistic bilevel optimum. 
\end{theorem}
For discussion on existence of pessimistic bilevel optimum and additional results on optimality conditions, the authors may refer to \cite{LuMiPi87,LoMo96,dempe02,dempe2014necessary,wiesemann2013pessimistic}.

{\color{black}
\subsection{Example}
Below we provide a simple example of a bilevel optimization problem \cite{my-ifac12} that arises in case of two firms in a  Stackelberg competition. The leader has complete knowledge about the follower's inverse-demand function and the cost function, and desires to maximize it's own profits by taking into account the actions of the follower firm. The two firms compete solely by choosing their production levels that maximize their profits ($\Pi_l$ and $\Pi_f$), and the follower acts only after observing the actions of the leader. Formally, this model can be presented as follows:
\begin{eqnarray} \label{eq:simple}
\max_{q_l,q_f} \quad&& \Pi_l = P(q_l,q_f) q_l - C_l(q_l) \\
\mbox{s.t.} \quad&& q_f \in \argmax_{q_f} \lbrace {\Pi_f = P(q_l,q_f) q_f - C_f(q_f)} \rbrace , \\
&& q_l, q_f \geq 0 ,
\end{eqnarray}
where $P(q_l,q_f)$ is the unit price of the goods sold, which depends on the todal production. The assumption is that at the optimum, all demand is satisfied. $C_l(\cdot)$ is the cost of production of the leader and $C_f(\cdot)$ is the cost of production of the follower. The variables in this model are the production levels of each firm $q_l$ and $q_f$. The leader sets its production level first, and then the follower chooses its production level based on the leader's decision. This simple model assumes homogeneity of the products manufactured by the firms.

By assuming that the firms produce and sell homogeneous goods, we may assume a single linear price function for both firms as an inverse demand function of the form
\begin{equation} \label{eq:simpleprice}
P(q_l,q_f) = \alpha - \beta (q_l,q_f) ,
\end{equation}

where $\alpha, \beta > 0$ are constants. Additionally, since costs often tend to increase with the amount of production, we assume convex quadratic cost functions for both firms to be of the form
\begin{align} \label{eq:simplecosts}
C(q_l) & = \delta_l q_l^2 + \gamma_l q_l + c_l, \\
C(q_f) & = \delta_f q_f^2 + \gamma_f q_f + c_f,
\end{align}

where $c_i$ denotes the fixed costs of the respective firm, and $\delta_i$ and $\gamma_i$ are positive constants. It is possible to solve the above model analytically, using the first order conditions of the lower level problem to reduce it to single level, and then using the first order conditions of the reduced problem. 

The optimal level of production of the leader $(q_l^*)$ and the follower $(q_f^*)$ in terms of the constants of the model is given as follows:
\begin{align}
q_l^* &= \frac{2(\beta + \delta_f)(\alpha - \gamma_l) - \beta(\alpha - \gamma_f)}{4(\beta + \delta_f)(\beta + \delta_l) - 2\beta^2}. \label{eq:optlead}\\
q_f^* &= \frac{\alpha - \gamma_f}{2(\beta + \delta_f)} - \frac{\beta(\alpha - \gamma_l)-{\displaystyle \frac{\beta^2(\alpha - \gamma_f)}{2(\beta + \delta_f)}}}{4(\beta + \delta_f)(\beta + \delta_l) - 2\beta^2}. \label{eq:optfol}
\end{align}
}

\section{Classical Approaches}
In this section, we provide a brief overview of the classical algorithms that have been proposed for bilevel optimization. Given the difficult nature of bilevel problems, it is not surprising that much of the classical literature considers bilevel problems that are mathematically well-behaved; i.e., contains functions that are linear, quadratic or convex. Strong assumptions like continuous differentiability and lower semi-continuity are quite common. A significant amount of attention has been given to linear bilevel optimization problems with continuous \cite{WeHs91,Be93} and combinatorial \cite{ViSaJu96} variables. For more complex bilevel problems, the readers may refer \cite{Co99,AnFr92}.

\subsection{Single-level Reduction}\label{sec:singleLevelReduction1}
When the lower level problem is convex and sufficiently regular, it is possible to replace the lower level optimization problem with its Karush-Kuhn-Tucker (KKT) conditions. The KKT conditions appear as Lagrangian and complementarity constraints, and reduce the overall bilevel optimization problem to a single-level constrained optimization problem. For example, the problem in Definition \ref{def:bilevel1} can be reduced to the following form, when the convexity and regularity conditions at the lower level are met:

\begin{align*}
\minimize_{x_u\in X_U, x_l\in X_L, \lambda}\quad & F(x_u,x_l) \\
\st\quad  & \\
 & \hspace{-12mm} G_k(x_u,x_l)\leq 0, k=1,\dots,K,\\
 & \hspace{-12mm} \nabla_{x_l} L(x_u,x_l,\lambda) = 0,\\
 & \hspace{-12mm} g_j(x_u,x_l)\leq 0, j=1,\dots,J,\\
 & \hspace{-12mm} \lambda_j g_j(x_u,x_l) = 0, j=1,\dots,J,\\
 & \hspace{-12mm} \lambda_j \geq 0, j=1,\dots,J,\\
\mbox{where}\quad\quad\quad  & \\
 & \hspace{-12mm} L(x_u,x_l,\lambda) = f(x_u,x_l) + \sum_{j=1}^{J} \lambda_j g_j(x_u,x_l).
\end{align*}
The above formulation, though a single level optimization task, is not necessarily simple to handle. The Lagrangian constraints can lead to non-convexities even when suitable convexity assumptions are made on all the objectives and constraints in the bilevel formulation. The complementarity condition, inherently being combinatorial, renders the single-level optimization problem as a mixed integer program.

Interestingly, for linear bilevel optimization problems, the Lagrangian constraint is also linear. Therefore, the single-level optimization problem is a mixed integer linear program. Approaches based on vertex enumeration \cite{BiKa84,ChFl92,TuMiVa93}, as well as branch-and-bound (B\&B) \cite{BaFa82,FoMc81} have been used to solve these problems. It is noteworthy that B\&B methods constitute an exponentially slow algorithm with the number of integer variables. But, B\&B approaches have been successfully applied to single-level reductions of linear-quadratic \cite{BaMo90} and quadratic-quadratic \cite{EdBa91,AlHoPa92} bilevel problems. An extended KKT approach has also been considered \cite{shi2005extended} for handling linear bilevel problems. 

\subsection{Descent Methods}
In addition to KKT based approaches, a number of descent methods have been proposed for solving bilevel optimization problems. A descent direction in bilevel optimization leads to decrease in upper level function value while keeping the new point feasible. Given that a point is considered feasible only if it is lower level optimal, finding the descent direction can be quite challenging. To resolve the problem, researchers have investigated ways to approximate the gradient of the upper level objective \cite{KoLa90} as well as considered formulation of auxiliary programs \cite{SaGa94,ViSaJu94} to determine the direction of descent.

\subsection{Penalty Function Methods}
In penalty function methods the bilevel optimization problem is handled by solving a series of unconstrained optimization problems. The unconstrained problem is generated by adding a penalty term that measures the extent of violation of the constraints. The penalty term often requires a parameter and takes the value zero for feasible points and positive (minimization) for infeasible points. For bilevel problems, the first attempt towards using a penalized approach was made by \cite{AiSh81,AiSh84}. They replaced the lower level problem by a penalized problem; however, the bilevel hierarchy was still maintained and the reduced problem was still difficult to solve. Later double penalty method was introduced in \cite{IsAi92}, where both upper and lower level objective functions were penalized. The problem was reduced into a single level task by replacing the penalized lower level problem with its KKT conditions, and then solving the single level formulation by penalization. 
In a number of studies, the lower level problem is directly replaced by its KKT conditions and then a penalized approach is used to solve the single level problem. Few studies where penalty function approach has been used for linear bilevel problems are \cite{WhAn93,lv2007penalty}. In \cite{WhAn93}, the authors convert the linear bilevel program into a penalized bilinear optimization problem, and then solve a series of bilinear problems to find the optimum. In \cite{lv2007penalty}, the authors reduce the linear bilevel program into single level using KKT conditions, and then append the complementary slackness condition to the upper level objective function with a penalty. The penalized problem is then handled using a series of linear programs.

\subsection{Trust-region Methods}
In trust-region methods, the algorithms approximate a certain region of the objective function with a model function. The region is expanded if the approximation is good, otherwise it is contracted. The first study using trust-region method to solve non-linear bilevel programs was presented in \cite{liu1998trust}, where the lower level problem had a convex objective function and linear constraints. However, no upper level constraints were considered. Later, a more general idea was proposed in \cite{marcotte2001trust}, where the authors locally approximate the bilevel program with a model involving a linear program at the upper level and a linear variational inequality at the lower level. Trust-region and line search ideas have been combined to approach the bilevel optimum over iterations. Similarly, in \cite{colson2005trust}, the authors approximate the bilevel program around an iterate with a model that itself is a linear-quadratic bilevel program. The authors propose to solve the linear-quadratic bilevel program using a mixed integer solver after reducing it to a single level problem using its lower level KKT conditions. Convergence is achieved by sequentially solving linear-quadratic bilevel models.

{\color{black} Next, we discuss about evolutionary algorithms for bilevel optimization. At this point, we would like to refer the readers to other review papers \cite{colson,dempe2003,ViCa94,kalashnikov2015bilevel,sinha2017evolutionary} and books \cite{bilevel-book,ShIsBa97,dempe02,talbi2013metaheuristics,dempe2015bilevel} on bilevel optimization.}

\section{Evolutionary Approaches}

\subsection{Nested Methods}
Nested evolutionary algorithms are a popular approach to handle bilevel problems, where lower level optimization problem is solved corresponding to each and every upper level member~\cite{my-caor14}. Though effective, nested strategies are computationally very expensive and not viable for large scale bilevel problems. Nested methods in the area of evolutionary algorithms have been used in primarily two ways. The first approach has been to use an evolutionary algorithm at the upper level and a classical algorithm at the lower level, while the second approach has been to utilize evolutionary algorithms at both levels. Of course, the choice between two approaches is determined by the complexity of the lower level optimization problem. 

One of the first evolutionary algorithms for solving bilevel optimization problems was proposed in the early 1990s. Mathieu et al. \cite{mathieu} used a nested approach with genetic algorithm at the upper level, and linear programming at the lower level. Another nested approach was proposed in \cite{yin-bilevel}, where the upper level was an evolutionary algorithm and the lower level was solved using Frank-Wolfe algorithm (reduced gradient method) for every upper level member. The authors demonstrated that the idea can be effectively utilized to solve non-convex bilevel optimization problems. 

Nested particle swarm optimization (PSO) was used in \cite{li06} to solve bilevel optimization problems. The effectiveness of the technique was shown on a number of standard test problems with small number of variables, but the computational expense of the nested procedure was not reported. A hybrid approach was proposed in \cite{li07b}, where  simplex-based crossover strategy was used at the upper level, and the lower level was solved using one of the classical approaches. The authors report the generations and population sizes required by the algorithm that can be used to compute the upper level function evaluations, but they do not explicitly report the total number of lower level function evaluations, which presumably is high. 

Differential evolution (DE) based approaches have also been used, for instance, in \cite{zhu2006hybrid}, authors used DE at the upper level and relied on the interior point algorithm at the lower level; similarly, in \cite{angelo13} authors have used DE at both levels. Authors have also combined two different specialized evolutionary algorithms to handle the two levels, for example, in \cite{angelo2015study} authors use an ant colony optimization to handle the upper level and differential evolution to handle the lower level in a transportation routing problem. Another nested approach utilizing ant colony algorithm for solving a bilevel model for production-distribution planning is \cite{calvete2011bilevel}. Scatter search algorithms have also been employed for solving production-distribution planning problems, for instance \cite{camacho2015heuristic}.

{\color{black}
Through a number of approaches involving evolutionary algorithms at one or both levels, the authors have demonstrated the ability of their methods in solving problems that might otherwise be difficult to handle using classical bilevel approaches. However, as already stated, most of these approaches are practically non-scalable. With increasing number of upper level variables, the number of lower level optimization tasks required to be solved increases exponentially. Moreover, if the lower level optimization problem itself is difficult to solve, numerous instances of such a problem cannot be solved, as required by these methods.
}

\subsection{Single-level Reduction}
The idea behind single-level reduction, in the context of evolutionary algorithms, is similar to the the discussions in Section \ref{sec:singleLevelReduction1}. A number of researchers in the area of evolutionary computation have also used the KKT conditions of the lower level to reduce the bilevel problem into a single-level problem. Most often, such an approach is able to solve problems that adhere to certain regularity conditions at the lower level because of the requirement of the KKT conditions. However, as the reduced single-level problem is solved with an evolutionary algorithm, usually the upper level objective function and constraints can be more general and not adhering to such regularities. For instance, one of the earliest papers using such an approach is by Hejazi et al. \cite{hejazi2002linear}, who reduced the linear bilevel problem to single-level and then used a genetic algorithm, where chromosomes emulate the vertex points, to solve the problem. Wang et al. \cite{wang05} reduced the bilevel problem into a single-level optimization problem using KKT conditions, and then utilized a constraint handling scheme to successfully solve a number of standard test problems. Their algorithm was able to handle non-differentiability at the upper level objective function, but not elsewhere. Later on, Wang et al. \cite{wang11} introduced an improved algorithm that was able to handle non-convex lower level problem and performed better than the previous approach \cite{wang05}. {\color{black}However, the number of function evaluations in both approaches remained quite high (requiring function evaluations to the tune of 100,000 for 2 to 5 variable bilevel problems).} In \cite{wang2008genetic}, the authors used a simplex-based genetic algorithm to solve linear-quadratic bilevel problems after reducing it to a single level task. More recently, Jiang et al. \cite{jiang2013application} reduced the bilevel optimization problem into a non-linear optimization problem with complementarity constraints, which is sequentially smoothed and solved with a PSO algorithm. Along similar lines of using lower level optimality conditions, Li \cite{li2015genetic} solved a fractional bilevel optimization problem by utilizing optimality results of the linear fractional lower level problem. In \cite{wan2013hybrid}, the authors embed the chaos search technique in PSO to solve single-level reduced problem.

\subsection{Metamodeling-based Methods}
Metamodeling-based solution methods are commonly used for optimization problems \cite{wang2007review}, where actual function evaluations are expensive. A meta-model or surrogate model is an approximation of the actual model that is relatively quicker to evaluate. Based on a small sample from the actual model, a surrogate model can be trained and used subsequently for optimization. Given that, for complex problems, it is hard to approximate the entire model with a small set of sample points, researchers often resort to iterative meta modeling techniques, where the actual model is approximated locally during iterations.

\begin{figure*}[hbt]
\epsfig{file=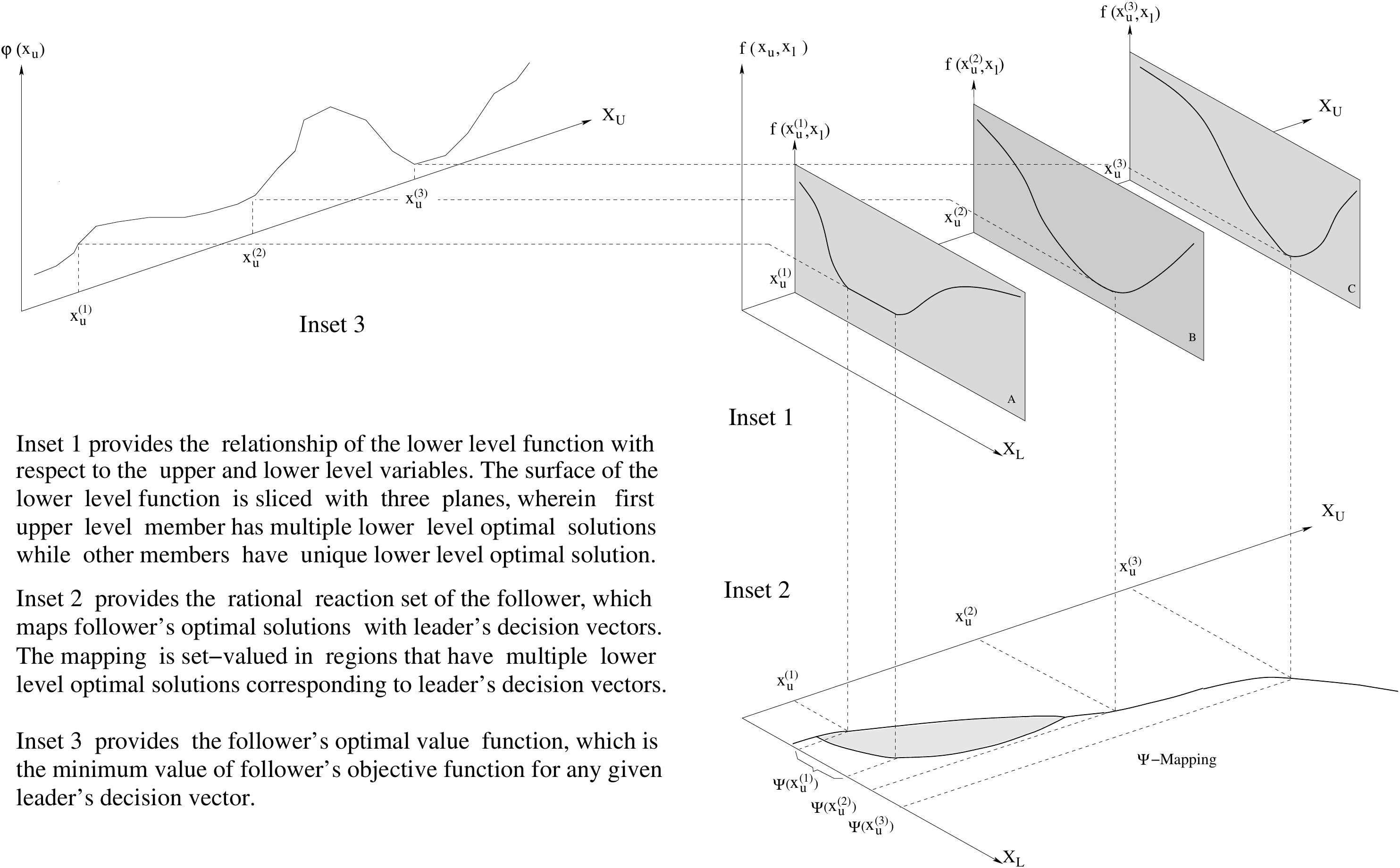,width=0.96\linewidth}
\begin{center}
\captionsetup{justification=centering}
\caption{Graphical representation of rational reaction set ($\Psi$) and lower level optimal value function ($\varphi$).}\label{fig:explain}
\end{center}
\end{figure*}

Bilevel optimization problems contain an inherent complexity that leads to a requirement of large number of evaluations to solve the problem. Metamodeling, when used with population-based algorithms, offers a viable means to handle bilevel optimization problems. 
In this subsection, we discuss four ways in which metamodeling can be applied to bilevel optimization. The discussion related to approximation of the rational reaction set and lower level optimal value function is a review of some recent work. However, before starting, we refer the readers to Figure~\ref{fig:explain}, which provides an understanding of these two mappings graphically for a hypothetical bilevel problem. We also provide a  brief discussion on approximating the bilevel problem with an auxiliary problem.  

\subsubsection{Reaction set mapping}
One of the approaches to solve bilevel optimization problems using evolutionary algorithms would be through iterative approximation of the reaction set mapping $\Psi$. If the $\Psi$-mapping (introduced in Table \ref{tab:notations}) in a bilevel optimization problem is known, it effectively reduces the problem to single level optimization. However, this mapping is seldom available; therefore, the approach could be to solve the lower level problem for a few upper level members and then utilize the lower level optimal solutions and corresponding upper level members to generate an approximate mapping $\hat{\Psi}$. It is noteworthy that approximating a set-valued $\Psi$-mapping offers its own challenges and is not a straightforward task. Assuming that an approximate mapping, $\hat{\Psi}$, can be generated, the following single level optimization problem can be solved for a few generations of the algorithm before deciding to further refine the reaction set.
\begin{align*}
\minimize_{x_u\in X_U, x_l\in X_L}\quad & F(x_u,x_l) \\
\st\quad  & \\
 & \hspace{-12mm} x_l \in \hat{\Psi}(x_u) \\
 & \hspace{-12mm} G_k(x_u,x_l)\leq 0, k=1,\dots,K
\end{align*}
{\color{black}Evolutionary algorithms that rely on this idea to solve bilevel optimization problems are \cite{my-bleaq-arxiv13,my-ejor16b,my-cec14,angelo2014differential}. In \cite{my-bleaq-arxiv13,my-cec14} authors have used quadratic approximation to approximate the local reaction set. This helps in saving lower level optimization calls when the approximation for the local reaction set is good. In case the approximations generated by the algorithm are not acceptable, the method defaults to a nested approach. It is noteworthy that a bilevel algorithm that uses a surrogate model for reaction set mapping may need not be limited to quadratic models but other models can also be used. }


\subsubsection{Optimal lower level value function}
Another way to use metamodeling would be through the approximation of the optimal value function $\varphi$. If the $\varphi$-mapping (introduced in Table \ref{tab:notations}) is known, the bilevel problem can once again be reduced to single level optimization problem as follows \cite{ye2010new},
\begin{align*}
\minimize_{x_u\in X_U, x_l\in X_L}\quad & F(x_u,x_l) \\
\st\quad  & \\
 & \hspace{-12mm} f(x_u,x_l) \le \varphi(x_u) \\
 & \hspace{-12mm} g_j(x_u,x_l)\leq 0, j=1,\dots,J\\
 & \hspace{-12mm} G_k(x_u,x_l)\leq 0, k=1,\dots,K.
\end{align*}
However, since the value function is seldom known, one can attempt to approximate this function using metamodeling techniques. The optimal value function is a single-valued mapping; therefore, approximating this function avoids the complexities associated with set-valued mapping. As described previously, an approximate mapping $\hat{\varphi}$, can be generated with the population members of an evolutionary algorithm and the following single level optimization problem can be solved with refinements at every few generations.
\begin{align*}
\minimize_{x_u\in X_U, x_l\in X_L}\quad & F(x_u,x_l) \\
\st\quad  & \\
 & \hspace{-12mm} f(x_u,x_l) \le \hat{\varphi}(x_u) \\
 & \hspace{-12mm} g_j(x_u,x_l)\leq 0, j=1,\dots,J\\
 & \hspace{-12mm} G_k(x_u,x_l)\leq 0, k=1,\dots,K.
\end{align*}
An evolutionary approach that relies on this idea can be found in \cite{my-cec16a,my-bleaq2-arxiv17}.

\subsubsection{Bypassing lower level problem}
Another way to use a meta-model in bilevel optimization would be completely by-pass the lower level problem, as follows:
\begin{align*}
\minimize_{x_u\in X_U}\quad & \hat{F}(x_u) \\
\st\quad  & \\
 & \hspace{-12mm} \hat{G}_k(x_u)\leq 0, k=1,\ldots,K.
\end{align*} 
Given that the optimal $x_l$ are essentially a function of $x_u$, it is possible to construct a single level problem by ignoring $x_l$ completely. However, the landscape for such a single level problem can be highly non-convex, disconnected and non-differentiable. Advanced metamodeling techniques might be required to use this approach, which may be beneficial for certain classes of bilevel problems. A training set for the metamodel can be constructed by solving few lower level problems for different $x_u$. Both upper level objective $F$ and constraint set ($G_k$) can then be meta-modeled using $x_u$ alone. Given the complex structure of such a single-level problem, it might be sensible to create such an approximation locally. We are currently pursuing such an approach using artificial neural network as the meta-modeling approach.

\subsubsection{Auxiliary bilevel meta-model}
Building up on the trust-region methods for solving bilevel optimization problems, it is possible to utilize the population members in an evolutionary algorithm to formulate auxiliary bilevel problem(s). The auxiliary bilevel problem(s) may be simple enough to be solved using faster specialized techniques. The population members could then be updated based on the obtained auxiliary solution(s). For the moment, there does not exist any evolutionary algorithm based on this idea, but it may be an interesting direction to pursue in the future.

{\color{black}
\section{Discrete Bilevel Optimization}
In this section, we would like to discuss the contributions made towards solving discrete bilevel optimization problems. The formulation of the bilevel problem remains the same as described in Definitions \ref{def:bilevel1} and \ref{def:bilevel2}, along with one or more variables at either of the levels being discrete. Presence of discrete variables can pose a variety of challenges depending upon, whether the discrete variables are present at upper level, lower level or both levels. In the classical literature, branch-and-bound and branch-and-cut are some of the commonly used techniques to handle discreteness in variables. Most of the work on discrete bilevel optimization employs an extension of these ideas from single-level optimization. To highlight the kind of complexities induced in the presence of discrete variables, we consider a simple linear bilevel problem described in \cite{vicente1996discrete,bilevel-book} to show how the inducible region changes based on the upper, lower or both level variables being discrete or continuous.

Consider the following lower level optimization problem which we use for identifying the inducible region of a bilevel problem:
\begin{align*}
\minimize_{x_l}\quad & x_l \\
\st\quad  & \\
 & \hspace{-12mm} x_u+x_l\leq 2,
\hspace{7mm} -x_u+x_l \leq 2,\\
 & \hspace{-12mm} 5x_u-4x_l\leq 10,
\hspace{2mm} -5x_u-4x_l\leq 10.
\end{align*}
For the above lower level problem there are four possible scenarios based on the variables being continuous or discrete:
\begin{enumerate}
\item Continuous-Continuous bilevel program: Consider $x_u \in \reals$ and $x_l \in \reals$
\item Discrete-Continuous bilevel program: Consider $x_u \in \mathbb{Z}$ and $x_l \in \reals$
\item Discrete-Discrete bilevel program: Consider $x_u \in \mathbb{Z}$ and $x_l \in \mathbb{Z}$
\item Continuous-Discrete bilevel program: Consider $x_u \in \reals$ and $x_l \in \mathbb{Z}$
\end{enumerate}
For each scenario the inducible region can be very different that has been shown through Figure~\ref{fig:inducibleRegion}. The figure clearly demonstrates how a discrete variable at any level of the problem can lead to a disconnected search space. Of course there could also be situations where each level has a mix of continuous and discrete variables. 
\begin{figure}
\begin{center}
\epsfig{file=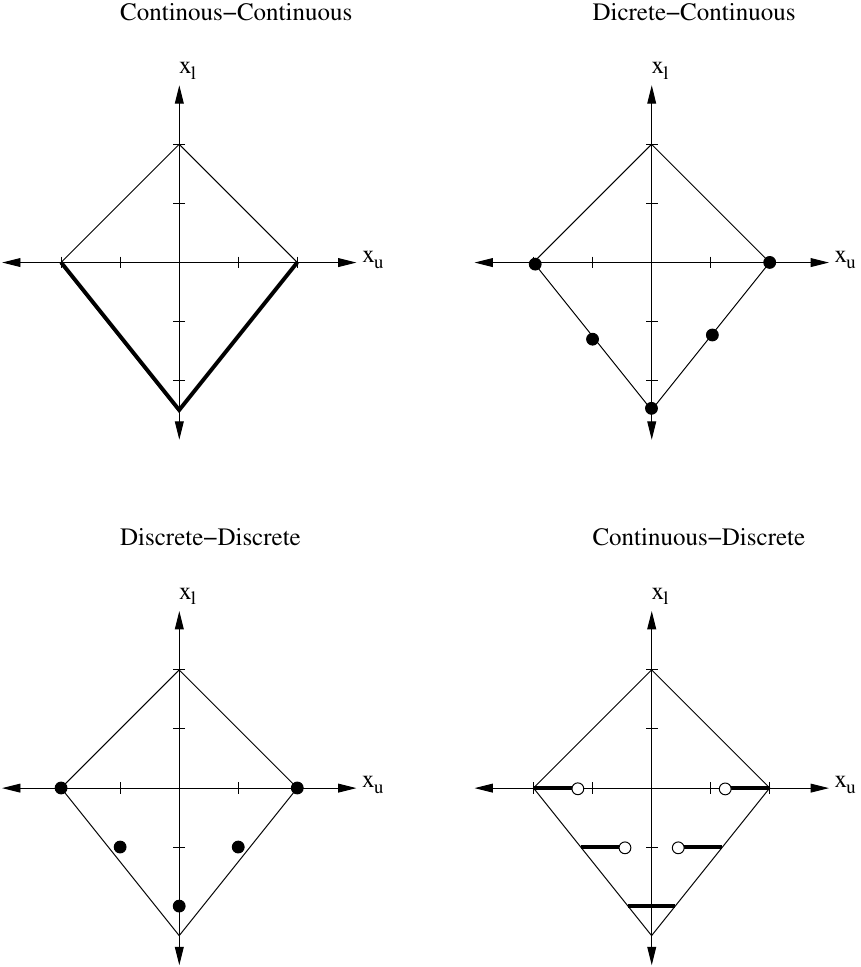,width=0.89\linewidth}
\end{center}
\caption{Inducible region in different cases when the upper and lower level variables belong to continuous or discrete sets.}
\label{fig:inducibleRegion}
\end{figure}

\subsection{Discrete Bilevel Optimization Survey}
One of the early works on discrete bilevel optimization was by Vicente et al. (1996) \cite{vicente1996discrete}, which focused on discrete linear bilevel programs, and analyzed the properties and existence of the optimal solution for different kinds of discretizations arising from the upper and lower level variables. The authors have shown in the paper that certain compactness conditions guarantee the existence of optimal solution in continous-continuous linear bilevel programs, discrete-continuous linear bilevel programs and discrete-discrete linear bilevel programs. The conditions are equivalent to stating that the inducible region is non-empty. However, the existence conditions in the case of continuous-discrete linear bilevel programs are not straightforward. For instance, the inducible region for the continuous-discrete linear bilevel problem in Figure~\ref{fig:inducibleRegion} is a non-compact set that may lead to non-existence of a bilevel optimal solution even when the inducible region is non-empty.

A few studies preceded the study by Vicente et al. (1996) \cite{vicente1996discrete}. For instance, in \cite{moore1990mixed} the authors solved mixed integer linear bilevel problems. The authors pointed out the difficulties involved in fathoming while solving mixed integer bilevel problems using traditional branch-and-bound techniques. Certain fathoming rules used in case of mixed integer linear programming, like fathoming when the relaxed subproblem is worse than the value of the incumbent or fathoming when the solution of the relaxed subproblem is feasible for the mixed integer problem, are not directy applicable to mixed integer linear bilevel problems. Therefore, the authors proposed a branch-and-bound approach involving stricter fathoming conditions. However, the algorithm has a nested structure and is not scalable beyond few integer variables, and to counter which the authors also proposed some heuristics. This study was followed by \cite{bard1992algorithm}, where the same authors solved discrete linear bilevel programs involving only binary variables using an implicit enumeration scheme. In this approach, the authors place a cut, similar to the one used by Bialas and Karwan \cite{bialas1984two} (for the continuous linear bilevel program), seeking incremental improvements in the upper level objective function. A cutting plane method utlizing the {\it Chv$\acute{a}$tal-Gomory} cut for the continuous-discrete bilevel program was proposed in \cite{dempe1996discrete}. Benders-decomposition based techniques have also been employed to solve bilevel problems with mixed integers at the upper level and continuous linear programs at the lower level. The original problem is decomposed into a master and a slave problem. Fixing the integer values converts the slave problem into a bilevel linear program which is solved by KKT based reduction techniques, and  the solution to the slave is utilized to create a cut for the master problem.  The algorithm switches between master and slave problems until the optimality criteria is met. Certain studies in this direction are  \cite{saharidis2009resolution,fontaine2014benders,caramia2015decomposition}.

Despite the attempts made towards algorithm development for discrete bilevel programs, the research is still open for new methods and ideas as none of the proposed techniques would scale well for problems with larger number of variables. Apart from a few nested approaches and KKT-based single level reduction approaches, to our best knowledge, there does not exist any algorithmic study involving evolutionary algorithms for mixed integer bilevel problems that attempt to solve the problem efficiently by utilizing its properties. However, given that the evolutionary approaches are, in particular, potent for handling difficulties such as discreteness and non-differentiabilities they offer a significant scope for solving discrete bilevel optimization problems. Some attempts towards mixed integer bilevel optimization using evolutionary methods are \cite{hecheng2008exponential,arroyo2009genetic,miandoabchi2011optimizing,legillon2012cobra,aksen2013matheuristic,camacho2014solving,handoko2015solving,chaabani2015co}.

There is no dearth of application problems involving discrete variables. While mixed integer bilevel problems are ubiquitous, even combinatorial bilevel programs find innumerable applications in the areas of network design, facility location, hub-and-spoke networks etc. In these areas, these problems are commonly studied in the context of interdiction, protection, robust design, competition and supply chain management, among others. Some of the related applications are highlighted in Section~\ref{sec:applications}.
}

\section{Multiobjective Bilevel Optimization}
In many of the practical problems, a leader and/or the follower might face multiple objectives. This gives rise to multiobjective bilevel optimization problems that we define below.
\begin{definition}\label{def:mbilevel1}
For the upper-level objective function $F:\reals^n\times\reals^m \to\reals^p$ and lower-level objective function $f:\reals^n\times\reals^m \to\reals^q$, the multi-objective bilevel problem is given by 
\begin{align*}
\minimizeQuotes_{x_u\in X_U,x_l\in X_L} & F(x_u,x_l) = (F_1(x_u,x_l),\ldots,F_p(x_u,x_l)) \\
\st &\\  & \hspace{-12mm}x_l\in \argmin_{x_l \in X_L} 
	\lbrace
		f(x_u,x_l) = (f_1(x_u,x_l),\ldots,f_q(x_u,x_l)): \\ & \hspace{7mm}g_j(x_u,x_l)\leq 0, j=1,\dots,J
	\rbrace\\
 & \hspace{-12mm}G_k(x_u,x_l)\leq 0, k=1,\dots,K
\end{align*}
where $G_k:\reals^n\times\reals^m \to \reals$, $k=1,\dots,K$ denote the upper level constraints, and $g_j:\reals^n\times\reals^m \to \reals$ represent the lower level constraints, respectively. The sets $X_U$ and $X_L$ in the definition may denote additional restrictions like integrality.
\end{definition}

The set-valued mapping in this case can be defined as follows and an equivalent definition can be written as in the single-objective case:
\begin{align*}
\Psi(x_u)=&\argmin_{x_l \in X_L}\{f(x_u,x_l) = (f_1(x_u,x_l),\ldots,f_q(x_u,x_l)) : \\ 
& \hspace{13mm} g_j(x_u,x_l)\leq 0, j=1,\dots,J\},
\end{align*}

\subsection{Optimistic vs Pessimistic}
The optimistic or pessimistic position becomes more prominent in multiobjective bilevel optimization. In the presence of multiple objectives at the lower level, the set-valued mapping $\Psi(\cdot)$ normally represents a set of Pareto-optimal solutions corresponding to any given $x_u$, which we refer as follower's Pareto-optimal frontier. A solution to the overall problem (with optimistic or pessimistic position) is expected to produce a trade-off frontier for the leader that we refer as the leader's Pareto-optimal frontier. From the perspective of the leader, it becomes important that what kind of position she seeks to take while solving the problem, as it determines that which solution(s) from the lower level frontier should be considered at the upper level.

Though optimistic positions have commonly been studied in classical \cite{eichfelder2} and evolutionary \cite{my-ecj10} literature in the context of multiobjective bilevel optimization; it is far from realism to expect that the follower will cooperate to an extent that she chooses any point from her Pareto-optimal frontier that is most suitable for the leader. This relies on the assumption that the follower is indifferent to the entire set of optimal solutions, and therefore decides to cooperate. The situation was entirely different in the single-objective case, where, in case of multiple optimal solutions, all the solutions offered an equal value to the follower. However, this can not be assumed in the multiobjective case. Solution to the optimistic formulation in multiobjective bilevel optimization leads to the best possible Pareto-optimal frontier that can be achieved by the leader. Similarly, solution to the pessimistic formulation leads to the worst possible Pareto-optimal frontier at the upper level.

\begin{figure}
\begin{center}
\epsfig{file=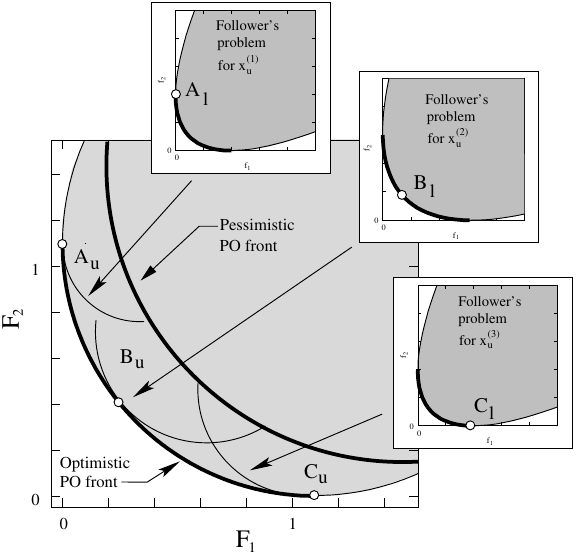,width=0.89\linewidth}
\end{center}
\caption{Leader's Pareto-optimal (PO) frontiers for optimistic and pessimistic positions. Few follower's Pareto-optimal (PO) frontiers are shown (in insets) along with their representations in the leader's objective space.}
\label{fig:multiobjective}
\end{figure}

If the value function or the choice function of the follower is known to the leader, it provides an information as to what kind of trade-off is preferred by the follower. A knowledge of such a function effectively, casually speaking, reduces the lower level optimization problem into a single-objective optimization task, where the value function may be directly optimized. The leader's Pareto-optimal frontier for such intermediate positions lies between the optimistic and the pessimistic frontiers. Figure \ref{fig:multiobjective} shows the optimistic and pessimistic frontiers for a hypothetical multiobjective bilevel problem with two objectives at upper and lower levels. Follower's frontier corresponding to $x_{u}^{1}$, $x_{u}^{2}$ and $x_{u}^{3}$, and her decisions $A_l$, $B_l$ and $C_l$ are shown in the insets. The corresponding representations of the follower's frontier and decisions ($A_u$, $B_u$ and $C_u$) in the leader's space are also shown.

{\color{black}
\subsection{Example}
Below, we provide a bilevel optimization problem involving design of a tax policy \cite{my-cec13}. The upper level in this example is the governemnt that wants to tax the lower level, a mining company, based on the pollution it causes to the environment. The government here has two objectives: the first objective is to maximize the revenues generated by the mining project, which may include the additional jobs, taxes, etc; and the second objective is to minimize the harm caused to the environment as a result of mining. Obviously, there is a trade-off between the two objectives, and the government as a decision maker needs to choose one of the preferred trade-off solutions. The mining company has a sole objective of maximizing its profit under the constraints set by the government. In this scenario, the government would like to have a tax structure such that it is able to maximize its own revenues in addition to being able to restrain the mining company from causing extensive damage to the environment. It is possible for the leader to optimally regulate the problem in its favour, provided that it has complete knowledge of the follower's strategies. The hierarchical optimization problem in this case can be formulated as follows:
\begin{align}
\max_{\tau,q} \quad &\boldF(q, \tau) = (R, - D) \label{eq:object} \\
\mbox{s.t.} \quad & q \in \argmax_{q}
\left\{ \begin{aligned}
\pi(q)& = p(q)q - c(q) - R\\
\pi(q)& \ge 0	\\
\end{aligned} \right\} \label{eq:constr1}\\
& q \ge 0, \tau \ge 0.
\end{align}

In (\ref{eq:object}), the first objective deals with the tax revenue, where $R = \tau q$; $\tau$ is the per unit tax imposed on the mine, and $q$ is the amount of metal extracted from the ore by the follower. The second objective denotes the environmental damage caused by the mine that the government ultimately wants to minimize. $D = kq$, where $k$ is the pollution coefficient signifying the negative impact of extraction on the environment. The damages are thus linear and scale proportionately with the amount of gold extracted from the earth since a larger base of operation implies larger environmental damage.

Equation (\ref{eq:constr1}) gives the profit of the mine, where $p(q)q$ (price function times amount of metal extracted) is the revenue function, and $c(q)$ is the extraction cost function followed by the additional tax levied on the mine. The mine is most likely to be a price taker when it comes to the price of gold and must base its mining decisions on the possible price paid by their customers. It would therefore be plausible to replace the price function for gold in the above equation by a constant. However, given the assumption that the mine can extract a large amount of ore, and subsequently gold, at one time, it might be possible for it to affect the price of gold slightly. Therefore, we assume the price function to be linear with a small slope. Extraction cost can be considered to be quadratic. Thus, we have the following model:
\begin{flalign}
\max_{\tau,q} & \quad F(q, \tau) = (\tau q, - kq) \label{eq:object2}\\
\mbox{s.t.} & \quad \notag\\
& q \in \argmax_{q} 
\left\{ \begin{aligned}
\pi(q)=&(\alpha - \beta q)q - \\ &(\delta q^{2} + \gamma q + \phi) - \tau q\\
\pi(q) \ge& \hspace{1mm} 0	\\
\end{aligned} \right\}\\
&q \ge 0, \tau \ge 0,
\end{flalign}
where $\alpha, \beta, \delta, \gamma, \phi$ are constants, and $\phi$ represents the fixed costs of setting up operations. The above problem can be solved analytically by taking a weighted sum of squares of the upper level objectives ($w\tau q - (1-w)kq: w \in [0,1]$). The optimal solution to the above problem is given as follows:
\begin{align}
& \tau^*(w) = \frac{\alpha - \gamma - k}{2} + \frac{k}{2w},\\
&q^*(w) = \frac{w(\alpha - \gamma)-(1-w)k}{4w(\beta + \delta)}\label{eq:production}.
\end{align}

We assume the parameters as $\alpha=100, \beta=1, \delta=1, \gamma=1$ and $\phi=0$. By varying the government's preference weights ($w$) in its domain, one can generate the entire Pareto-optimal solutions for the leader. Note that a very high taxation (or weight to the environmental objective) may lead to no production at the lower level, for instance, we find that $w < 0.01$ does not lead to any production. The Pareto frontier generated using weights $0.01 \le w \le 1$ has been provided in Figure \ref{fig:analytical} for the above model.
\begin{figure}
	\includegraphics[width=\linewidth]{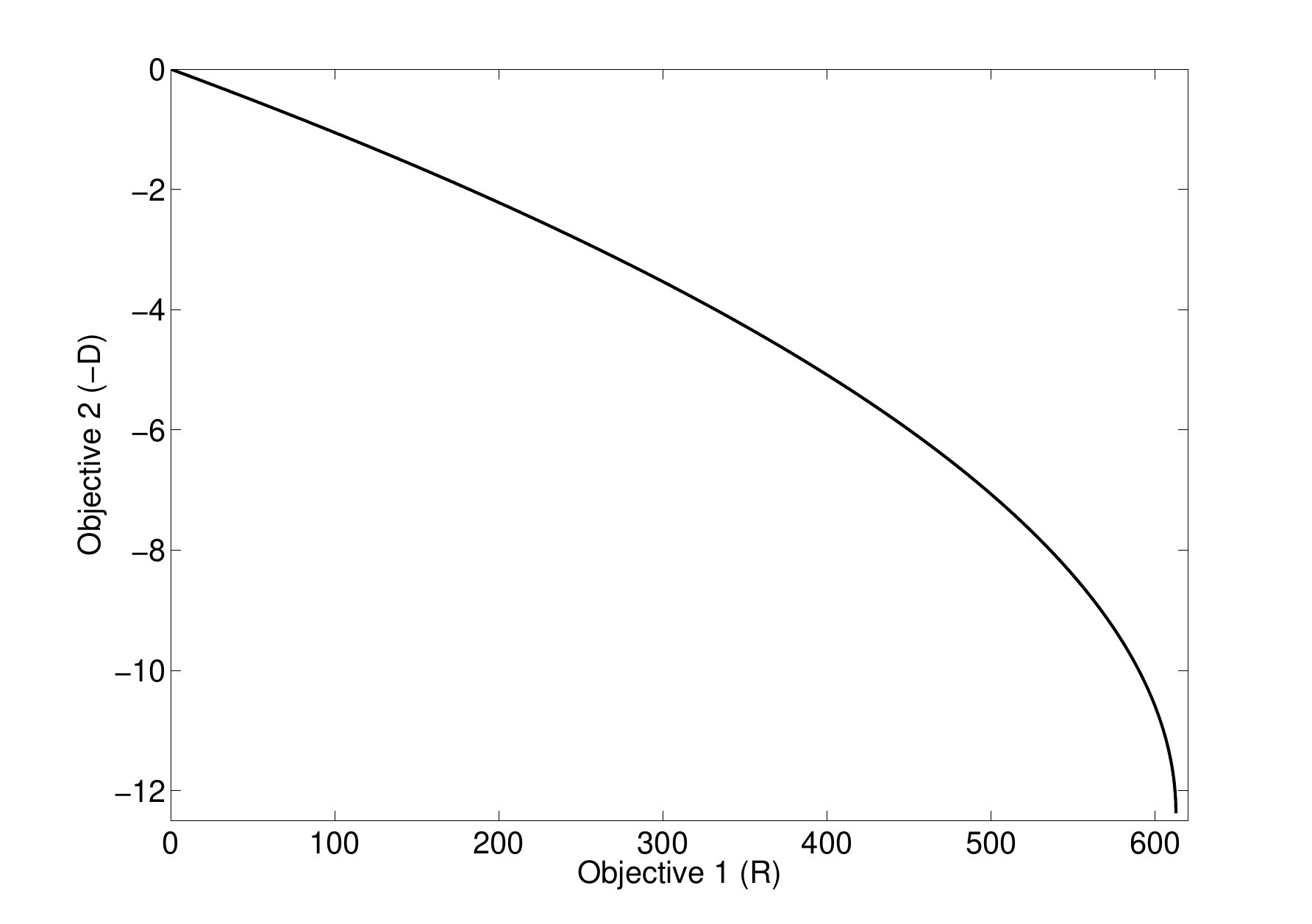}
	\caption[Pareto-optimal Front]{Pareto-optimal frontier for the government showing the trade-off between tax revenues and environmental pollution.}
	\label{fig:analytical}
\end{figure}
}

\subsection{Multiobjective Bilevel Optimization Survey}
There exists a significant amount of work on single objective bilevel optimization; however, little has been done on multi-objective bilevel optimization primarily because of the computational and decision making complexities that these problems offer. For results on optimality conditions in multiobjective bilevel optimization, the readers may refer to \cite{gadhi2012necessary,ye2011necessary,bank1983non}. On the methodology side, Eichfelder \cite{eichfelder,eichfelder2} solved simple multi-objective bilevel problems using a classical approach. The lower level problems in these studies have been solved using a numerical optimization technique, and the upper level problem is handled using an adaptive exhaustive search method. This makes the solution procedure computationally demanding and non-scalable to large-scale problems. In another study, Shi and Xia \cite{shi-xia} used $\epsilon$-constraint method at both levels of multi-objective bilevel problem to convert the problem into an $\epsilon$-constraint bilevel problem. The $\epsilon$-parameter is elicited from the decision maker, and the problem is solved by replacing the lower level constrained optimization problem with its KKT conditions.

One of the first studies, utilizing an evolutionary approach for multiobjective bilevel optimization was by Yin \cite{yin-bilevel}. The study involved multiple objectives at the upper lever, and a single objective at the lower level. The study suggested a nested genetic algorithm, and applied it on a transportation planning and management problem. Later Halter and Mostaghim \cite{halter-sanaz} used a particle swarm optimization (PSO) based nested strategy to solve a multi-component chemical system. The lower level problem in their application was linear for which they used a specialized linear multi-objective PSO approach. Recently, a hybrid bilevel  evolutionary multi-objective optimization algorithm coupled with local search was proposed in \cite{my-ecj10} (For earlier versions, refer \cite{my-cec09a,my-ifac09,my-emo09,my-mcdm09}). In the paper, the authors handled non-linear as well as discrete bilevel problems with relatively larger number of variables. The study also provided a suite of test problems for bilevel multi-objective optimization. 

There has been some work done on decision making aspects at upper and lower levels. For example, in \cite{my-emo11} an optimistic version of multiobjective bilevel optimization, involving interaction with the upper level decision maker, has been solved. The approach leads to the most preferred point at the upper level instead of the entire Pareto-frontier. Since multi-objective bilevel optimization is computationally expensive, such an approach was justified as it led to enormous savings in computational expense. Studies that have considered decision making at the lower level include \cite{my-emo15,my-ieeetec16}. In \cite{my-emo15}, the authors have replaced the lower level with a value function that effectively reduces the lower level problem to single-objective optimization task. In \cite{my-ieeetec16}, the follower's value function is known with uncertainty, and the authors propose a strategy to handle such problems. Other work related to bilevel multi-objective optimization can be found in \cite{pieume09,pramnik11,linnala12,ruuska12,zhang12}.

\section{Applications}\label{sec:applications}
Bilevel optimization commonly appears in many practical problems. They are often encountered in the fields of economics, transportation, engineering and management, among others. The following list will provide an insight to the readers on the relevance of these problems to practice.



\begin{figure*}[hbtp]
\begin{minipage}[t]{0.48\linewidth}
\begin{center}
\includegraphics[height=1.90in]{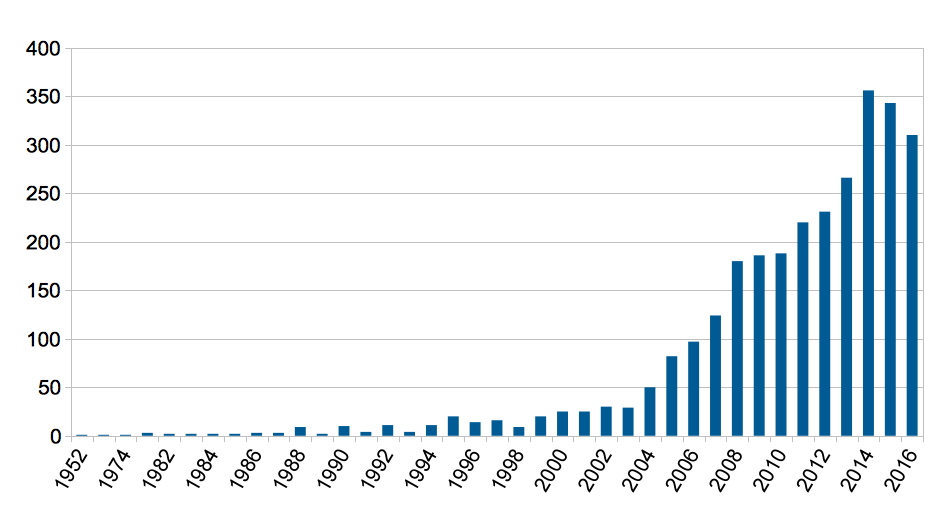}
\end{center}
\caption{Interest on bilevel optimization over time.}
\label{fig:all_growth}
\end{minipage}\hfill
\begin{minipage}[t]{0.48\linewidth}
\begin{center}
\includegraphics[height=1.90in]{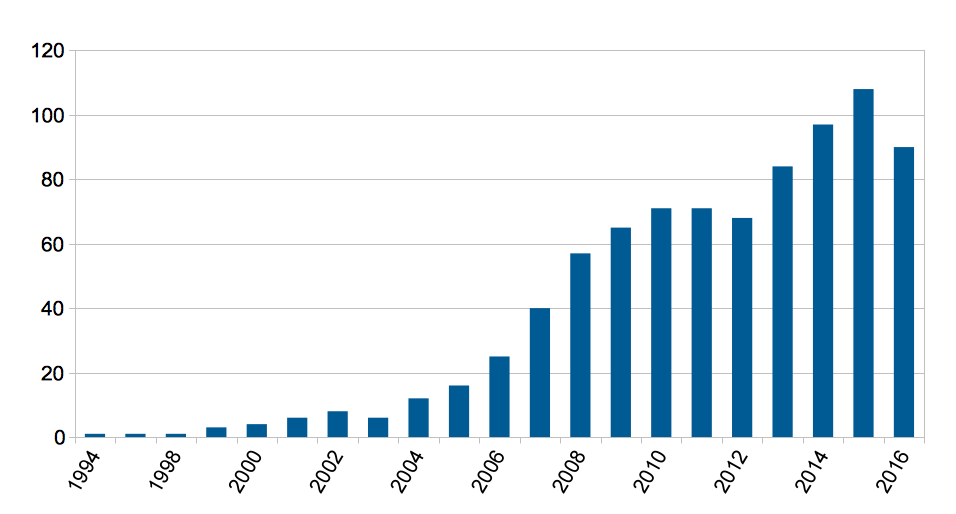}
\end{center}
\caption{Interest on evolutionary bilevel optimization over time.}
\label{fig:ea_growth}
\end{minipage}
\end{figure*}

\begin{enumerate}

\item Toll Setting Problem: Toll-setting problem is essentially a part of network problems. In this problem, there is an authority that wants to optimize the tolls for a network of roads. The authority acts as a leader and the network users act as followers. 
Papers on toll-setting problem and its multi-objective extensions can be found in \cite{brotcorne01,migdalas95,constantin95,labbe1998bilevel,yin2002multiobjective,marcotte2004bilevel,kalashnikov2010comparison,wang2014bilevel,my-cec15a,fan2015optimal,gonzalez2015scatter,kalashnikov2016heuristic}. Bilevel optimization is quite commonly used in network design problems. Instead of going through specific problems, we refer the readers to a variety of applications of bilevel optimization to the area of network design \cite{LeBo86,BeBlBoLe92,MaMa92,yang1998models,yin2000genetic,ceylan2004traffic,yamada2009designing,chen2010stochastic,gallo2010meta,fan2011bi,mesbah2011optimization,xu2012nonlinear,ren2013integrated,camacho2015genetic}.


\item Environmental Economics: Bilevel optimization commonly appears in environmental economics, where an authority wants to tax an organization or individual that is polluting the environment as a result of its operations. Finding an optimal level of tax that offers a compromise between revenues and pollution results in a bilevel optimization problem with the regulator as the leader and the polluting entity as a follower \cite{amouzegar1999determining,my-cec13,my-cec15c,bostian2015valuing,whittaker2016spatial}.


\item Chemical Industry: In chemical industries, the chemists often face a bilevel optimization problem where they have to decide upon the conditions (state variables and quantity of reactants) for the reaction to achieve optimal output. While optimizing the output is the upper level problem, the lower level appears as an equilibrium condition, which is an entropy functional minimization problem. Such applications of bilevel optimization can be found in \cite{smith82,clark1990bilevel,raghunathan2003mathematical,halter2006bilevel}


\item Optimal design: Bilevel problems are very common in structural optimization or optimal shape design. For instance, in structural optimization one often requires to minimize the weight or cost of a structure as an upper level objective with the decision variables as shape of the structure, choice of materials, amount of material etc. The constraints at the upper level involve bounds on displacements, stresses and contact forces whose values are determined by solving the potential energy minimization problem at the lower level. The equilibrium condition in many optimal shape design problems appears in the form of variational inequalities which require the overall problem to formulated as a two level task. For optimal design applications the readers may refer to \cite{herskovits2000contact,christiansen2001stochastic,bendsoe95,kovcvara1997topology, kocvara1995solution}.

\item Defense applications: Bilevel optimization has a number of applications in the defense sector \cite{bracken1974defense}, for example attacker-defender Stackelberg games \cite{israeli2002shortest,brown2006defending,scaparra2008bilevel,o2011designing,an13,ramamoorthy2016hub,ramamoorthy2017hub}.  Specifically, some recent applications include planning the prepositioning of defensive missile interceptors to counter an attack threat \cite{brown05}, interdicting nuclear weapons project \cite{brown09}, homeland security applications \cite{wein09,lowe06}, location problems \cite{aksen2012bilevel}. The bilevel problem, while offending, involves maximizing the damage caused to the opponent by taking into account the optimal reactions of the opponent. Conversely, while defending, the bilevel problem involves minimizing the maximum damage that an attacker can cause.

\item Facility Location: Facility location problems may take the form of a Stackelberg game if a firm, while locating its facility, decides to account the actions of its competitors. For instance, in \cite{kuccukaydin2011competitive} the authors study the scenario where a firm enters a market by locating new facilities, and its competitor reacts by adjusting the attractiveness of its existing facilities. Another study considers location of logistics distribution centers by minimizing the planners' cost at the upper level and customers' cost at the lower level \cite{sun2008bi}. Other applications of bilevel optimization to facility location problem may be found in
\cite{jin2007bi,uno2008evolutionary,sun08,alekseeva2009hybrid,calvete2013efficient,camacho2014solving,panin2014bilevel,camacho2014solving,maric2014metaheuristic,caramia2015decomposition,maldonado2016analyzing}.

\item Inverse Optimal control: Inverse optimal control problems are essentially bilevel in nature \cite{mombaur2010human,albrecht2011imitating,johnson2013inverse,my-cec16b} with wide applications in robotics, computer vision, communication theory and remote sensing to name a few.
One of the major challenges in control theory is deriving the performance index or reward function which fits best on a given dataset. Such tasks lie in the category of inverse optimal control theory, where one solicits the calculation of the cause based on the given result. Such a requirement necessitates solving a parameter estimation problem with an optimal control problem.

\item Machine learning: Most of the machine learning and evolutionary optimization techniques often involve a number of parameters. A proper choice of these parameters has a substantial effect on the accuracy and efficiency of the approach. Tuning of these parameters is often achieved using brute force strategies, such as grid search and random search. A bilevel formulation of this problem allows for systematic and more efficient search when the number of parameters are large. Some of the approaches that have acknowledged the bilevel nature of this problem are \cite{bennett2006model,bennett2008bilevel,my-gecco14,liang2015evolutionary}.

\item Principal-agent problems: Principal-agent problem \cite{laffont2009theory} is a classical problem in economics, where a principal (leader) sub-contracts a job to an agent (follower). Given that the agent prefers to act in his own interests rather than those of the principal, it becomes important for the principal to have an incentive scheme that aligns the interests of the agent with the principal. Design of such contracts appears as a bilevel optimization problem. In real life, principal-agent relationships are commonly found in doctor-patient, senior management-lower management, employer-employee, corporate board-shareholders and politician-voters scenarios. Studies that the readers may refer to are \cite{van1993principal,garen1994executive,Bilevel-linear,cecchini2013solving,xu2007supply}.
\end{enumerate}

\section{Interest over time}
In this section, we perform a text analysis of papers on bilevel programming (and Stackelberg games) that are indexed in SCOPUS. To begin with, we analyze the volume of publications every year on bilevel programming since 1950s to present, and then closely look at the themes within bilevel programming that have contributed to the growth over the years. The themes were discovered using a non-parametric Bayesian approach \cite{teh2006hierarchical}, which clusters the documents together based on similarities. The documents may probabilistically belong to multiple clusters at the same time.

Figure \ref{fig:all_growth} shows how the interest on bilevel programming has been growing at a slow pace until early 2000s and then picked-up significantly at the middle of the previous decade. The studies on bilevel programming using evolutionary algorithms appeared for the first time during the mid 1990s that took another decade to pick up to the extent that almost $10\%$ of all studies on bilevel optimization utilize evolutionary methods. Figure \ref{fig:all_growth} shows the growth of evolutionary methods in the context of bilevel optimization from the 1990s to present. While the early papers on bilevel programming (pre-2000) were mainly focused on solution methods and optimality conditions, the growth in the post-2000 period was fueled by papers on applications of bilevel programming.

\begin{figure}[t]
\begin{subfigure}[t]{0.41\linewidth}
\begin{center}
\epsfig{file=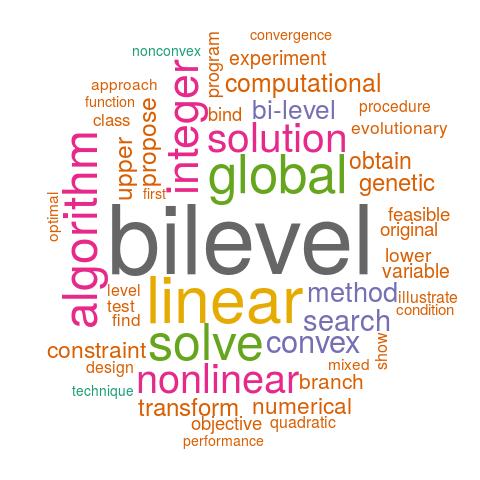,width=\linewidth}
\end{center}
\label{fig:cloud0}
\end{subfigure}\hfill
\begin{subfigure}[t]{0.59\linewidth}
\begin{center}
\epsfig{file=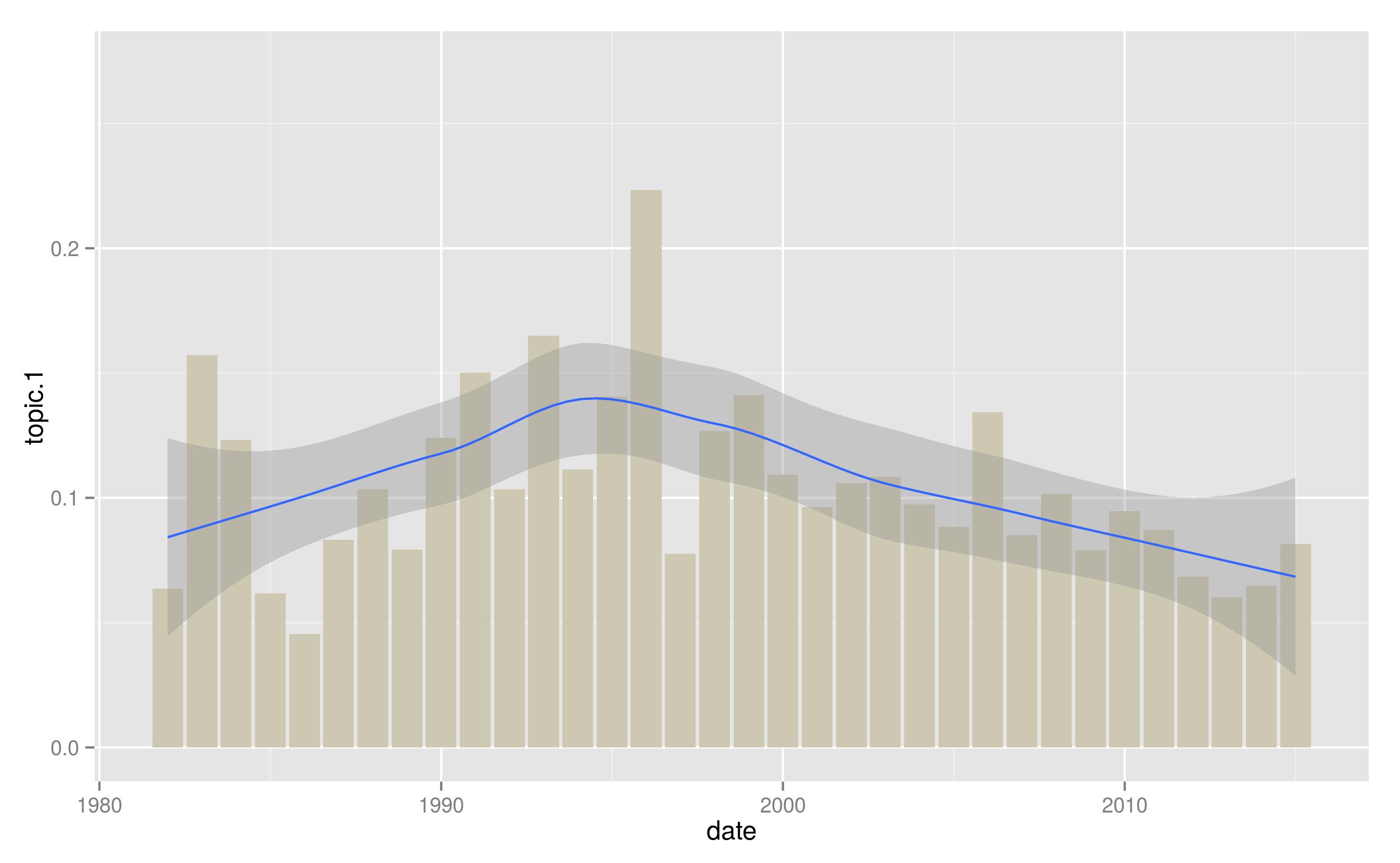,width=\linewidth} 
\end{center}
\label{fig:interest0}
\end{subfigure}
\caption{Topic: Methods.}
\label{fig:first}

\begin{subfigure}[t]{0.41\linewidth}
\begin{center}
\epsfig{file=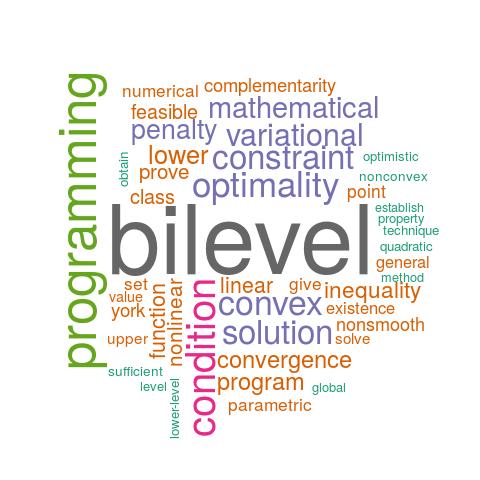,width=\linewidth}
\end{center}
\label{fig:cloud1}
\end{subfigure}\hfill
\begin{subfigure}[t]{0.59\linewidth}
\begin{center}
\epsfig{file=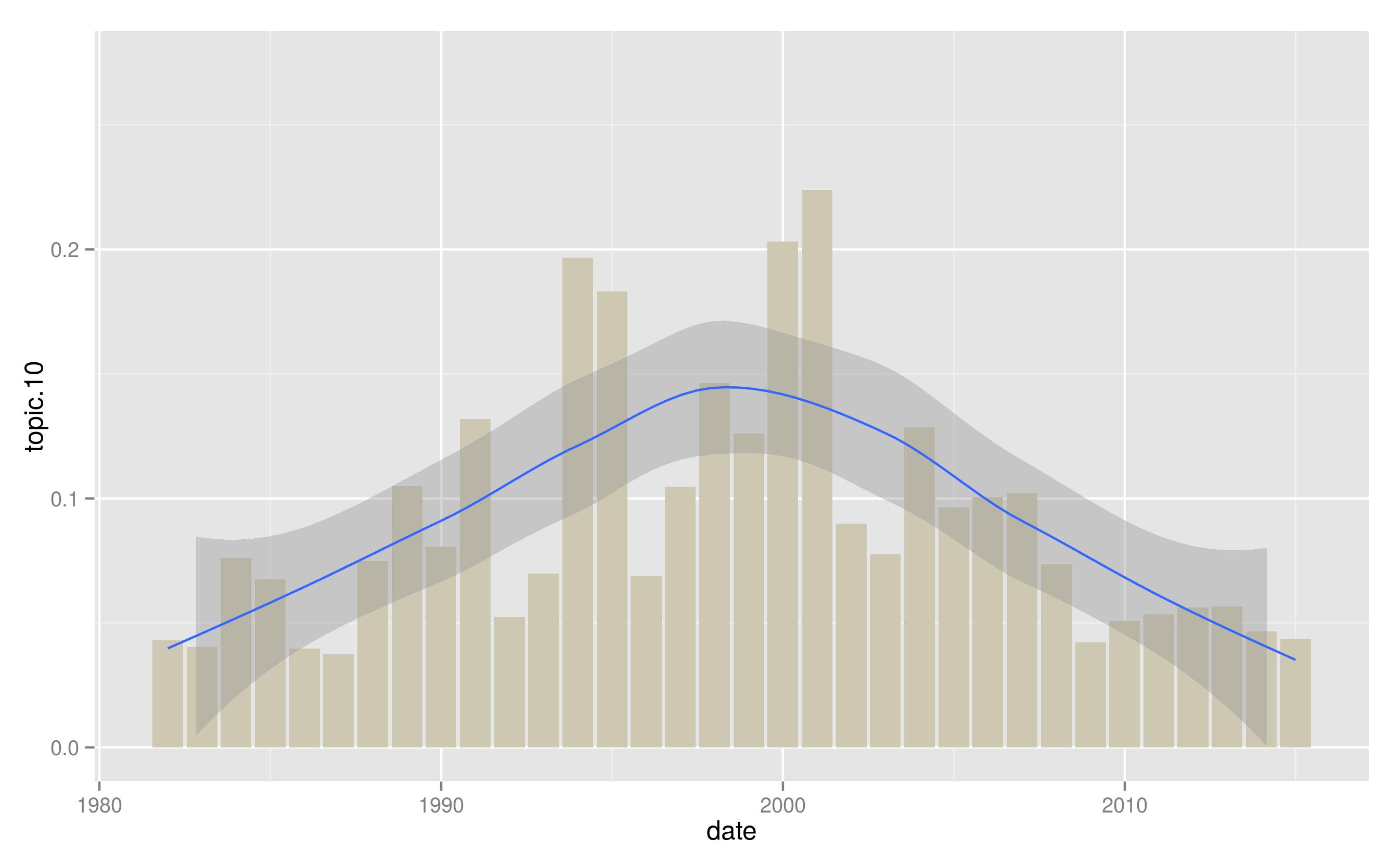,width=\linewidth} 
\end{center}
\label{fig:interest1}
\end{subfigure}
\caption{Topic: Optimality conditions.}

\begin{subfigure}[t]{0.41\linewidth}
\begin{center}
\epsfig{file=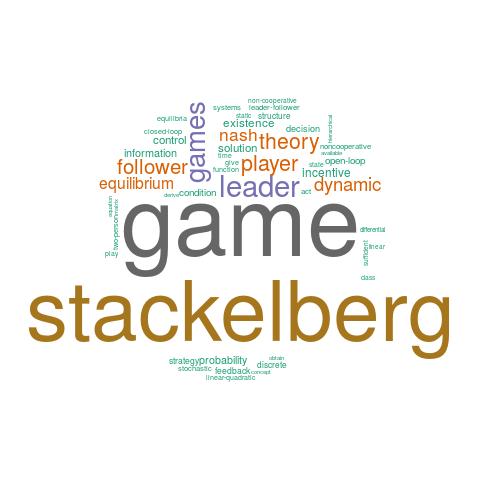,width=\linewidth}
\end{center}
\label{fig:cloud2}
\end{subfigure}\hfill
\begin{subfigure}[t]{0.59\linewidth}
\begin{center}
\epsfig{file=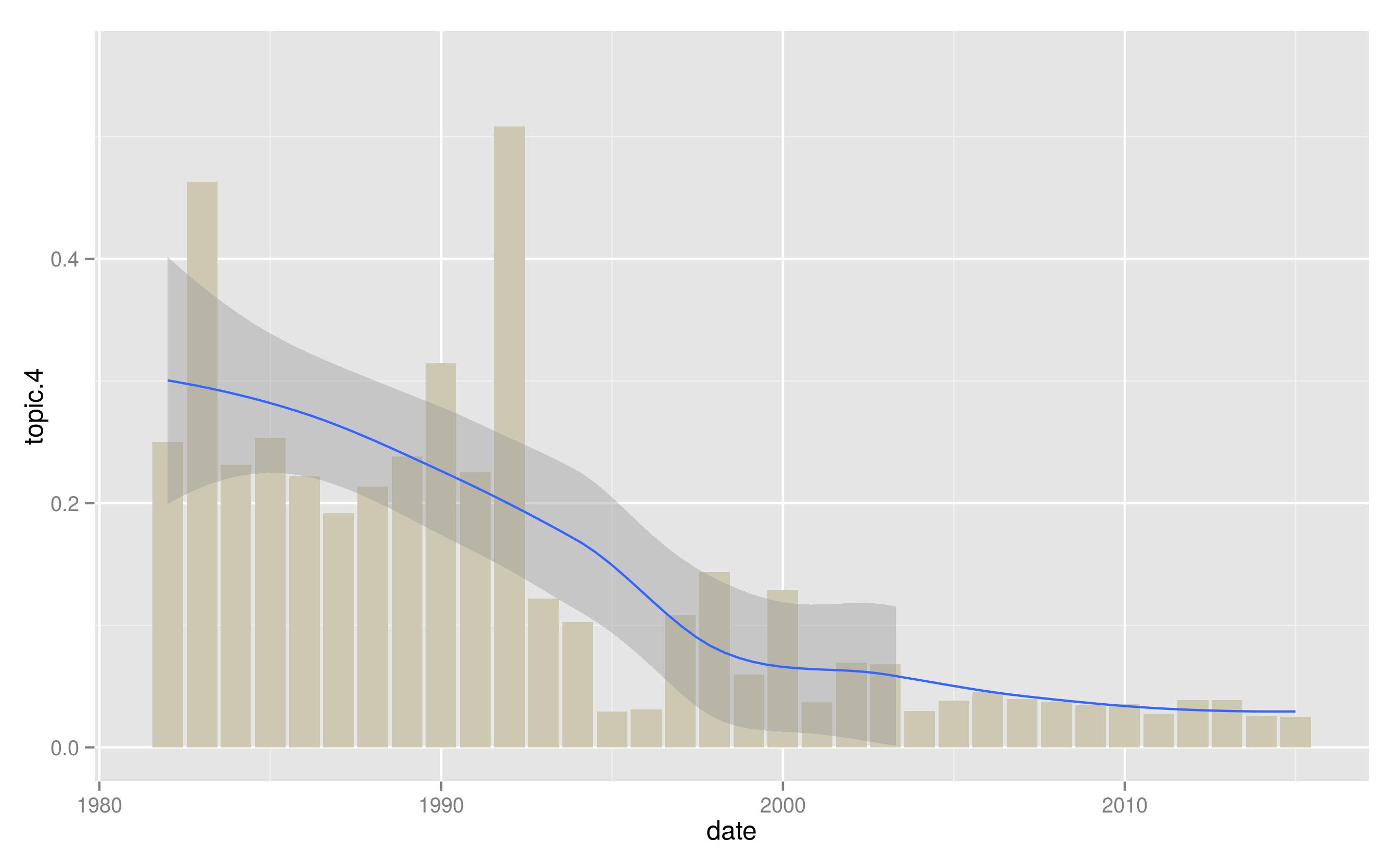,width=\linewidth} 
\end{center}
\label{fig:interest2}
\end{subfigure}
\caption{Topic: Classical game theory.}

\begin{subfigure}[t]{0.41\linewidth}
\begin{center}
\epsfig{file=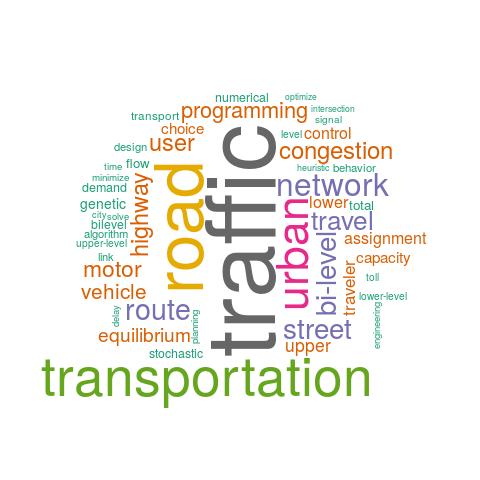,width=\linewidth}
\end{center}
\label{fig:cloud3}
\end{subfigure}\hfill
\begin{subfigure}[t]{0.59\linewidth}
\begin{center}
\epsfig{file=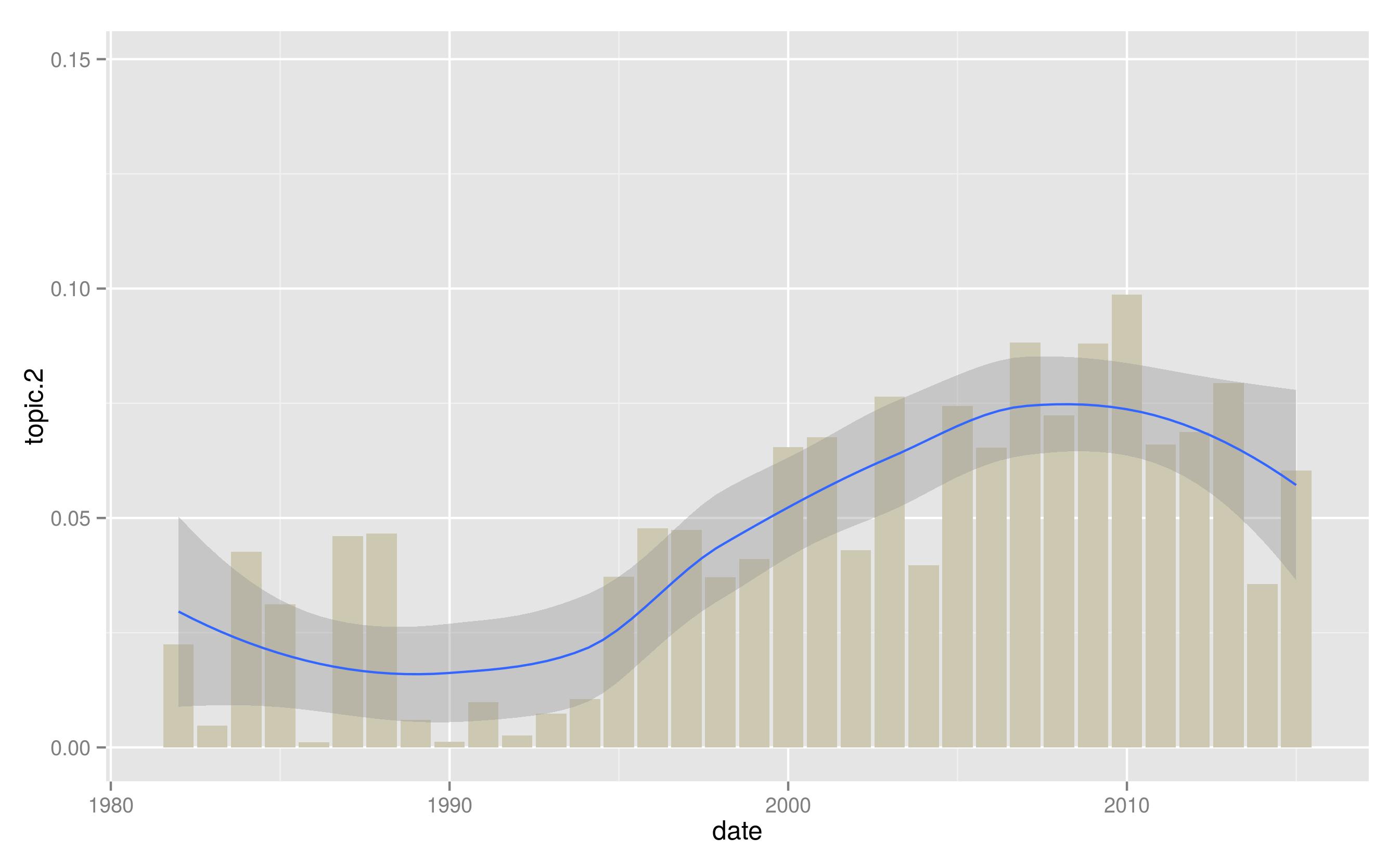,width=\linewidth} 
\end{center}
\label{fig:interest3}
\end{subfigure}
\caption{Topic: Network design.}
\end{figure}

\begin{figure}[t]
\begin{subfigure}[t]{0.41\linewidth}
\begin{center}
\epsfig{file=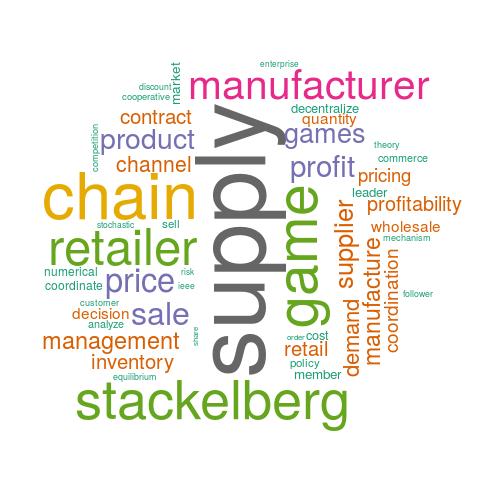,width=\linewidth}
\end{center}
\label{fig:cloud4}
\end{subfigure}\hfill
\begin{subfigure}[t]{0.59\linewidth}
\begin{center}
\epsfig{file=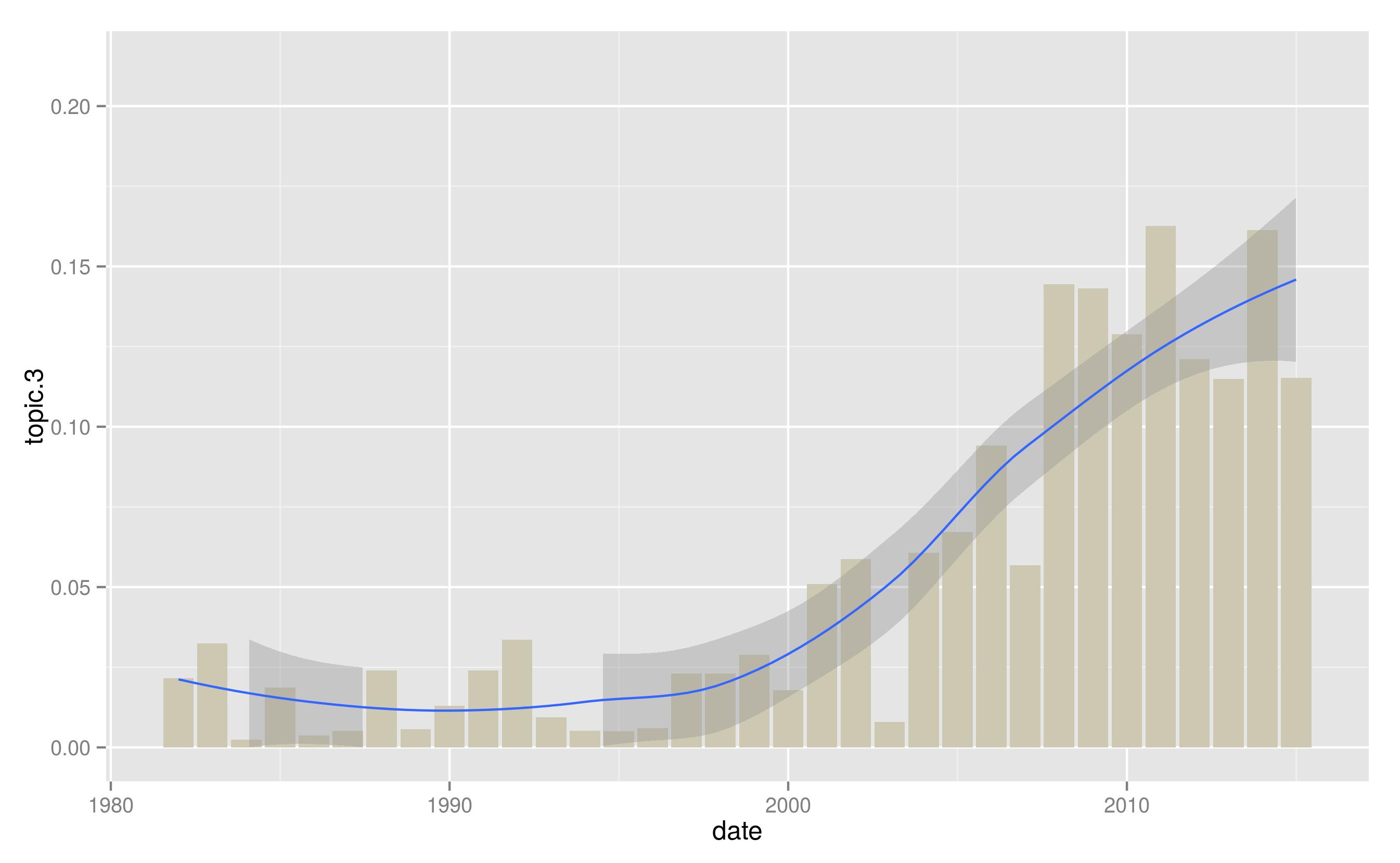,width=\linewidth} 
\end{center}
\label{fig:interest4}
\end{subfigure}
\caption{Topic: Supply chain applications.}

\begin{subfigure}[t]{0.41\linewidth}
\begin{center}
\epsfig{file=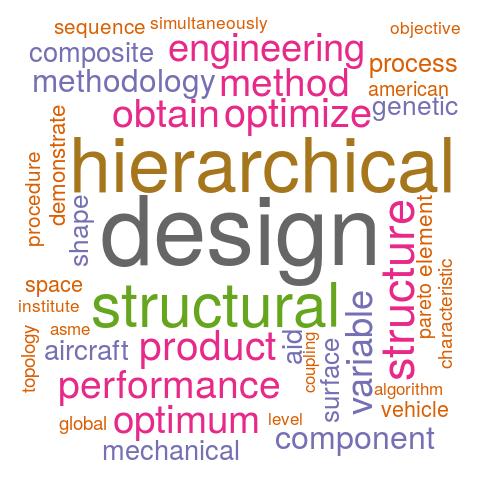,width=\linewidth}
\end{center}
\label{fig:cloud5}
\end{subfigure}\hfill
\begin{subfigure}[t]{0.59\linewidth}
\begin{center}
\epsfig{file=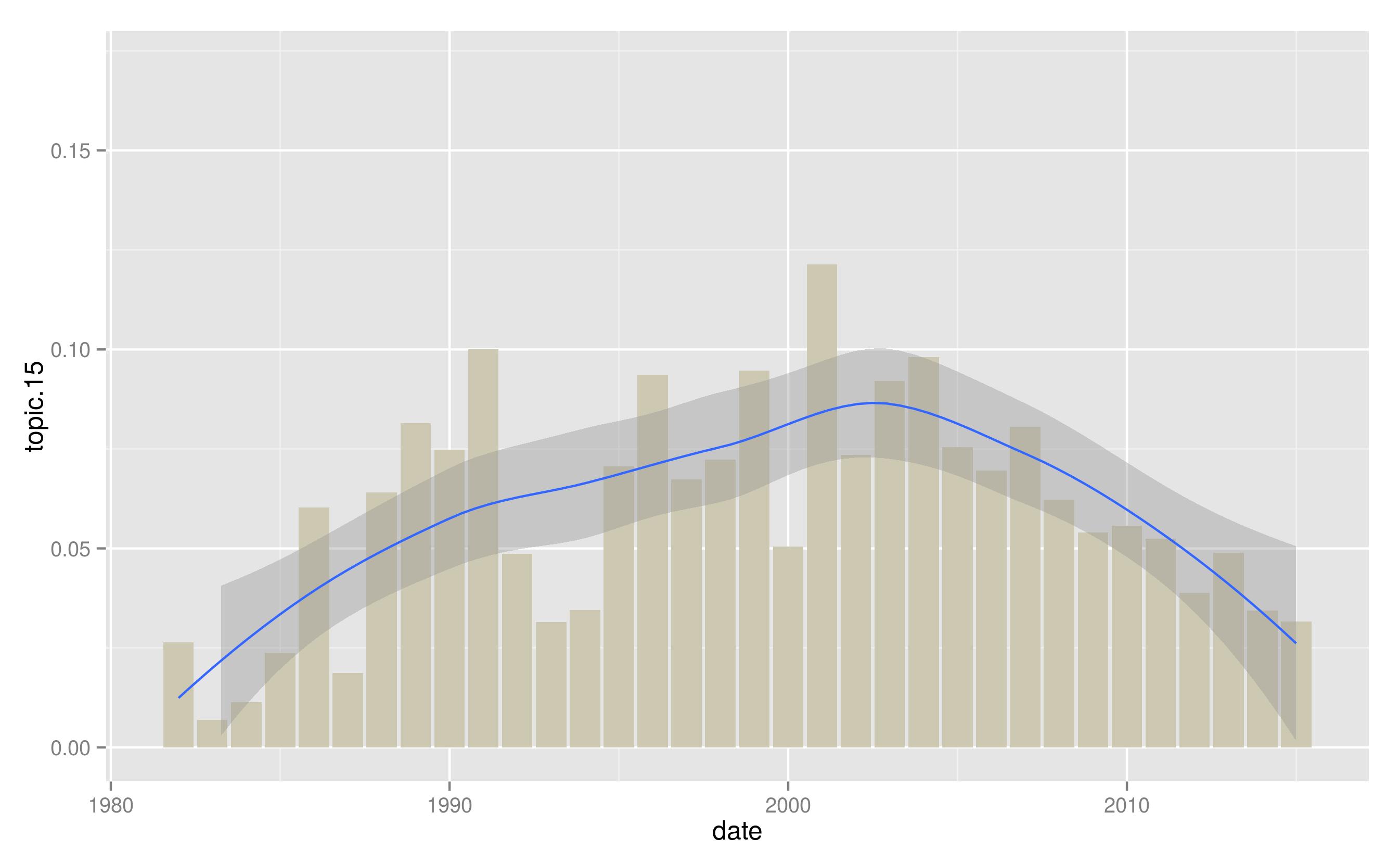,width=\linewidth} 
\end{center}
\label{fig:interest5}
\end{subfigure}
\caption{Topic: Optimal design applications.}

\begin{subfigure}[t]{0.41\linewidth}
\begin{center}
\epsfig{file=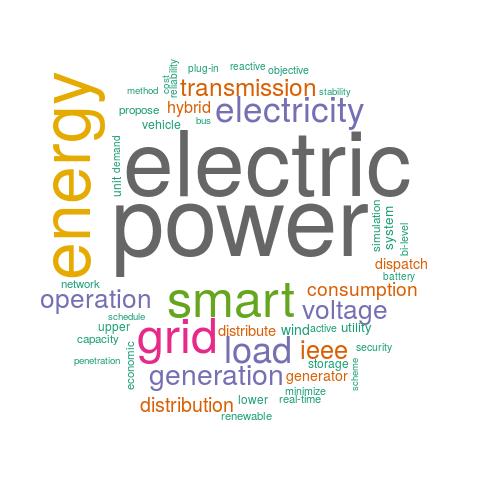,width=\linewidth}
\end{center}
\label{fig:cloud6}
\end{subfigure}\hfill
\begin{subfigure}[t]{0.59\linewidth}
\begin{center}
\epsfig{file=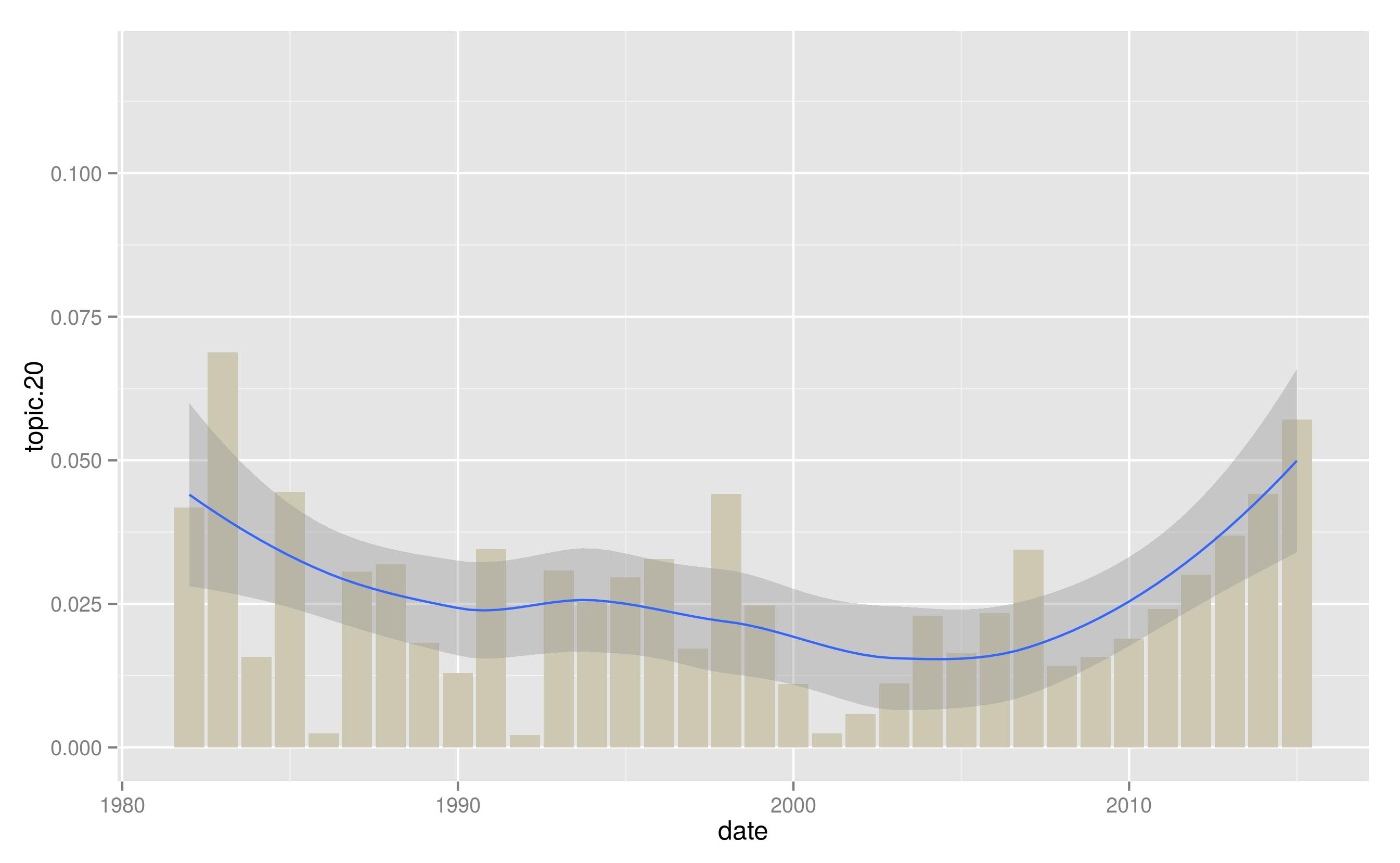,width=\linewidth} 
\end{center}
\label{fig:interest6}
\end{subfigure}
\caption{Topic: Electricity transmission applications.}

\begin{subfigure}[t]{0.41\linewidth}
\begin{center}
\epsfig{file=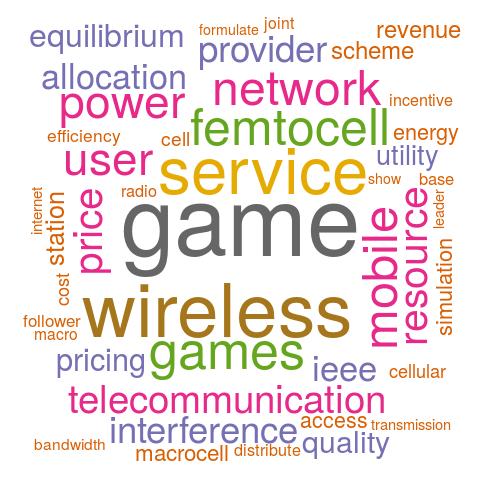,width=\linewidth}
\end{center}
\label{fig:cloud7}
\end{subfigure}\hfill
\begin{subfigure}[t]{0.59\linewidth}
\begin{center}
\epsfig{file=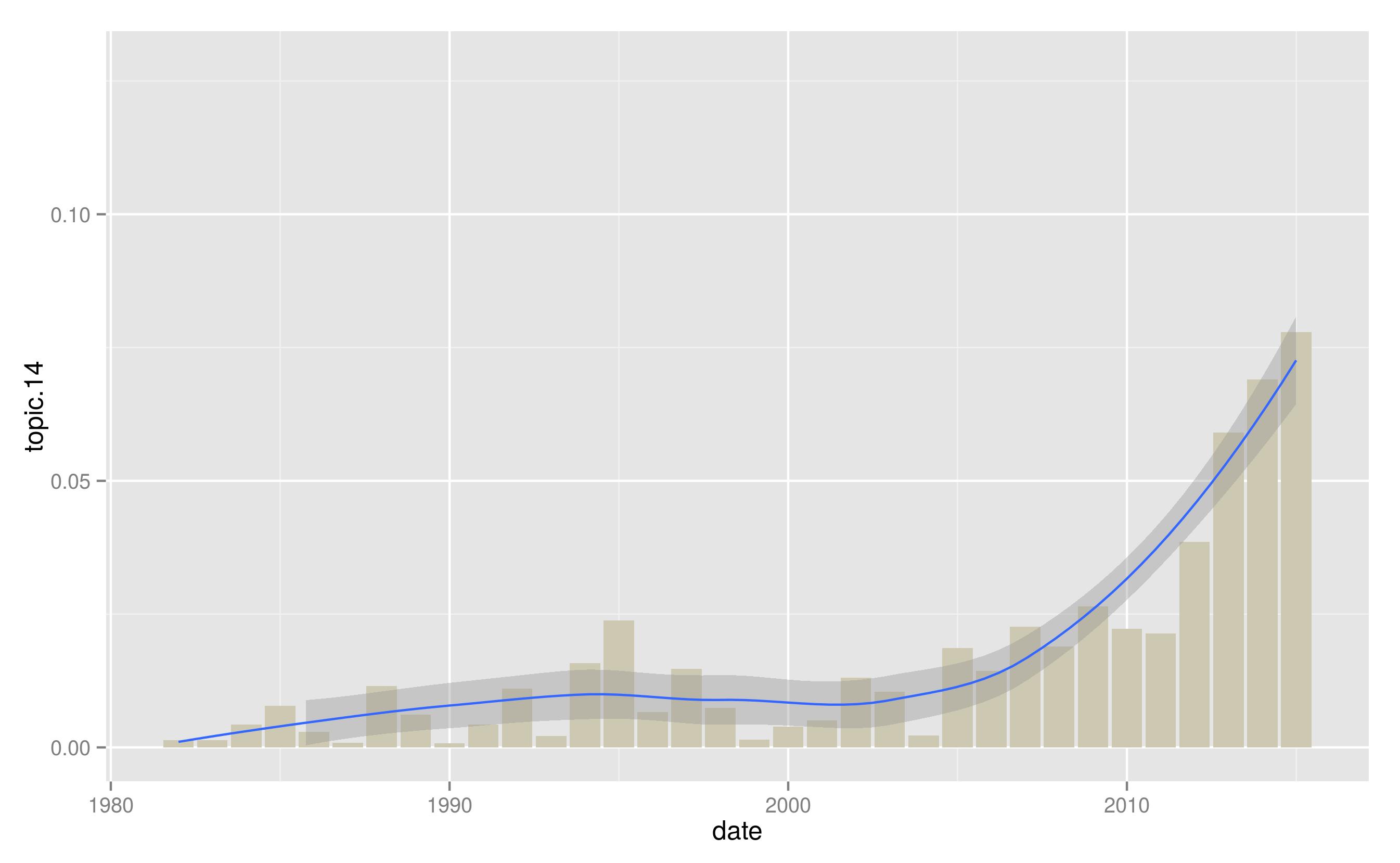,width=\linewidth} 
\end{center}
\label{fig:interest7}
\end{subfigure}
\caption{Topic: Telecommunication applications.}
\end{figure}

\begin{figure}[t]
\begin{subfigure}[t]{0.41\linewidth}
\begin{center}
\epsfig{file=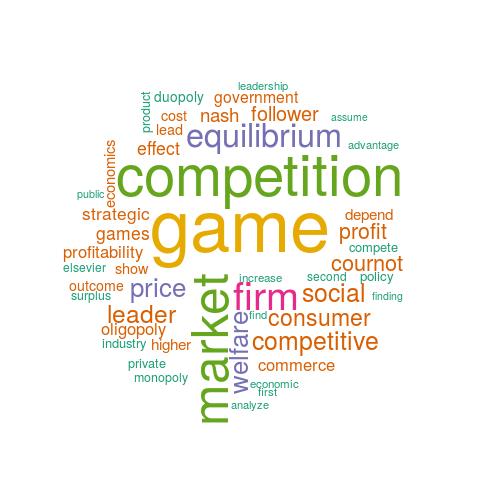,width=\linewidth}
\end{center}
\label{fig:cloud8}
\end{subfigure}\hfill
\begin{subfigure}[t]{0.59\linewidth}
\begin{center}
\epsfig{file=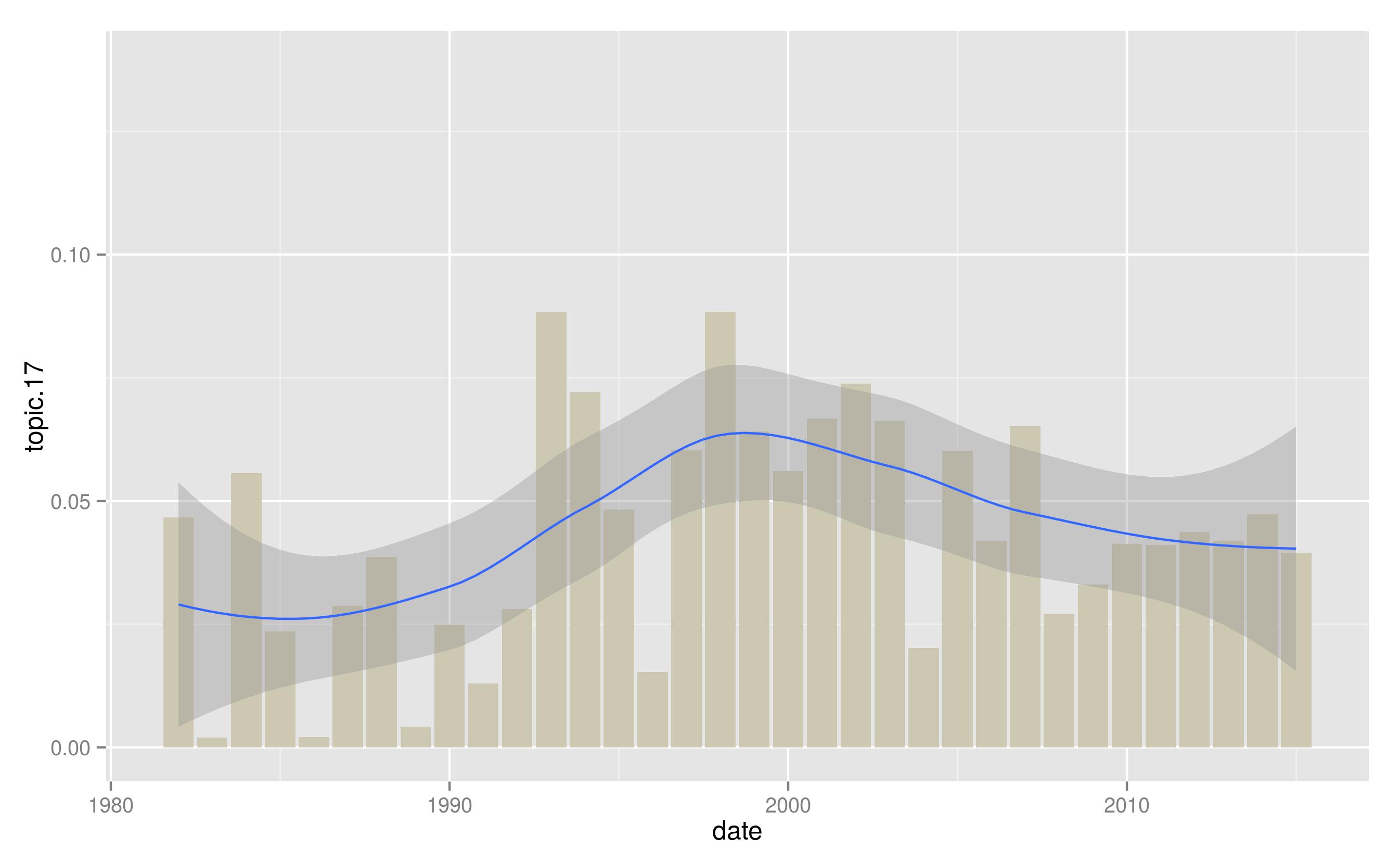,width=\linewidth} 
\end{center}
\label{fig:interest8}
\end{subfigure}
\caption{Topic: Business applications.}

\begin{subfigure}[t]{0.41\linewidth}
\begin{center}
\epsfig{file=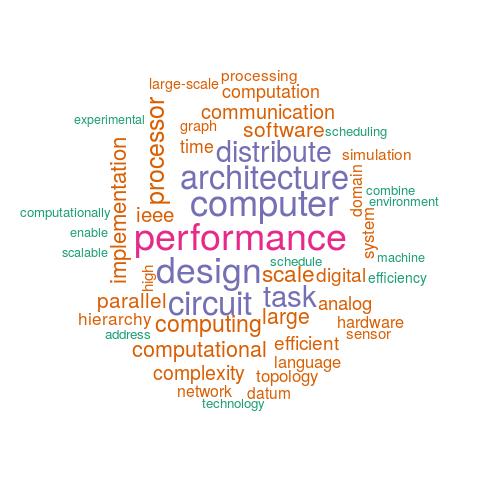,width=\linewidth}
\end{center}
\label{fig:cloud9}
\end{subfigure}\hfill
\begin{subfigure}[t]{0.59\linewidth}
\begin{center}
\epsfig{file=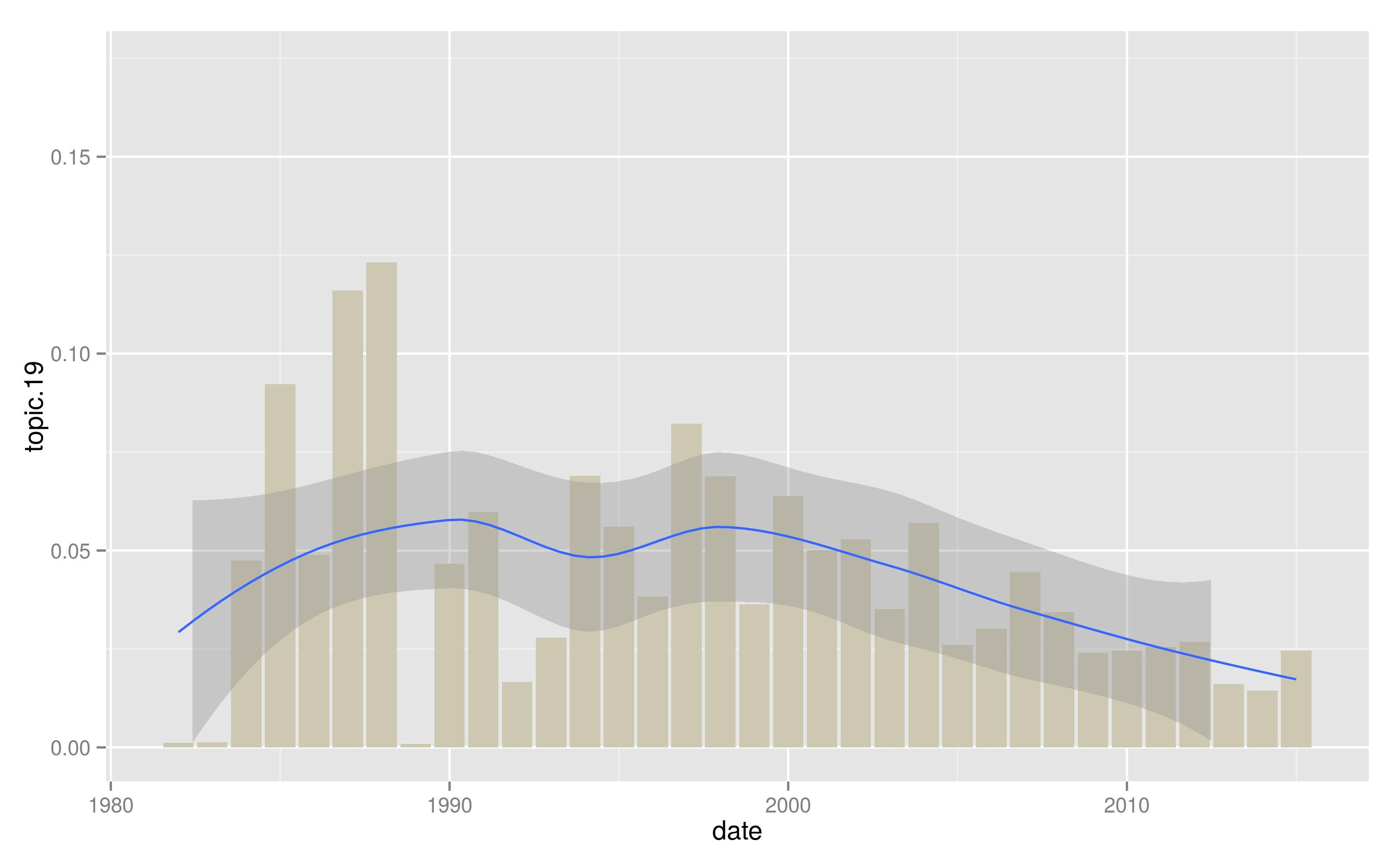,width=\linewidth} 
\end{center}
\label{fig:interest9}
\end{subfigure}
\caption{Topic: Computer architecture and circuit design.}

\begin{subfigure}[t]{0.41\linewidth}
\begin{center}
\epsfig{file=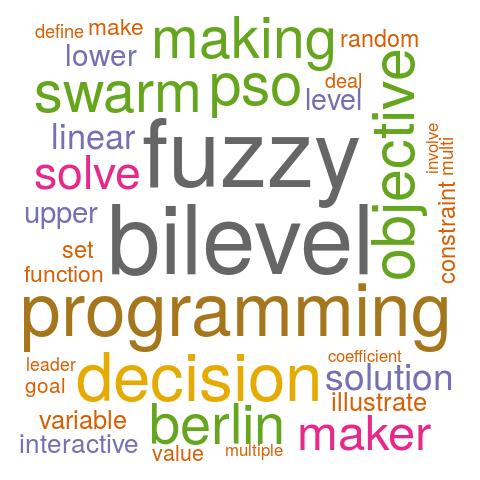,width=\linewidth}
\end{center}
\label{fig:cloud10}
\end{subfigure}\hfill
\begin{subfigure}[t]{0.59\linewidth}
\begin{center}
\epsfig{file=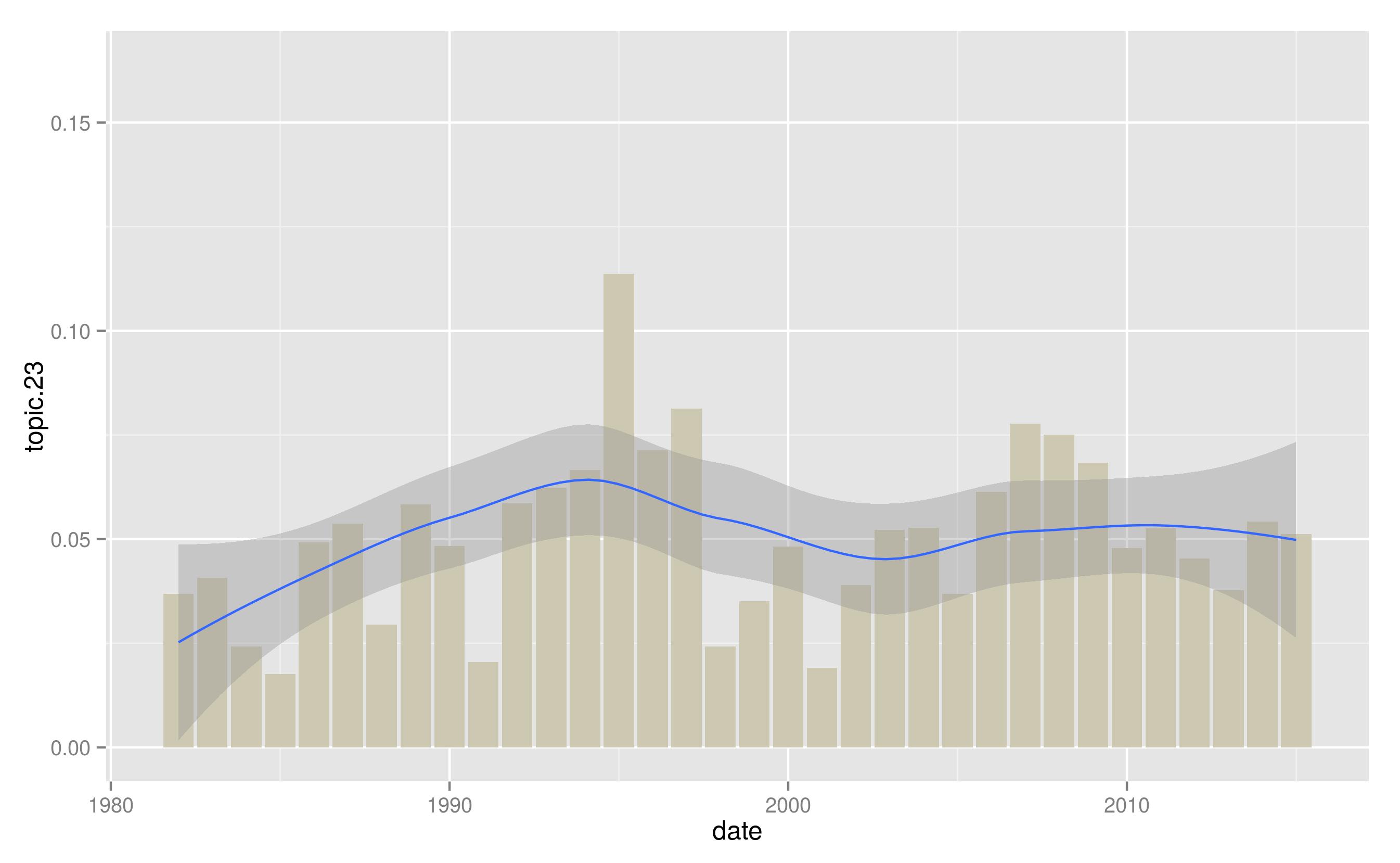,width=\linewidth} 
\end{center}
\label{fig:interest10}
\end{subfigure}
\caption{Topic: Hierarchical decision making applications.}

\begin{subfigure}[t]{0.41\linewidth}
\begin{center}
\epsfig{file=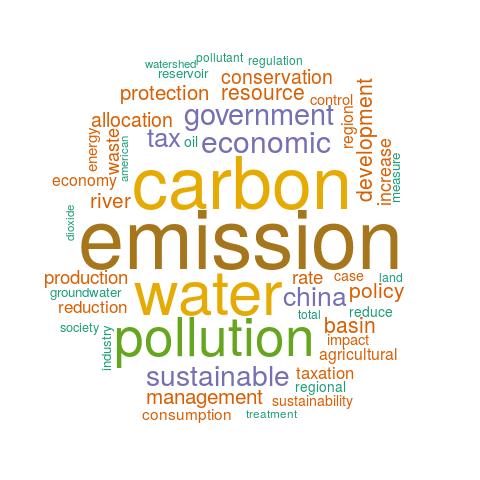,width=\linewidth}
\end{center}
\label{fig:cloud11}
\end{subfigure}\hfill
\begin{subfigure}[t]{0.59\linewidth}
\begin{center}
\epsfig{file=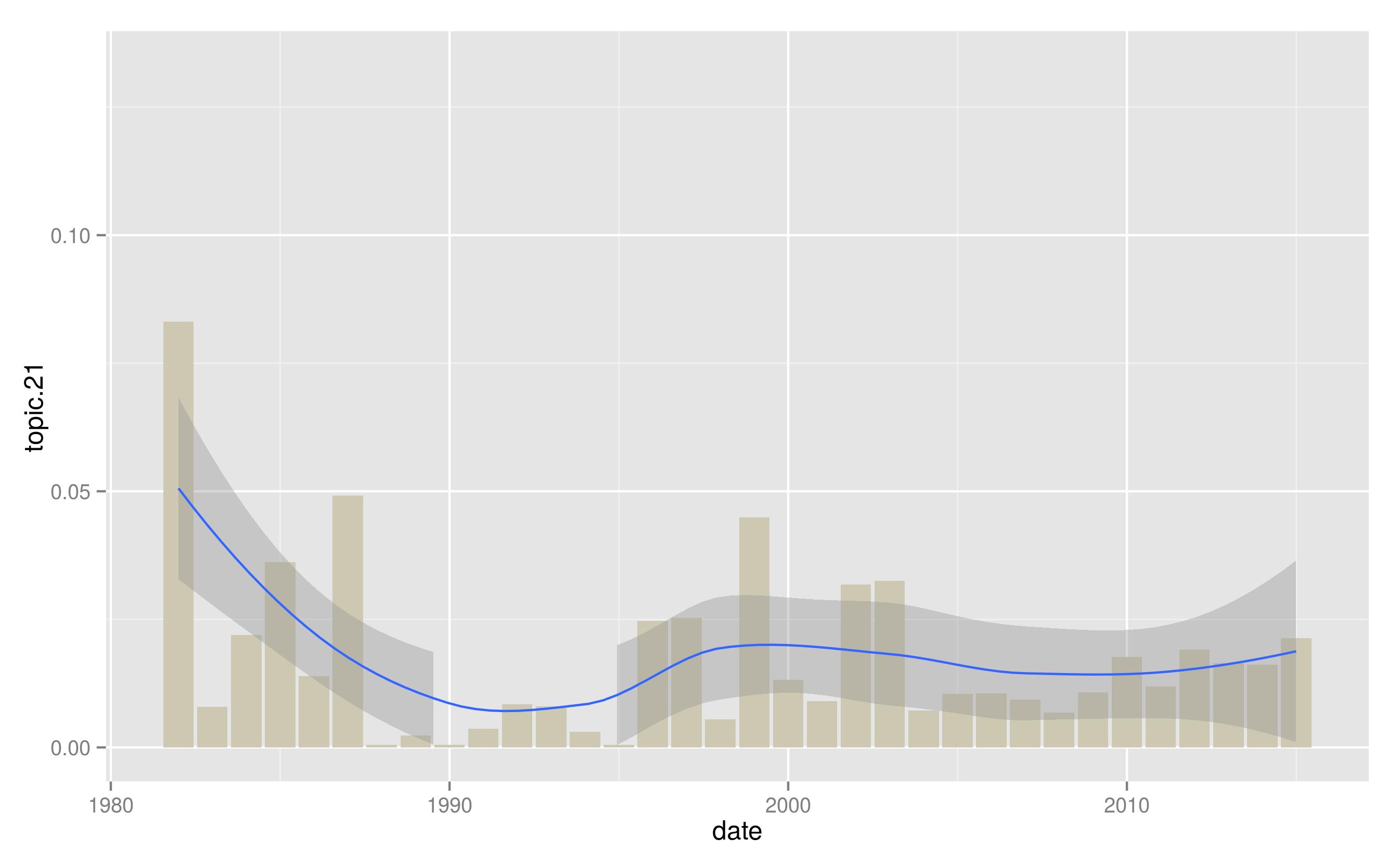,width=\linewidth} 
\end{center}
\label{fig:interest11}
\end{subfigure}
\caption{Topic: Environment applications.}
\end{figure}

\begin{figure}[t]
\begin{subfigure}[t]{0.41\linewidth}
\begin{center}
\epsfig{file=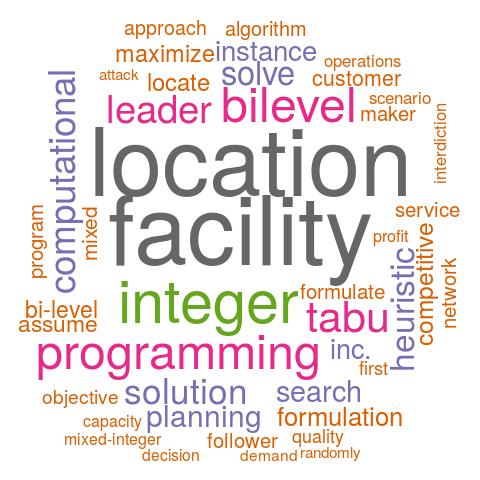,width=\linewidth}
\end{center}
\label{fig:cloud12}
\end{subfigure}\hfill
\begin{subfigure}[t]{0.59\linewidth}
\begin{center}
\epsfig{file=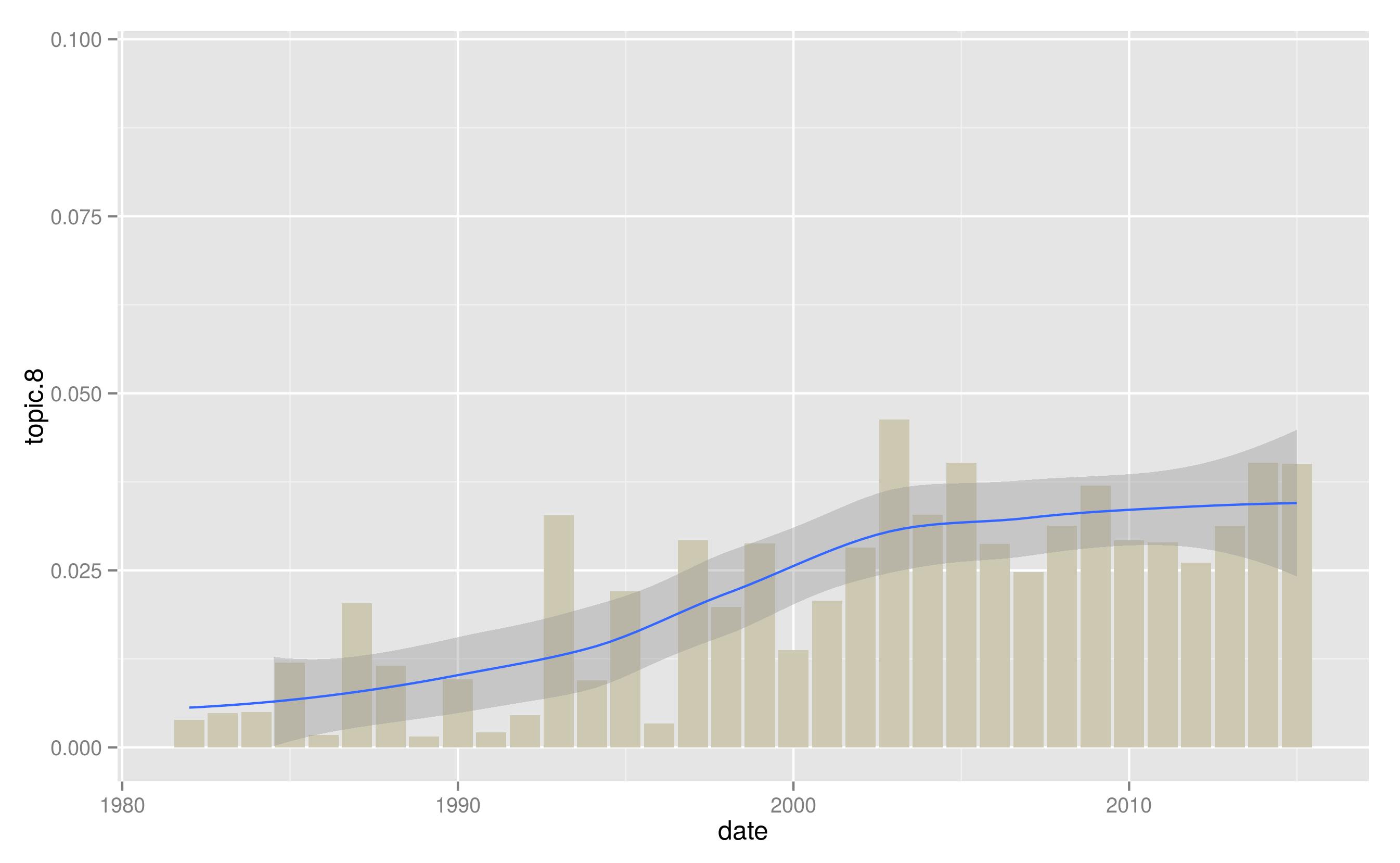,width=\linewidth} 
\end{center}
\label{fig:interest12}
\end{subfigure}
\caption{Topic: Facility location applications.}

\begin{subfigure}[t]{0.41\linewidth}
\begin{center}
\epsfig{file=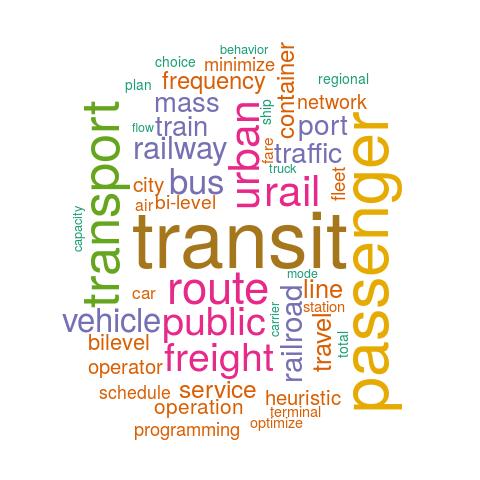,width=\linewidth}
\end{center}
\label{fig:cloud13}
\end{subfigure}\hfill
\begin{subfigure}[t]{0.59\linewidth}
\begin{center}
\epsfig{file=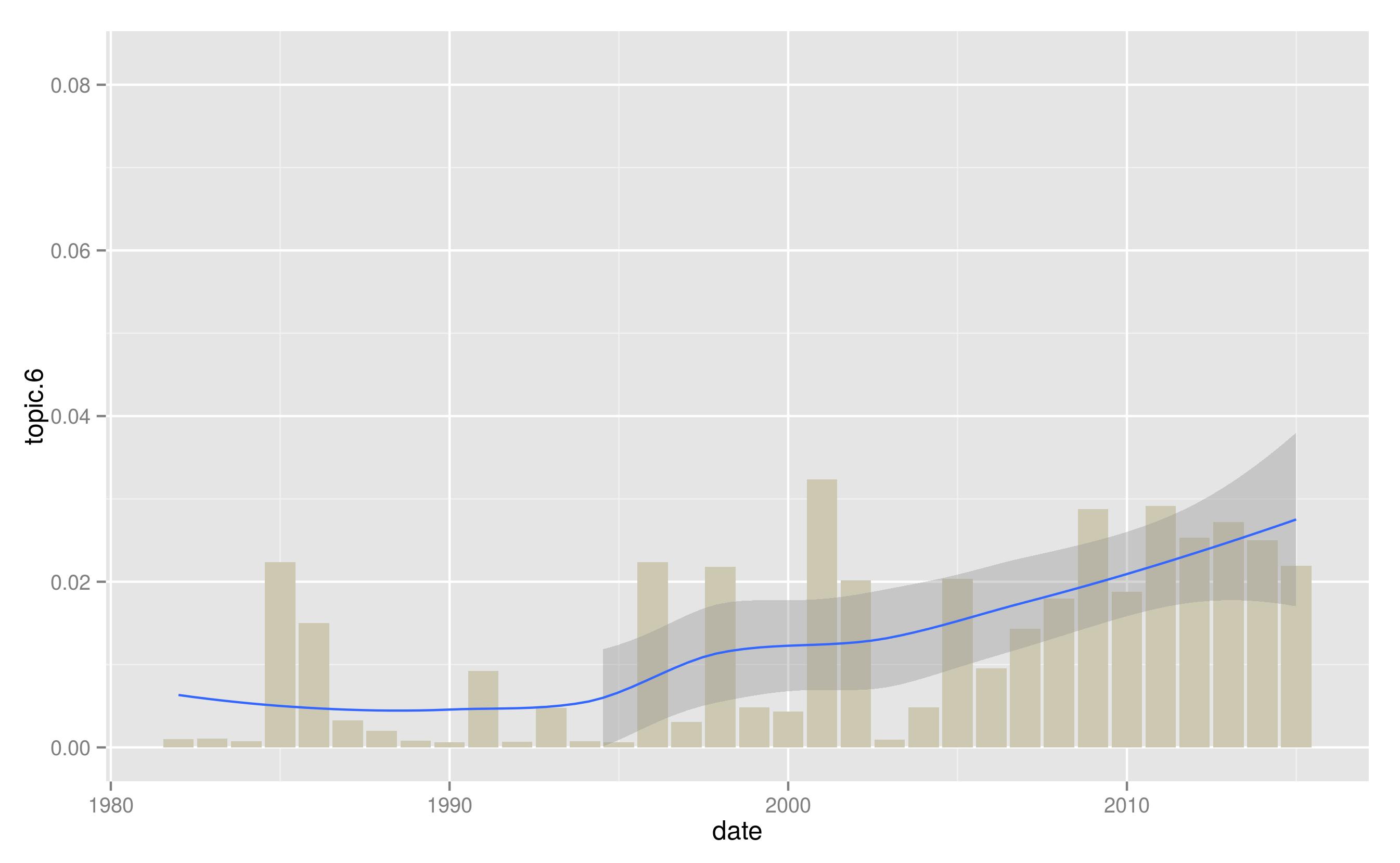,width=\linewidth} 
\end{center}
\label{fig:interest13}
\end{subfigure}
\caption{Topic: Railway applications.}

\begin{subfigure}[t]{0.41\linewidth}
\begin{center}
\epsfig{file=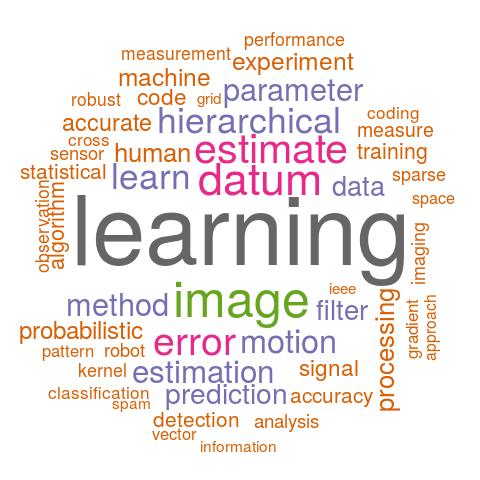,width=\linewidth}
\end{center}
\label{fig:cloud14}
\end{subfigure}\hfill
\begin{subfigure}[t]{0.59\linewidth}
\begin{center}
\epsfig{file=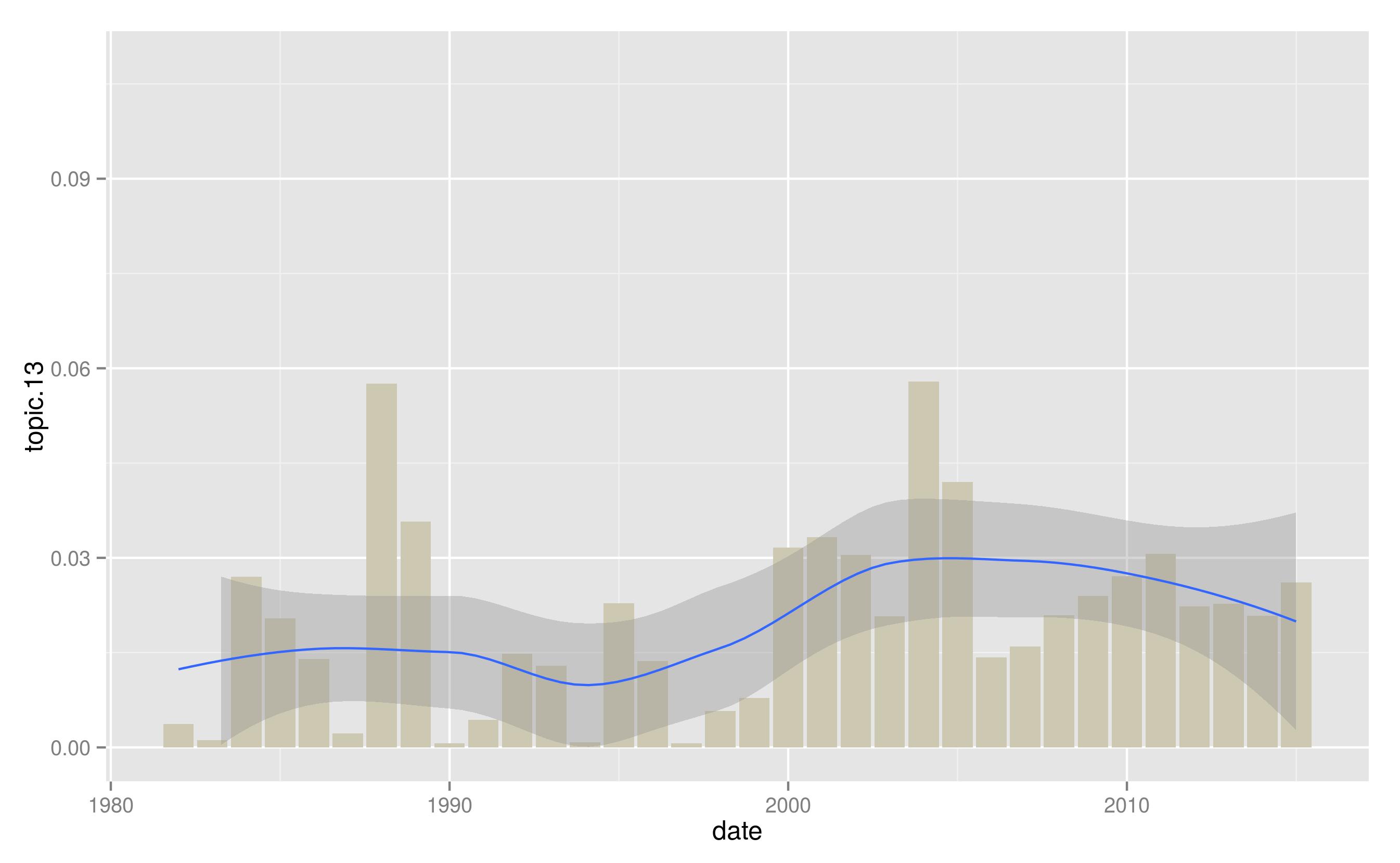,width=\linewidth} 
\end{center}
\label{fig:interest14}
\end{subfigure}
\caption{Topic: Machine learning applications.}

\begin{subfigure}[t]{0.41\linewidth}
\begin{center}
\epsfig{file=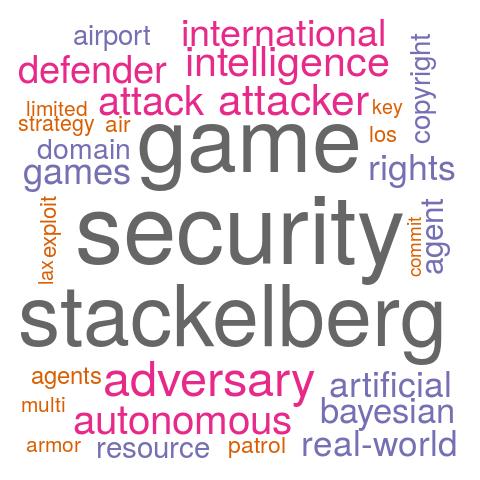,width=\linewidth}
\end{center}
\label{fig:cloud15}
\end{subfigure}\hfill
\begin{subfigure}[t]{0.59\linewidth}
\begin{center}
\epsfig{file=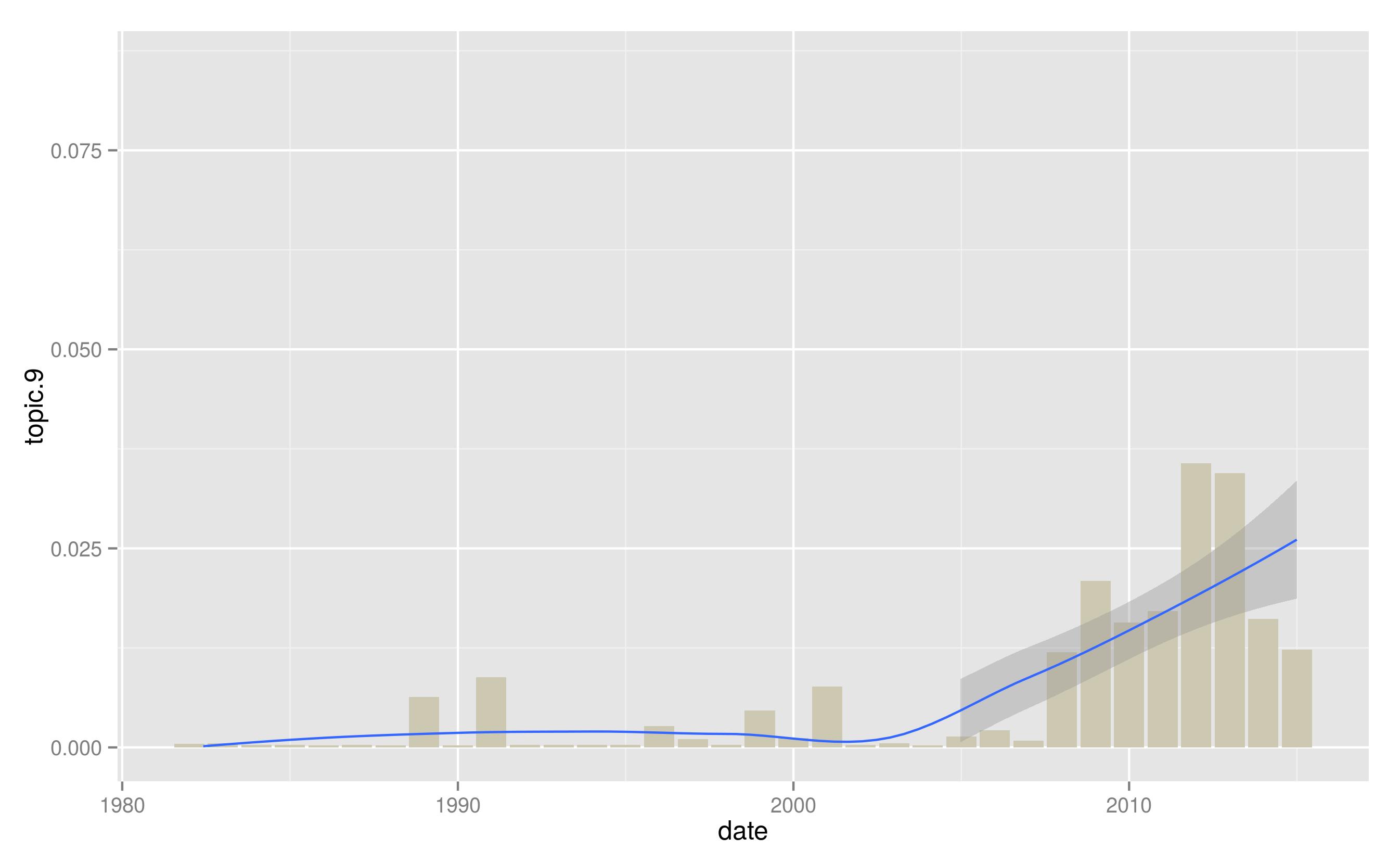,width=\linewidth} 
\end{center}
\label{fig:interest15}
\end{subfigure}
\caption{Topic: Defense applications.}
\label{fig:last}
\end{figure}

To identify the themes that have contributed towards the literature on bilevel optimization we used topic models. ``Topic models are algorithms for discovering the main themes that pervade a large and otherwise unstructured collection of documents.'' \cite{blei2012probabilistic} These models can be used to organize unstructured collections as well as develop insights from large text databases that made it suitable for our purposes. The results from the topic model for the documents retrieved from SCOPUS are given in Figures \ref{fig:first}-\ref{fig:last}. The figures consist of themes that have a volume at least to the order of around $1\%$. Along with the identification of themes, the analysis helped in determining the attention received by particular themes at any point in time. 

Each figure contains a word cloud in the inset that describes the theme and volume of papers as the other inset that describes the number of papers published on that theme over the years. Interestingly, we observe that number of papers on classical bilevel methods and optimality conditions peaked during 1995-2000. Since 2000 a number of bilevel applications picked-up, for example, we see a growth in supply chain applications, electricity transmission applications, telecommunication applications, facility location applications, railway applications and machine learning applications. Defense applications that appeared to be minimally present before 2000 show a significant presence later. Network design, optimal design and business applications do not show a trend but represent the highest volume of applications on bilevel optimization.




\section{Conclusions and Future Research Directions}

In this concluding section, we will raise a few perspectives that have not yet received much attention but may offer interesting directions for future research. Apart from metamodeling based techniques to solve bilevel problems, we would like to highlight the importance of being able to account for different forms of uncertainties that are often encountered when solving practical problems. Another interesting direction is concerned with scalability of bilevel algorithms and ability to leverage distributed computing platforms to handle large scale problems.

{\color{black}
\subsection{Metamodeling-based Algorithms}
Though we have already highlighted the importance of metamodeling-based methods in an earlier section, we have decided to discuss it once again because of its potential in solving practical bilevel optimization problems. Bilevel optimization problems inherit a number of mappings, and any knowledge about the structure of these mappings can simplify the solution procedure extensively. In this paper, we have highlighted approaches that are based on approximating the reaction set mapping and the optimal lower level value function mapping. Knowing one of these mappings reduces any bilevel optimization problem to a single-level optimization problem. For solving large scale bilevel problems, one has to exploit the structure and properties of bilevel problems that are essentially contained in these mappings. Other ways of utilizing metamodeling for bilevel problems that we have discussed are: approximating the bilevel problem by bypassing the lower level problem completely; and utilizing auxiliary bilevel models to locally approximate a bilevel optimization problem. Only few studies in the context of evolutionary algorithms utilize such a strategies and offer opportunities for future contributions.

\subsection{Multiobjective Bilevel Optimization and Decision Making}
Multiobjective bilevel optimization has received only lukewarm interest from researchers. A number of issues, like decision interaction between the two levels and uncertainties in decision making, remain unexplored. It might be of interest for researchers working in the area of multi-criteria decision making and multiobjective evolutionary optimization to explore how two levels of decision makers interact to arrive at a compromising or an equilibrium solution in different situations. Similarly, the notion of uncertain decision-maker's preferences at one or both levels has also not received enough attention~\cite{my-ieeetec16}. In the field of Decision Analysis, however, plenty of research has been carried out to extend traditional frameworks of decision making such as Expected Utility Theory~\cite{neumann-morgenstern} and Multi-attribute Utility Theory~\cite{keeney-raiffa} to account for uncertain preferences on, for instance, the trade-offs among multiple decision objectives or the risk-attitude. However, preferential uncertainty in bilevel optimization problems still requires development of theory as well as methods to account for decision behavior in a hierarchical setting. The approaches that have been proposed so far are still very preliminary and require substantial future research.
}

\subsection{Bilevel Optimization under Variable Uncertainty}
Another important research topic is concerned with the inherent uncertainty of decision variables. This poses several challenges for the existing more deterministic optimization frameworks that may fail to find solutions that are robust and sufficiently close to the optimal solutions. The fact that bilevel problems have nested optimization tasks makes the search of robust solutions substantially more challenging compared with single-level optimization problems. A few preliminary ideas for handling variable uncertainty have already been suggested~\cite{my-cec15d}, but the required algorithm side innovations that would make these problems accessible to practitioners are still missing.

\subsection{Scaling of Evolutionary Bilevel Algorithms}
Bilevel problems are well-known for being highly computationally intensive already before considering any types of uncertainties. One of the promising directions for handling larger bilevel problems could be the use of the recent distributed computing platforms such as Apache Spark project. The programming model of Spark is quite different from the Hadoop MapReduce framework, and it has managed to overcome many of the earlier limitations. In particular, its current form has turned out to support the use and development of iterative algorithms quite well. Therefore, it may be interesting to investigate whether this novel platform will be able to offer a helping hand for researchers and practitioners who seek to solve bigger bilevel problems faster.

Though considerable progress has already been made during the last few years, evolutionary bilevel optimization is still a relatively young field with numerous opportunities for both computational as well as theoretical innovations. The growing availability of algorithms is also opening the field to more applied research, and we believe that in the future we are likely to see a considerable amount of novel applications. 

\section{Acknowledgments}
Ankur Sinha and Pekka Malo would like to acknowledge the support provided by Liikesivistysrahasto and Helsinki School of Economics Foundation. K. Deb acknowledges the support
from NSF Beacon Center for the study of evolution in action at
MSU under Cooperative Agreement No. DBI-0939454. Any
opinions, findings, and conclusions or recommendations
expressed in this material are those of the author(s) and do not
necessarily reflect the views of the National Science Foundation.


\begin{thebibliography}{100}

\bibitem{AiSh81}
E.~Aiyoshi and K.~Shimizu.
\newblock Hierarchical decentralized systems and its new solution by a barrier
  method.
\newblock {\em {IEEE} Transactions on Systems, Man, and Cybernetics},
  6:444--449, 1981.

\bibitem{AiSh84}
E.~Aiyoshi and K.~Shimizu.
\newblock {A solution method for the static constrained Stackelberg problem via
  penalty method}.
\newblock {\em {IEEE} Transactions on Automatic Control}, 29:1111--1114, 1984.

\bibitem{aksen2012bilevel}
Deniz Aksen and Necati Aras.
\newblock A bilevel fixed charge location model for facilities under imminent
  attack.
\newblock {\em Computers \& Operations Research}, 39(7):1364--1381, 2012.

\bibitem{aksen2013matheuristic}
Deniz Aksen and Necati Aras.
\newblock A matheuristic for leader-follower games involving facility
  location-protection-interdiction decisions.
\newblock In {\em Metaheuristics for Bi-level Optimization}, pages 115--151.
  Springer, 2013.

\bibitem{AlHoPa92}
F.~Al-Khayyal, R.~Horst, and P.~Pardalos.
\newblock Global optimization of concave functions subject to quadratic
  constraints: an application in nonlinear bilevel programming.
\newblock {\em Annals of Operations Research}, 34:125--147, 1992.

\bibitem{albrecht2011imitating}
Sebastian Albrecht, K~Ramirez-Amaro, Federico Ruiz-Ugalde, David Weikersdorfer,
  M~Leibold, Michael Ulbrich, and Michael Beetz.
\newblock Imitating human reaching motions using physically inspired
  optimization principles.
\newblock In {\em Humanoid Robots (Humanoids), 2011 11th IEEE-RAS International
  Conference on}, pages 602--607. IEEE, 2011.

\bibitem{alekseeva2009hybrid}
E.~Alekseeva, N.~Kochetova, Yu~Kochetov, and A.~Plyasunov.
\newblock A hybrid memetic algorithm for the competitive p-median problem.
\newblock {\em IFAC Proceedings Volumes}, 42(4):1533--1537, 2009.

\bibitem{amouzegar1999determining}
Mahyar~A. Amouzegar and Khosrow Moshirvaziri.
\newblock Determining optimal pollution control policies: An application of
  bilevel programming.
\newblock {\em European Journal of Operational Research}, 119(1):100--120,
  1999.

\bibitem{an13}
B.~An, F.~Ord\'{o}\~{n}ez, M.~Tambe, E.~Shieh, R.~Yang, C.~Baldwin,
  J.~DiRenzo~III, K.~Moretti, B.~Maule, and G.~Meyer.
\newblock {A Deployed Quantal Response-Based Patrol Planning System for the
  U.S. Coast Guard}.
\newblock {\em Interfaces}, 43(5):400--420, 2013.

\bibitem{AnFr92}
G.~Anandalingam and T.~Friesz.
\newblock Hierarchical optimization: an introduction.
\newblock {\em Annals of Operations Research}, 34:1--11, 1992.

\bibitem{angelo13}
J.~Angelo, E.~Krempser, and H.~Barbosa.
\newblock Differential evolution for bilevel programming.
\newblock In {\em Proceedings of the 2013 Congress on Evolutionary Computation
  (CEC-2013)}. IEEE Press, 2013.

\bibitem{angelo2015study}
Jaqueline~S. Angelo and Helio J.~C. Barbosa.
\newblock A study on the use of heuristics to solve a bilevel programming
  problem.
\newblock {\em International Transactions in Operational Research}, 2015.

\bibitem{angelo2014differential}
Jaqueline~S Angelo, Eduardo Krempser, and Helio~JC Barbosa.
\newblock Differential evolution assisted by a surrogate model for bilevel
  programming problems.
\newblock In {\em Evolutionary Computation (CEC), 2014 IEEE Congress on}, pages
  1784--1791. IEEE, 2014.

\bibitem{arroyo2009genetic}
J.~M. Arroyo and F.~J. Fern{\'a}ndez.
\newblock A genetic algorithm approach for the analysis of electric grid
  interdiction with line switching.
\newblock In {\em Intelligent System Applications to Power Systems, 2009.
  ISAP'09. 15th International Conference on}, pages 1--6. IEEE, 2009.

\bibitem{bank1983non}
Bernd Bank, J{\"u}rgen Guddat, Diethard Klatte, Bernd Kummer, and Klaus Tammer.
\newblock {\em Non-linear parametric optimization}.
\newblock Birkh{\"a}user Basel, 1983.

\bibitem{BaFa82}
J.~Bard and J.~Falk.
\newblock An explicit solution to the multi-level programming problem.
\newblock {\em Computers and Operations Research}, 9:77--100, 1982.

\bibitem{BaMo90}
J.~Bard and J.~Moore.
\newblock A branch and bound algorithm for the bilevel programming problem.
\newblock {\em {SIAM} Journal on Scientific and Statistical Computing},
  11:281--292, 1990.

\bibitem{bilevel-book}
J.~F. Bard.
\newblock {\em Practical Bilevel Optimization: Algorithms and Applications}.
\newblock The Netherlands: Kluwer, 1998.

\bibitem{bard1992algorithm}
Jonathan~F. Bard and James~T. Moore.
\newblock An algorithm for the discrete bilevel programming problem.
\newblock {\em Naval Research Logistics (NRL)}, 39(3):419--435, 1992.

\bibitem{Be93}
O.~Ben-Ayed.
\newblock Bilevel linear programming.
\newblock {\em Computers and Operations Research}, 20:485--501, 1993.

\bibitem{BeBlBoLe92}
O.~Ben-Ayed, C.~Blair, D.~Boyce, and L.~LeBlanc.
\newblock Construction of a real-world bilevel linear programming model of the
  highway design problem.
\newblock {\em Annals of Operations Research}, 34:219--254, 1992.

\bibitem{bendsoe95}
Martin~P Bendsoe.
\newblock {\em Optimization of structural topology, shape, and material},
  volume 414.
\newblock Springer, 1995.

\bibitem{bennett2006model}
Kristin~P. Bennett, Jing Hu, Xiaoyun Ji, Gautam Kunapuli, and Jong-Shi Pang.
\newblock Model selection via bilevel optimization.
\newblock In {\em Neural Networks, 2006. IJCNN'06. International Joint
  Conference on}, pages 1922--1929. IEEE, 2006.

\bibitem{bennett2008bilevel}
Kristin~P. Bennett, Gautam Kunapuli, Jing Hu, and Jong-Shi Pang.
\newblock Bilevel optimization and machine learning.
\newblock In {\em Computational Intelligence: Research Frontiers}, pages
  25--47. Springer, 2008.

\bibitem{BiKa84}
W.~Bialas and M.~Karwan.
\newblock Two-level linear programming.
\newblock {\em Management Science}, 30:1004--1020, 1984.

\bibitem{bialas1984two}
Wayne~F. Bialas and Mark~H. Karwan.
\newblock Two-level linear programming.
\newblock {\em Management science}, 30(8):1004--1020, 1984.

\bibitem{blei2012probabilistic}
David~M. Blei.
\newblock Probabilistic topic models.
\newblock {\em Communications of the ACM}, 55(4):77--84, 2012.

\bibitem{my-cec15c}
M.~Bostian, A.~Sinha, Gerald Whittaker, and Bradley Barnhart.
\newblock Incorporating data envelopment analysis solution methods into bilevel
  multi-objective optimization.
\newblock In {\em 2015 IEEE Congress on Evolutionary Computation (CEC-2015)}.
  IEEE Press, 2015.

\bibitem{bostian2015valuing}
Moriah Bostian, Gerald Whittaker, Brad Barnhart, Rolf F{\"a}re, and Shawna
  Grosskopf.
\newblock Valuing water quality tradeoffs at different spatial scales: An
  integrated approach using bilevel optimization.
\newblock {\em Water Resources and Economics}, 11:1--12, 2015.

\bibitem{BrMc73}
J.~Bracken and J.~McGill.
\newblock Mathematical programs with optimization problems in the constraints.
\newblock {\em Operations Research}, 21:37--44, 1973.

\bibitem{bracken1974defense}
Jerome Bracken and James~T. McGill.
\newblock Defense applications of mathematical programs with optimization
  problems in the constraints.
\newblock {\em Operations Research}, 22(5):1086--1096, 1974.

\bibitem{brotcorne01}
Luce Brotcorne, Martine Labbe, Patrice Marcotte, and Gilles Savard.
\newblock A bilevel model for toll optimization on a multicommodity
  transportation network.
\newblock {\em Transportation Science}, 35(4):345--358, 2001.

\bibitem{brown05}
G.~Brown, M.~Carlyle, D.~Diehl, J.~Kline, and K.~Wood.
\newblock {A Two-Sided Optimization for Theater Ballistic Missile Defense}.
\newblock {\em Operations Research}, 53(5):745--763, 2005.

\bibitem{brown09}
G.~Brown, M.~Carlyle, R.~C. Harney, E.M. Skroch, and K.~Wood.
\newblock {Interdicting a Nuclear-Weapons Project}.
\newblock {\em Operations Research}, 57(4):866--877, 2009.

\bibitem{brown2006defending}
Gerald Brown, Matthew Carlyle, Javier Salmer{\'o}n, and Kevin Wood.
\newblock Defending critical infrastructure.
\newblock {\em Interfaces}, 36(6):530--544, 2006.

\bibitem{calvete2013efficient}
Herminia~I. Calvete, Carmen Gal{\'e}, and Jos{\'e}~A. Iranzo.
\newblock An efficient evolutionary algorithm for the ring star problem.
\newblock {\em European Journal of Operational Research}, 231(1):22--33, 2013.

\bibitem{calvete2011bilevel}
Herminia~I. Calvete, Carmen Gal{\'e}, and Mar{\'\i}a-Jos{\'e} Oliveros.
\newblock Bilevel model for production--distribution planning solved by using
  ant colony optimization.
\newblock {\em Computers \& Operations Research}, 38(1):320--327, 2011.

\bibitem{camacho2014solving}
Jos{\'e}-Fernando Camacho-Vallejo, {\'A}lvaro~Eduardo Cordero-Franco, and
  Rosa~G. Gonz{\'a}lez-Ram{\'\i}rez.
\newblock Solving the bilevel facility location problem under preferences by a
  stackelberg-evolutionary algorithm.
\newblock {\em Mathematical Problems in Engineering}, 2014, 2014.

\bibitem{camacho2015genetic}
Jos{\'e}-Fernando Camacho-Vallejo, Julio Mar-Ortiz, Francisco L{\'o}pez-Ramos,
  and Ricardo~Pedraza Rodr{\'\i}guez.
\newblock A genetic algorithm for the bi-level topological design of local area
  networks.
\newblock {\em PloS one}, 10(6):e0128067, 2015.

\bibitem{camacho2015heuristic}
Jos{\'e}-Fernando Camacho-Vallejo, Rafael Mu{\~n}oz-S{\'a}nchez, and
  Jos{\'e}~Luis Gonz{\'a}lez-Velarde.
\newblock A heuristic algorithm for a supply chain's production-distribution
  planning.
\newblock {\em Computers \& Operations Research}, 61:110--121, 2015.

\bibitem{caramia2015decomposition}
Massimiliano Caramia and Renato Mari.
\newblock A decomposition approach to solve a bilevel capacitated facility
  location problem with equity constraints.
\newblock {\em Optimization Letters}, pages 1--23, 2015.

\bibitem{cecchini2013solving}
Mark Cecchini, Joseph Ecker, Michael Kupferschmid, and Robert Leitch.
\newblock Solving nonlinear principal-agent problems using bilevel programming.
\newblock {\em European Journal of Operational Research}, 230(2):364--373,
  2013.

\bibitem{ceylan2004traffic}
Halim Ceylan and Michael G.~H. Bell.
\newblock Traffic signal timing optimisation based on genetic algorithm
  approach, including drivers\d5 routing.
\newblock {\em Transportation Research Part B: Methodological}, 38(4):329--342,
  2004.

\bibitem{chaabani2015co}
Abir Chaabani, Slim Bechikh, and Lamjed~Ben Said.
\newblock A co-evolutionary decomposition-based algorithm for bi-level
  combinatorial optimization.
\newblock In {\em 2015 IEEE Congress on Evolutionary Computation (CEC)}, pages
  1659--1666. IEEE, 2015.

\bibitem{chen2010stochastic}
Anthony Chen, Juyoung Kim, Seungjae Lee, and Youngchan Kim.
\newblock Stochastic multi-objective models for network design problem.
\newblock {\em Expert Systems with Applications}, 37(2):1608--1619, 2010.

\bibitem{ChFl92}
Y.~Chen and M.~Florian.
\newblock On the geometry structure of linear bilevel programs: a dual
  approach.
\newblock Technical Report CRT-867, Centre de Recherche sur les Transports,
  1992.

\bibitem{christiansen2001stochastic}
Snorre Christiansen, Michael Patriksson, and Laura Wynter.
\newblock Stochastic bilevel programming in structural optimization.
\newblock {\em Structural and multidisciplinary optimization}, 21(5):361--371,
  2001.

\bibitem{clark1990bilevel}
Peter~A. Clark and Arthur~W. Westerberg.
\newblock Bilevel programming for steady-state chemical process design-i.
  fundamentals and algorithms.
\newblock {\em Computers \& Chemical Engineering}, 14(1):87--97, 1990.

\bibitem{Co99}
B.~Colson.
\newblock Mathematical programs with equilibrium constraints and nonlinear
  bilevel programming problems.
\newblock Technical report, Master's thesis, Department of Mathematics, FUNDP,
  Namur, Belgium, 1999.

\bibitem{colson}
B.~Colson, P.~Marcotte, and G.~Savard.
\newblock An overview of bilevel optimization.
\newblock {\em Annals of Operational Research}, 153:235--256, 2007.

\bibitem{colson2005trust}
Beno{\^\i}t Colson, Patrice Marcotte, and Gilles Savard.
\newblock A trust-region method for nonlinear bilevel programming: algorithm
  and computational experience.
\newblock {\em Computational Optimization and Applications}, 30(3):211--227,
  2005.

\bibitem{constantin95}
Isabelle Constantin and Michael Florian.
\newblock Optimizing frequencies in a transit network: a nonlinear bi-level
  programming approach.
\newblock {\em International Transactions in Operational Research}, 2(2):149 --
  164, 1995.

\bibitem{my-cec15d}
K.~Deb, Z.~Lu, and A.~Sinha.
\newblock Handling decision variable uncertainty in bilevel optimization
  problems.
\newblock In {\em 2015 IEEE Congress on Evolutionary Computation (CEC-2015)}.
  IEEE Press, 2015.

\bibitem{my-cec09a}
K.~Deb and A.~Sinha.
\newblock Constructing test problems for bilevel evolutionary multi-objective
  optimization.
\newblock In {\em 2009 IEEE Congress on Evolutionary Computation (CEC-2009)},
  pages 1153--1160. IEEE Press, 2009.

\bibitem{my-mcdm09}
K.~Deb and A.~Sinha.
\newblock An evolutionary approach for bilevel multi-objective problems.
\newblock In {\em Cutting-Edge Research Topics on Multiple Criteria Decision
  Making, Communications in Computer and Information Science}, volume~35, pages
  17--24. Berlin, Germany: Springer, 2009.

\bibitem{my-emo09}
K.~Deb and A.~Sinha.
\newblock Solving bilevel multi-objective optimization problems using
  evolutionary algorithms.
\newblock In {\em Evolutionary Multi-Criterion Optimization (EMO-2009)}, pages
  110--124. Berlin, Germany: Springer-Verlag, 2009.

\bibitem{my-ecj10}
K.~Deb and A.~Sinha.
\newblock An efficient and accurate solution methodology for bilevel
  multi-objective programming problems using a hybrid evolutionary-local-search
  algorithm.
\newblock {\em Evolutionary Computation Journal}, 18(3):403--449, 2010.

\bibitem{dempe2003}
S.~Dempe.
\newblock Annotated bibliography on bilevel programming and mathematical
  programs with equilibrium constraints.
\newblock {\em Optimization}, 52(3):339--359, 2003.

\bibitem{dempe2007new}
S.~Dempe, J.~Dutta, and B.S. Mordukhovich.
\newblock New necessary optimality conditions in optimistic bilevel
  programming.
\newblock {\em Optimization}, 56(5-6):577--604, 2007.

\bibitem{dempe1996discrete}
Stephan Dempe.
\newblock {\em Discrete bilevel optimization problems}.
\newblock Citeseer, 1996.

\bibitem{dempe02}
Stephan Dempe.
\newblock {\em {Foundations of Bilevel Programming}}.
\newblock Kluwer Academic Publishers, Secaucus, NJ, USA, 2002.

\bibitem{dempe2015bilevel}
Stephan Dempe, Vyacheslav Kalashnikov, Gerardo~A. P{\'e}rez-Vald{\'e}s, and
  Nataliya Kalashnykova.
\newblock Bilevel programming problems.
\newblock {\em Energy Systems. Springer, Berlin}, 2015.

\bibitem{dempe2014necessary}
Stephan Dempe, Boris~S. Mordukhovich, and Alain~Bertrand Zemkoho.
\newblock Necessary optimality conditions in pessimistic bilevel programming.
\newblock {\em Optimization}, 63(4):505--533, 2014.

\bibitem{deng1998complexity}
Xiaotie Deng.
\newblock Complexity issues in bilevel linear programming.
\newblock In {\em Multilevel optimization: Algorithms and applications}, pages
  149--164. Springer, 1998.

\bibitem{EdBa91}
T.~Edmunds and J.~Bard.
\newblock Algorithms for nonlinear bilevel mathematical programming.
\newblock {\em {IEEE} Transactions on Systems, Man, and Cybernetics},
  21:83--89, 1991.

\bibitem{eichfelder}
G.~Eichfelder.
\newblock Solving nonlinear multiobjective bilevel optimization problems with
  coupled upper level constraints.
\newblock Technical Report Preprint No. 320, Preprint-Series of the Institute
  of Applied Mathematics, Univ. Erlangen-Nornberg, Germany, 2007.

\bibitem{eichfelder2}
Gabriele Eichfelder.
\newblock Multiobjective bilevel optimization.
\newblock {\em Mathematical Programming}, 123(2):419--449, June 2010.

\bibitem{fan2015optimal}
Wei Fan.
\newblock Optimal congestion pricing toll design for revenue maximization:
  comprehensive numerical results and implications.
\newblock {\em Canadian Journal of Civil Engineering}, 42(8):544--551, 2015.

\bibitem{fan2011bi}
Wei Fan and Randy Machemehl.
\newblock Bi-level optimization model for public transportation network
  redesign problem: Accounting for equity issues.
\newblock {\em Transportation Research Record: Journal of the Transportation
  Research Board}, 2263:151--162, 2011.

\bibitem{fontaine2014benders}
Pirmin Fontaine and Stefan Minner.
\newblock Benders decomposition for discrete--continuous linear bilevel
  problems with application to traffic network design.
\newblock {\em Transportation Research Part B: Methodological}, 70:163--172,
  2014.

\bibitem{FoMc81}
J.~Fortuny-Amat and B.~McCarl.
\newblock A representation and economic interpretation of a two-level
  programming problem.
\newblock {\em Journal of the Operational Research Society}, 32:783--792, 1981.

\bibitem{my-ifac12}
A.~Frantsev, A.~Sinha, and P.~Malo.
\newblock Finding optimal strategies in multi-period stackelberg games using an
  evolutionary framework.
\newblock In {\em IFAC Workshop on Control Applications of Optimization
  (IFAC-2009)}. Elsevier, 2012.

\bibitem{gadhi2012necessary}
N~Gadhi and Stephan Dempe.
\newblock Necessary optimality conditions and a new approach to multiobjective
  bilevel optimization problems.
\newblock {\em Journal of Optimization Theory and Applications},
  155(1):100--114, 2012.

\bibitem{gallo2010meta}
Mariano Gallo, Luca D'Acierno, and Bruno Montella.
\newblock A meta-heuristic approach for solving the urban network design
  problem.
\newblock {\em European Journal of Operational Research}, 201(1):144--157,
  2010.

\bibitem{garen1994executive}
John~E. Garen.
\newblock Executive compensation and principal-agent theory.
\newblock {\em Journal of Political Economy}, pages 1175--1199, 1994.

\bibitem{gonzalez2015scatter}
Jos{\'e}~Luis Gonz{\'a}lez-Velarde, Jos{\'e}-Fernando Camacho-Vallejo, and
  Gabriel Pinto~Serrano.
\newblock A scatter search algorithm for solving a bilevel optimization model
  for determining highway tolls.
\newblock {\em Computaci{\'o}n y Sistemas}, 19(1):05--16, 2015.

\bibitem{halter-sanaz}
W.~Halter and S.~Mostaghim.
\newblock Bilevel optimization of multi-component chemical systems using
  particle swarm optimization.
\newblock In {\em Proceedings of World Congress on Computational Intelligence
  (WCCI-2006)}, pages 1240--1247, 2006.

\bibitem{halter2006bilevel}
Werner Halter and Sanaz Mostaghim.
\newblock Bilevel optimization of multi-component chemical systems using
  particle swarm optimization.
\newblock In {\em Evolutionary Computation, 2006. CEC 2006. IEEE Congress on},
  pages 1240--1247. IEEE, 2006.

\bibitem{handoko2015solving}
Stephanus~Daniel Handoko, Lau~Hoong Chuin, Abhishek Gupta, Ong~Yew Soon,
  Heng~Chen Kim, and Tan~Puay Siew.
\newblock Solving multi-vehicle profitable tour problem via knowledge adoption
  in evolutionary bi-level programming.
\newblock In {\em 2015 IEEE Congress on Evolutionary Computation (CEC)}, pages
  2713--2720. IEEE, 2015.

\bibitem{HaJaSa92}
P.~Hansen, B.~Jaumard, and G.~Savard.
\newblock New branch-and-bound rules for linear bilevel programming.
\newblock {\em {SIAM} Journal on Scientific and Statistical Computing},
  13:1194--1217, 1992.

\bibitem{HaPa88}
P.~Harker and J.-S. Pang.
\newblock Existence of optimal solutions to mathematical programs with
  equilibrium constraints.
\newblock {\em Operations Research Letters}, 7:61--64, 1988.

\bibitem{hecheng2008exponential}
Li~Hecheng and Wang Yuping.
\newblock Exponential distribution-based genetic algorithm for solving
  mixed-integer bilevel programming problems.
\newblock {\em Journal of Systems Engineering and Electronics},
  19(6):1157--1164, 2008.

\bibitem{hejazi2002linear}
S.~Reza Hejazi, Azizollah Memariani, G.~Jahanshahloo, and Mohammad~Mehdi
  Sepehri.
\newblock Linear bilevel programming solution by genetic algorithm.
\newblock {\em Computers \& Operations Research}, 29(13):1913--1925, 2002.

\bibitem{herskovits2000contact}
J.~Herskovits, A.~Leontiev, G.~Dias, and G.~Santos.
\newblock Contact shape optimization: a bilevel programming approach.
\newblock {\em Structural and multidisciplinary optimization}, 20(3):214--221,
  2000.

\bibitem{IsAi92}
Y.~Ishizuka and E.~Aiyoshi.
\newblock Double penalty method for bilevel optimization problems.
\newblock {\em Annals of Operations Research}, 34:73--88, 1992.

\bibitem{israeli2002shortest}
Eitan Israeli and R.~Kevin Wood.
\newblock Shortest-path network interdiction.
\newblock {\em Networks}, 40(2):97--111, 2002.

\bibitem{jiang2013application}
Yan Jiang, Xuyong Li, Chongchao Huang, and Xianing Wu.
\newblock Application of particle swarm optimization based on chks smoothing
  function for solving nonlinear bilevel programming problem.
\newblock {\em Applied Mathematics and Computation}, 219(9):4332--4339, 2013.

\bibitem{jin2007bi}
Qin Jin and Shi Feng.
\newblock Bi-level simulated annealing algorithm for facility location.
\newblock {\em Systems Engineering}, 2:007, 2007.

\bibitem{johnson2013inverse}
Mark Johnson, Navid Aghasadeghi, and Timothy Bretl.
\newblock Inverse optimal control for deterministic continuous-time nonlinear
  systems.
\newblock In {\em Decision and Control (CDC), 2013 IEEE 52nd Annual Conference
  on}, pages 2906--2913. IEEE, 2013.

\bibitem{kalashnikov2010comparison}
Vyacheslav Kalashnikov, Fernando Camacho, Ronald Askin, and Nataliya
  Kalashnykova.
\newblock Comparison of algorithms for solving a bi-level toll setting problem.
\newblock {\em International Journal of Innovative Computing, Information and
  Control}, 6(8):3529--3549, 2010.

\bibitem{kalashnikov2015bilevel}
Vyacheslav~V. Kalashnikov, Stephan Dempe, Gerardo~A. P{\'e}rez-Vald{\'e}s,
  Nataliya~I. Kalashnykova, and Jos{\'e}-Fernando Camacho-Vallejo.
\newblock Bilevel programming and applications.
\newblock {\em Mathematical Problems in Engineering}, 2015, 2015.

\bibitem{kalashnikov2016heuristic}
Vyacheslav~V. Kalashnikov, Roberto~Carlos Herrera~Maldonado, Jos{\'e}-Fernando
  Camacho-Vallejo, and Nataliya~I. Kalashnykova.
\newblock A heuristic algorithm solving bilevel toll optimization problems.
\newblock {\em The International Journal of Logistics Management},
  27(1):31--51, 2016.

\bibitem{keeney-raiffa}
R.~L. Keeney and H.~Raiffa.
\newblock {\em Decisions with Multiple Objectives: {P}references and Value
  Tradeoffs}.
\newblock New York: Wiley, 1976.

\bibitem{kovcvara1997topology}
Michal Ko{\v{c}}vara.
\newblock Topology optimization with displacement constraints: a bilevel
  programming approach.
\newblock {\em Structural optimization}, 14(4):256--263, 1997.

\bibitem{kocvara1995solution}
Michal Kocvara and Jifi~V. Outrata.
\newblock On the solution of optimum design problems with variational
  inequalities.
\newblock {\em Recent Advances in Nonsmooth Optimization}, pages 172--192,
  1995.

\bibitem{KoLa90}
C.~Kolstad and L.~Lasdon.
\newblock Derivative evaluation and computational experience with large bilevel
  mathematical programs.
\newblock {\em Journal of Optimization Theory and Applications}, 65:485--499,
  1990.

\bibitem{kuccukaydin2011competitive}
Hande K{\"u}{\c{c}}{\"u}kaydin, Necati Aras, and I~Kuban Alt{\i}nel.
\newblock Competitive facility location problem with attractiveness adjustment
  of the follower: A bilevel programming model and its solution.
\newblock {\em European Journal of Operational Research}, 208(3):206--220,
  2011.

\bibitem{labbe1998bilevel}
Martine Labb{\'e}, Patrice Marcotte, and Gilles Savard.
\newblock A bilevel model of taxation and its application to optimal highway
  pricing.
\newblock {\em Management science}, 44(12-part-1):1608--1622, 1998.

\bibitem{laffont2009theory}
Jean-Jacques Laffont and David Martimort.
\newblock {\em The theory of incentives: the principal-agent model}.
\newblock Princeton university press, 2009.

\bibitem{LeBo86}
L.~Leblanc and D.~Boyce.
\newblock A bilevel programming algorithm for exact solution of the network
  design problem with user-optimal flows.
\newblock {\em Transportation Research}, 20 B:259--265, 1986.

\bibitem{legillon2012cobra}
Fran{\c{c}}ois Legillon, Arnaud Liefooghe, and El-Ghazali Talbi.
\newblock Cobra: A cooperative coevolutionary algorithm for bi-level
  optimization.
\newblock In {\em 2012 IEEE Congress on Evolutionary Computation}, pages 1--8.
  IEEE, 2012.

\bibitem{li2015genetic}
Hecheng Li.
\newblock A genetic algorithm using a finite search space for solving
  nonlinear/linear fractional bilevel programming problems.
\newblock {\em Annals of Operations Research}, pages 1--16, 2015.

\bibitem{li07b}
Hecheng Li and Yuping Wang.
\newblock A hybrid genetic algorithm for solving nonlinear bilevel programming
  problems based on the simplex method.
\newblock {\em International Conference on Natural Computation}, 4:91--95,
  2007.

\bibitem{li06}
Xiangyong Li, Peng Tian, and Xiaoping Min.
\newblock A hierarchical particle swarm optimization for solving bilevel
  programming problems.
\newblock In Leszek Rutkowski, Ryszard Tadeusiewicz, Lotfi~A. Zadeh, and
  Jacek~M. Zurada, editors, {\em Artificial Intelligence and Soft Computing -
  ICAISC 2006}, volume 4029 of {\em Lecture Notes in Computer Science}, pages
  1169--1178. Springer Berlin Heidelberg, 2006.

\bibitem{liang2015evolutionary}
Jason Liang and Risto Miikkulainen.
\newblock Evolutionary bilevel optimization for complex control tasks.
\newblock In {\em Proceedings of the 17th Annual Genetic and Evolutionary
  Computation Conference (GECCO 2015)}. New York: ACM Press, 2015.

\bibitem{LiMo95}
M.~Lignola and J.~Morgan.
\newblock Topological existence and stability for stackelberg problems.
\newblock {\em Journal of Optimization Theory and Applications}, 84:145--169,
  1995.

\bibitem{LiMo02}
Maria~Beatrice Lignola and Jacqueline Morgan.
\newblock Existence of solutions to bilevel variational problems in banach
  spaces.
\newblock In {\em Equilibrium Problems: Nonsmooth Optimization and Variational
  Inequality Models}, pages 161--174. Springer, 2001.

\bibitem{linnala12}
Mikko Linnala, Elina Madetoja, Henri Ruotsalainen, and Jari
  H{\"a}m{\"a}l{\"a}inen.
\newblock Bi-level optimization for a dynamic multiobjective problem.
\newblock {\em Engineering Optimization}, 44(2):195--207, 2012.

\bibitem{liu1998trust}
Guoshan Liu, Jiye Han, and Shouyang Wang.
\newblock A trust region algorithm for bilevel programing problems.
\newblock {\em Chinese science bulletin}, 43(10):820--824, 1998.

\bibitem{LoMo96}
P.~Loridan and J.~Morgan.
\newblock {Weak via strong Stackelberg problems}.
\newblock {\em Journal of Global Optimization}, 8:263--287, 1996.

\bibitem{lowe06}
J.~Lowe.
\newblock {Homeland Security: Operations Research Initiatives and
  Applications}.
\newblock {\em Interfaces}, 36(6):483--485, 2006.

\bibitem{LuMiPi87}
R.~Lucchetti, F.~Mignanego, and G.~Pieri.
\newblock {Existence theorem of equilibrium points in Stackelberg games with
  constraints}.
\newblock {\em Optimization}, 18:857--866, 1987.

\bibitem{lv2007penalty}
Yibing Lv, Tiesong Hu, Guangmin Wang, and Zhongping Wan.
\newblock A penalty function method based on kuhn--tucker condition for solving
  linear bilevel programming.
\newblock {\em Applied Mathematics and Computation}, 188(1):808--813, 2007.

\bibitem{maldonado2016analyzing}
Sayuri Maldonado-Pinto, Martha-Selene Casas-Ram{\'\i}rez, and Jos{\'e}-Fernando
  Camacho-Vallejo.
\newblock Analyzing the performance of a hybrid heuristic for solving a bilevel
  location problem under different approaches to tackle the lower level.
\newblock {\em Mathematical Problems in Engineering}, 2016, 2016.

\bibitem{MaMa92}
P.~Marcotte and G.~Marquis.
\newblock Efficient implementation of heuristics for the continuous network
  design problem.
\newblock {\em Annals of Operations Research}, 34:163--176, 1992.

\bibitem{marcotte2004bilevel}
Patrice Marcotte, Gilles Savard, and Fr{\'e}d{\'e}ric Semet.
\newblock A bilevel programming approach to the travelling salesman problem.
\newblock {\em Operations Research Letters}, 32(3):240--248, 2004.

\bibitem{marcotte2001trust}
Patrice Marcotte, Gilles Savard, and D.~L. Zhu.
\newblock A trust region algorithm for nonlinear bilevel programming.
\newblock {\em Operations research letters}, 29(4):171--179, 2001.

\bibitem{maric2014metaheuristic}
Miroslav Mari{\'c}, Zorica Stanimirovi{\'c}, Nikola Milenkovi{\'c}, and
  Aleksandar {\DJ}eni{\'c}.
\newblock Metaheuristic approaches to solving large-scale bilevel uncapacitated
  facility location problem with clients'preferences.
\newblock {\em Yugoslav Journal of Operations Research ISSN: 0354-0243 EISSN:
  2334-6043}, 25(3), 2014.

\bibitem{mathieu}
R.~Mathieu, L.~Pittard, and G.~Anandalingam.
\newblock Genetic algorithm based approach to bi-level linear programming.
\newblock {\em Operations Research}, 28(1):1--21, 1994.

\bibitem{mesbah2011optimization}
Mahmoud Mesbah, Majid Sarvi, and Graham Currie.
\newblock Optimization of transit priority in the transportation network using
  a genetic algorithm.
\newblock {\em Intelligent Transportation Systems, IEEE Transactions on},
  12(3):908--919, 2011.

\bibitem{miandoabchi2011optimizing}
Elnaz Miandoabchi and Reza~Zanjirani Farahani.
\newblock Optimizing reserve capacity of urban road networks in a discrete
  network design problem.
\newblock {\em Advances in Engineering Software}, 42(12):1041--1050, 2011.

\bibitem{migdalas95}
Athanasios Migdalas.
\newblock Bilevel programming in traffic planning: Models, methods and
  challenge.
\newblock {\em Journal of Global Optimization}, 7(4):381--405, 1995.

\bibitem{mombaur2010human}
Katja Mombaur, Anh Truong, and Jean-Paul Laumond.
\newblock From human to humanoid locomotion\d1an inverse optimal control
  approach.
\newblock {\em Autonomous robots}, 28(3):369--383, 2010.

\bibitem{moore1990mixed}
James~T. Moore and Jonathan~F. Bard.
\newblock The mixed integer linear bilevel programming problem.
\newblock {\em Operations research}, 38(5):911--921, 1990.

\bibitem{neumann-morgenstern}
John~Von Neumann and Oskar Morgenstern.
\newblock {\em Theory of games and economic behavior}.
\newblock Princeton, NJ. Princeton University Press., 1946.

\bibitem{o2011designing}
Jesse~R. O'Hanley and Richard~L. Church.
\newblock Designing robust coverage networks to hedge against worst-case
  facility losses.
\newblock {\em European Journal of Operational Research}, 209(1):23--36, 2011.

\bibitem{Ou93}
J.~Outrata.
\newblock {Necessary optimality conditions for Stackelberg problems}.
\newblock {\em Journal of Optimization Theory and Applications}, 76:305--320,
  1993.

\bibitem{panin2014bilevel}
Artem~Aleksandrovich Panin, M.~G. Pashchenko, and Aleksandr~Vladimirovich
  Plyasunov.
\newblock Bilevel competitive facility location and pricing problems.
\newblock {\em Automation and Remote Control}, 75(4):715--727, 2014.

\bibitem{pieume09}
C.~O. Pieume, L.~P. Fotso, and P.~Siarry.
\newblock Solving bilevel programming problems with multicriteria optimization
  techniques.
\newblock {\em OPSEARCH}, 46(2):169--183, 2009.

\bibitem{pramnik11}
Surapati Pramanik and Partha~Pratim Dey.
\newblock Bi-level multi-objective programming problem with fuzzy parameters.
\newblock {\em International Journal of Computer Applications}, 30(10):13--20,
  September 2011.
\newblock Published by Foundation of Computer Science, New York, USA.

\bibitem{raghunathan2003mathematical}
Arvind~U. Raghunathan and Lorenz~T. Biegler.
\newblock Mathematical programs with equilibrium constraints (mpecs) in process
  engineering.
\newblock {\em Computers \& Chemical Engineering}, 27(10):1381--1392, 2003.

\bibitem{ramamoorthy2016hub}
P.~Ramamoorthy, S.~Jayaswal, A.~Sinha, and N.~Vidyarthi.
\newblock Hub interdiction \& hub protection problems: Model formulations \&
  exact solution methods. {(No. WP 2016-10-01)}.
\newblock Technical report, Indian Institute of Management Ahmedabad, 2016.

\bibitem{ramamoorthy2017hub}
P.~Ramamoorthy, S.~Jayaswal, A.~Sinha, and N.~Vidyarthi.
\newblock Hub-and-spoke network design under the risk of interdiction. {(No. WP
  2017-05-01)}.
\newblock Technical report, Indian Institute of Management Ahmedabad, 2017.

\bibitem{ren2013integrated}
Gang Ren, Zhengfeng Huang, Yang Cheng, Xing Zhao, and Yong Zhang.
\newblock An integrated model for evacuation routing and traffic signal
  optimization with background demand uncertainty.
\newblock {\em Journal of Advanced Transportation}, 47(1):4--27, 2013.

\bibitem{ruuska12}
S.~Ruuska and K.~Miettinen.
\newblock Constructing evolutionary algorithms for bilevel multiobjective
  optimization.
\newblock In {\em Evolutionary Computation (CEC), 2012 IEEE Congress on}, pages
  1--7, june 2012.

\bibitem{saharidis2009resolution}
Georges~K. Saharidis and Marianthi~G. Ierapetritou.
\newblock Resolution method for mixed integer bi-level linear problems based on
  decomposition technique.
\newblock {\em Journal of Global Optimization}, 44(1):29--51, 2009.

\bibitem{SaGa94}
G.~Savard and J.~Gauvin.
\newblock The steepest descent direction for the nonlinear bilevel programming
  problem.
\newblock {\em Operations Research Letters}, 15:275--282, 1994.

\bibitem{scaparra2008bilevel}
Maria~P. Scaparra and Richard~L. Church.
\newblock A bilevel mixed-integer program for critical infrastructure
  protection planning.
\newblock {\em Computers \& Operations Research}, 35(6):1905--1923, 2008.

\bibitem{shi2005extended}
Chenggen Shi, Jie Lu, and Guangquan Zhang.
\newblock An extended kuhn--tucker approach for linear bilevel programming.
\newblock {\em Applied Mathematics and Computation}, 162(1):51--63, 2005.

\bibitem{shi-xia}
X.~Shi. and H.~S. Xia.
\newblock Model and interactive algorithm of bi-level multi-objective
  decision-making with multiple interconnected decision makers.
\newblock {\em Journal of Multi-Criteria Decision Analysis}, 10(1):27--34,
  2001.

\bibitem{ShIsBa97}
K.~Shimizu, Y.~Ishizuka, and J.~F. Bard.
\newblock {\em Nondifferentiable and two-level mathematical programming}.
\newblock Dordrecht: Kluwer Academic, 1997.

\bibitem{my-emo11}
A.~Sinha.
\newblock Bilevel multi-objective optimization problem solving using
  progressively interactive evolutionary algorithm.
\newblock In {\em Proceedings of the Sixth International Conference on
  Evolutionary Multi-Criterion Optimization (EMO-2011)}, pages 269--284.
  Berlin, Germany: Springer-Verlag, 2011.

\bibitem{my-ifac09}
A.~Sinha and K.~Deb.
\newblock Towards understanding evolutionary bilevel multi-objective
  optimization algorithm.
\newblock In {\em IFAC Workshop on Control Applications of Optimization
  (IFAC-2009)}, volume~7. Elsevier, 2009.

\bibitem{my-bleaq2-arxiv17}
A.~Sinha, Z.~Lu, K.~Deb, and P.~Malo.
\newblock Bilevel optimization based on iterative approximation of mappings.
\newblock {\em arXiv preprint arXiv:1702.03394}, 2017.

\bibitem{my-bleaq-arxiv13}
A.~Sinha, P.~Malo, and K.~Deb.
\newblock Efficient evolutionary algorithm for single-objective bilevel
  optimization.
\newblock {\em arXiv preprint arXiv:1303.3901}, 2013.

\bibitem{my-cec14}
A.~Sinha, P.~Malo, and K.~Deb.
\newblock An improved bilevel evolutionary algorithm based on quadratic
  approximations.
\newblock In {\em 2014 IEEE Congress on Evolutionary Computation (CEC-2014)},
  pages 1870--1877. IEEE Press, 2014.

\bibitem{my-emo15}
A.~Sinha, P.~Malo, and K.~Deb.
\newblock Towards understanding bilevel multi-objective optimization with
  deterministic lower level decisions.
\newblock In {\em Proceedings of the Eighth International Conference on
  Evolutionary Multi-Criterion Optimization (EMO-2015)}. Berlin, Germany:
  Springer-Verlag, 2015.

\bibitem{my-cec15a}
A.~Sinha, P.~Malo, and K.~Deb.
\newblock Transportation policy formulation as a multi-objective bilevel
  optimization problem.
\newblock In {\em 2015 IEEE Congress on Evolutionary Computation (CEC-2015)}.
  IEEE Press, 2015.

\bibitem{my-cec16a}
A.~Sinha, P.~Malo, and K.~Deb.
\newblock Solving optimistic bilevel programs by iteratively approximating
  lower level optimal value function.
\newblock In {\em 2016 IEEE Congress on Evolutionary Computation (CEC-2016)}.
  IEEE Press, 2016.

\bibitem{my-ejor16b}
A.~Sinha, P.~Malo, and K.~Deb.
\newblock Evolutionary algorithm for bilevel optimization using approximations
  of the lower level optimal solution mapping.
\newblock {\em European Journal of Operational Research}, 2016 (In press).

\bibitem{sinha2017evolutionary}
A.~Sinha, P.~Malo, and K.~Deb.
\newblock Evolutionary bilevel optimization: {A}n introduction and recent
  advances.
\newblock In {\em Recent Advances in Evolutionary Multi-objective
  Optimization}, pages 71--103. Springer, 2017.

\bibitem{my-ieeetec16}
A.~Sinha, P.~Malo, K.~Deb, P.~Korhonen, and J.~Wallenius.
\newblock Solving bilevel multi-criterion optimization problems with lower
  level decision uncertainty.
\newblock {\em IEEE Transactions on Evolutionary Computation}, 20(2):199--217,
  2016.

\bibitem{my-cec13}
A.~Sinha, P.~Malo, A.~Frantsev, and K.~Deb.
\newblock Multi-objective stackelberg game between a regulating authority and a
  mining company: A case study in environmental economics.
\newblock In {\em 2013 IEEE Congress on Evolutionary Computation (CEC-2013)}.
  IEEE Press, 2013.

\bibitem{my-caor14}
A.~Sinha, P.~Malo, A.~Frantsev, and K.~Deb.
\newblock Finding optimal strategies in a multi-period multi-leader-follower
  stackelberg game using an evolutionary algorithm.
\newblock {\em Computers \& Operations Research}, 41:374--385, 2014.

\bibitem{my-gecco14}
A.~Sinha, P.~Malo, P.~Xu, and K.~Deb.
\newblock A bilevel optimization approach to automated parameter tuning.
\newblock In {\em Proceedings of the 16th Annual Genetic and Evolutionary
  Computation Conference (GECCO 2014)}. New York: ACM Press, 2014.

\bibitem{smith82}
W.~R. Smith and R.~W. Missen.
\newblock {\em Chemical Reaction Equilibrium Analysis: Theory and Algorithms}.
\newblock John Wiley \& Sons, New York, 1982.

\bibitem{St52}
H.~Stackelberg.
\newblock {\em The theory of the market economy}.
\newblock Oxford University Press, New York, Oxford, 1952.

\bibitem{sun2008bi}
Huijun Sun, Ziyou Gao, and Jianjun Wu.
\newblock A bi-level programming model and solution algorithm for the location
  of logistics distribution centers.
\newblock {\em Applied Mathematical Modelling}, 32(4):610--616, 2008.

\bibitem{sun08}
Huijun Sun, Ziyou Gao, and Jianjun Wu.
\newblock A bi-level programming model and solution algorithm for the location
  of logistics distribution centers.
\newblock {\em Applied Mathematical Modelling}, 32(4):610 -- 616, 2008.

\bibitem{my-cec16b}
V.~Suryan, A.~Sinha, , P.~Malo, and K.~Deb.
\newblock Handling inverse optimal control problems using evolutionary bilevel
  optimization.
\newblock In {\em 2016 IEEE Congress on Evolutionary Computation (CEC-2016)}.
  IEEE Press, 2016.

\bibitem{talbi2013metaheuristics}
El-Ghazali Talbi.
\newblock {\em Metaheuristics for bi-level optimization}, volume 482.
\newblock Springer, 2013.

\bibitem{teh2006hierarchical}
Yee~Whye Teh, Michael~I Jordan, Matthew~J Beal, and David~M Blei.
\newblock Hierarchical dirichlet processes.
\newblock {\em Journal of the american statistical association}, 101(476),
  2006.

\bibitem{TuMiVa93}
H.~Tuy, A.~Migdalas, and P.~V\"{a}rbrand.
\newblock A global optimization approach for the linear two-level program.
\newblock {\em Journal of Global Optimization}, 3:1--23, 1993.

\bibitem{uno2008evolutionary}
Takeshi Uno, Hideki Katagiri, and Kosuke Kato.
\newblock An evolutionary multi-agent based search method for stackelberg
  solutions of bilevel facility location problems.
\newblock {\em International Journal of Innovative Computing, Information and
  Control}, 4(5):1033--1042, 2008.

\bibitem{van1993principal}
A.~Van~Ackere.
\newblock The principal/agent paradigm: characterizations and computations.
\newblock {\em European Journal of Operations Research}, 70:83--103, 1993.

\bibitem{ViCa94}
L.~Vicente and P.~Calamai.
\newblock Bilevel and multilevel programming: a bibliography review.
\newblock {\em Journal of Global Optimization}, 5:291--306, 1994.

\bibitem{ViSaJu94}
L.~Vicente, G.~Savard, and J.~J\'udice.
\newblock Descent approaches for quadratic bilevel programming.
\newblock {\em Journal of Optimization Theory and Applications}, 81:379--399,
  1994.

\bibitem{ViSaJu96}
L.~Vicente, G.~Savard, and J.~J\'udice.
\newblock The discrete linear bilevel programming problem.
\newblock {\em Journal of Optimization Theory and Applications}, 89:597--614,
  1996.

\bibitem{vicente1996discrete}
Luis Vicente, Gilles Savard, and J.~Judice.
\newblock Discrete linear bilevel programming problem.
\newblock {\em Journal of optimization theory and applications},
  89(3):597--614, 1996.

\bibitem{wan2013hybrid}
Zhongping Wan, Guangmin Wang, and Bin Sun.
\newblock A hybrid intelligent algorithm by combining particle swarm
  optimization with chaos searching technique for solving nonlinear bilevel
  programming problems.
\newblock {\em Swarm and Evolutionary Computation}, 8:26--32, 2013.

\bibitem{wang2007review}
G.~Gary Wang and S.~Shan.
\newblock Review of metamodeling techniques in support of engineering design
  optimization.
\newblock {\em Journal of Mechanical Design}, 129(4):370--380, 2007.

\bibitem{wang2008genetic}
Guangmin Wang, Zhongping Wan, Xianjia Wang, and Yibing Lv.
\newblock Genetic algorithm based on simplex method for solving
  linear-quadratic bilevel programming problem.
\newblock {\em Computers \& Mathematics with Applications}, 56(10):2550--2555,
  2008.

\bibitem{wang2014bilevel}
Judith Y.~T. Wang, Matthias Ehrgott, Kim~N. Dirks, and Abhishek Gupta.
\newblock A bilevel multi-objective road pricing model for economic,
  environmental and health sustainability.
\newblock {\em Transportation Research Procedia}, 3:393--402, 2014.

\bibitem{wang05}
Y.~Wang, Y.~C. Jiao, and H.~Li.
\newblock An evolutionary algorithm for solving nonlinear bilevel programming
  based on a new constraint-handling scheme.
\newblock {\em IEEE Transactions on Systems, Man, and Cybernetics, Part {C}:
  {A}pplications and Reviews}, 32(2):221--232, 2005.

\bibitem{wang11}
Yuping Wang, Hong Li, and Chuangyin Dang.
\newblock A new evolutionary algorithm for a class of nonlinear bilevel
  programming problems and its global convergence.
\newblock {\em INFORMS Journal on Computing}, 23(4):618--629, 2011.

\bibitem{wein09}
L.~Wein.
\newblock {Homeland Security: From Mathematical Models to Policy
  Implementation: The 2008 Philip McCord Morse Lecture}.
\newblock {\em Operations Research}, 57(4):801--811, 2009.

\bibitem{WeHs91}
U.~Wen and S.~Hsu.
\newblock Linear bi-level programming problems - a review.
\newblock {\em Journal of the Operational Research Society}, 42:125--133, 1991.

\bibitem{WhAn93}
D.~White and G.~Anandalingam.
\newblock A penalty function approach for solving bi-level linear programs.
\newblock {\em Journal of Global Optimization}, 3:397--419, 1993.

\bibitem{whittaker2016spatial}
Gerald Whittaker, Rolf F{\"a}re, Shawna Grosskopf, Bradley Barnhart, Moriah
  Bostian, George Mueller-Warrant, and Stephen Griffith.
\newblock Spatial targeting of agri-environmental policy using bilevel
  evolutionary optimization.
\newblock {\em Omega}, 2016.

\bibitem{wiesemann2013pessimistic}
Wolfram Wiesemann, Angelos Tsoukalas, Polyxeni-Margarita Kleniati, and
  Ber{\c{c}} Rustem.
\newblock Pessimistic bilevel optimization.
\newblock {\em SIAM Journal on Optimization}, 23(1):353--380, 2013.

\bibitem{xu2012nonlinear}
Jiuping Xu, Yan Tu, and Ziqiang Zeng.
\newblock A nonlinear multiobjective bilevel model for minimum cost network
  flow problem in a large-scale construction project.
\newblock {\em Mathematical Problems in Engineering}, 2012, 2012.

\bibitem{xu2007supply}
Qing Xu, Dao-li Zhu, and Shan-liang Li.
\newblock The supply chain optimal contract design under asymmetrical
  information [j].
\newblock {\em Systems Engineering-Theory \& Practice}, 4:003, 2007.

\bibitem{yamada2009designing}
Tadashi Yamada, Bona~Frazila Russ, Jun Castro, and Eiichi Taniguchi.
\newblock Designing multimodal freight transport networks: a heuristic approach
  and applications.
\newblock {\em Transportation Science}, 43(2):129--143, 2009.

\bibitem{yang1998models}
Hai Yang and Michael~G. H.~Bell.
\newblock Models and algorithms for road network design: a review and some new
  developments.
\newblock {\em Transport Reviews}, 18(3):257--278, 1998.

\bibitem{ye2011necessary}
Jane~J Ye.
\newblock Necessary optimality conditions for multiobjective bilevel programs.
\newblock {\em Mathematics of Operations Research}, 36(1):165--184, 2011.

\bibitem{ye2010new}
Jane~J. Ye and Daoli Zhu.
\newblock New necessary optimality conditions for bilevel programs by combining
  the mpec and value function approaches.
\newblock {\em SIAM Journal on Optimization}, 20(4):1885--1905, 2010.

\bibitem{yin-bilevel}
Y.~Yin.
\newblock Genetic algorithm based approach for bilevel programming models.
\newblock {\em Journal of Transportation Engineering}, 126(2):115--120, 2000.

\bibitem{yin2000genetic}
Yafeng Yin.
\newblock Genetic-algorithms-based approach for bilevel programming models.
\newblock {\em Journal of Transportation Engineering}, 126(2):115--120, 2000.

\bibitem{yin2002multiobjective}
Yafeng Yin.
\newblock Multiobjective bilevel optimization for transportation planning and
  management problems.
\newblock {\em Journal of advanced transportation}, 36(1):93--105, 2002.

\bibitem{Bilevel-linear}
G.~Zhang, J.~Liu, and T.~Dillon.
\newblock Decentralized multi-objective bilevel decision making with fuzzy
  demands.
\newblock {\em Knowledge-Based Systems}, 20:495--507, 2007.

\bibitem{zhang12}
T.~Zhang, T.~Hu, Y.~Zheng, and X.~Guo.
\newblock An improved particle swarm optimization for solving bilevel
  multiobjective programming problem.
\newblock {\em Journal of Applied Mathematics}, 2012.

\bibitem{zhu2006hybrid}
Xiaobo Zhu, Qian Yu, and Xianjia Wang.
\newblock A hybrid differential evolution algorithm for solving nonlinear
  bilevel programming with linear constraints.
\newblock In {\em Cognitive Informatics, 2006. ICCI 2006. 5th IEEE
  International Conference on}, volume~1, pages 126--131. IEEE, 2006.

\end{thebibliography}

\end{document}